\documentclass[11pt,onecolumn]{article}
\usepackage[left=1in,top=1in,right=1in,bottom=1in]{geometry}
\usepackage{graphicx,caption,subfigure,epstopdf,grffile,hyperref,authblk,topcapt,enumerate}
\usepackage{amssymb,amsmath,amsthm,xcolor,mathrsfs,diagbox}
\numberwithin{equation}{section}
\graphicspath{{./Figures/}}	
\DeclareMathOperator{\eps}{\varepsilon}

\newcommand{\tl}[1]{\tilde{#1}}
        \newcommand{\mc}[1]{\mathcal{#1}}
        \newcommand{\mb}[1]{\mathbb{#1}}
                \newcommand{\p}{\partial}

\title{\vspace{-20pt}Delayed Hopf bifurcation
and space-time buffer curves
in the Complex Ginzburg-Landau equation}
\author[1]{Ryan Goh}
\author[1]{Tasso J. Kaper}
\author[2]{Theodore Vo}
\affil[1]{\footnotesize Department of Mathematics and Statistics, Boston University, Boston, MA 02215, USA}
\affil[2]{\footnotesize School of Mathematics, Monash University, Clayton, Victoria 3800, Australia}

\begin{document}
\bibliographystyle{abbrv}
\maketitle
%=========================================================================================

%------------------------------------------------------------
\vspace{-10pt}
\begin{abstract}	\label{sec:abstract}
\noindent In this article, the recently-discovered
phenomenon of delayed Hopf bifurcations 
(DHB) in reaction-diffusion PDEs
is analyzed in the cubic Complex Ginzburg-Landau equation,
as an equation in its own right,
with a slowly-varying parameter.
We begin by using the classical asymptotic methods of stationary phase
and steepest descents
to show that solutions which have approached 
the attracting quasi-steady state (QSS)
before the Hopf bifurcation remain near that state for long
times after the instantaneous Hopf bifurcation 
and the QSS has become repelling.
In the complex time plane,
the phase function
of the linear PDE has a saddle point,
and the Stokes and anti-Stokes lines
are central to the asymptotics.
The nonlinear terms 
are treated 
by applying an iterative method
to the mild form of the PDE
given by perturbations 
about the linear particular solution.
This tracks the closeness of solutions
near the attracting and repelling QSS.
Next, we show that
beyond a key Stokes line through the saddle
there is a curve in the space-time plane
along which the particular solution
of the linear PDE 
ceases to be exponentially small,
causing the solution of the nonlinear PDE
to diverge from the repelling QSS
and exhibit large-amplitude oscillations. 
This curve is called the space-time buffer curve.
The homogeneous solution 
also stops being exponentially small
in a spatially dependent manner,
as determined also by the initial time.
Hence, a competition arises
between these two solutions,
as to which one ceases to be exponentially small first,
and this competition
governs spatial dependence of the DHB.
We find four different cases of DHB,
depending on the outcomes of the competition,
and we quantify to leading order how these depend
on the main system parameters,
including the Hopf frequency,
initial time, initial data,
source terms,
and diffusivity.
Examples are presented 
for each case,
with source terms that are
uni-modal, 
spatially-periodic, 
smooth step function, 
and algebraically-growing.
Also, rich 
spatio-temporal dynamics
are observed in the post-DHB oscillations.
Finally, it is shown that large-amplitude source terms
can be designed 
so that solutions spend
substantially longer times
near the repelling QSS,
and hence region-specific control
over the delayed onset of oscillations 
can be achieved.

\vspace{8pt}\noindent \textbf{Keywords}\qquad 
Slow passage through Hopf bifurcation,
dynamic bifurcation in PDEs,
spatially-inhomogeneous onset of oscillations,
hard onset of oscillations,
stationary phase method,
steepest descents method,
complex Ginzburg-Landau equation,
reaction-diffusion equations

\end{abstract}
%------------------------------------------------------------

%---------------------------------------------------------------------------------
\section{Introduction}		\label{sec:intro}
%---------------------------------------------------------------------------------

In applied mathematics, physics, biology,
and in many other areas of science and engineering,
the phenomenon of delayed Hopf bifurcation (DHB)
is a central feature 
of analytic ordinary differential equations (ODEs)
in which a parameter passes slowly through a Hopf point.
ODE models and experimental
examples arise in chemistry and pattern formation
\cite{BDN1997,ER1990,ERHG1991,KA1996,K1998,SM1996},
nonlinear mechanical oscillators and generalised Rayleigh oscillators
\cite{B1993,PDL2011,PSBT2016},
electrical engineering \cite{HXJBK2016,WYCXB2019},
fluid dynamics and geophysics \cite{AWVC2012,EKKV2017,HPK2008,MS1990},
neuroscience \cite{BC2011,BBKS1995,BB2018,BKG2010,DHCG1998, 
GDSSAH2007,GJK2001,I2001,RB1988,RT2002,
SW2000,S1993,UWDWS2015}, 
cardiac models \cite{K2016}, and 
the Kaldor model in business \cite{GW1994}.

In DHB for ODEs, the key system parameter 
passes slowly in time
through a Hopf bifurcation value
at which the stable equilibrium becomes unstable,
yet the solutions remain near the repelling equilibrium
for long times,
of length $\mathcal{O}(1)$ in the slow time,
after the Hopf point.
As a result, 
in the super-critical case,
the attendant (post-DHB) onset of oscillations
is a hard onset, 
with solutions
jumping rapidly away from the unstable equilibrium
to the stable limit cycle, which by that time
already has a large amplitude.
DHB has been studied in analytic ODEs
for more than 50 years,
going back at least to the seminal work \cite{S1973}.
The theory is further developed in
\cite{BER1989,ER1990,HKSW2016,HE1993,N1987,N1988,S1993},
and is also presented in recent monographs
\cite{K2015,W2020}.
Moreover, many of the above applications
have been inspired by \cite{BER1989,N1987,N1988}.

DHB has also been studied \cite{BB2018}
in large systems of ODEs.
There, the FitzHugh-Nagumo and Hodgkin-Huxley cable equations 
are studied 
with a slowly varying Neumann boundary condition 
at one end of the spatial domain 
and a zero-flux condition at the other end.
The spatial variable in these partial differential equations
(PDEs) is discretised 
(with centered finite differences for the Laplacian),
and the WKB method is used 
to analyze DHB
in the large system of ODEs.

Recently, it was discovered 
that the phenomenon of DHB
also occurs in reaction-diffusion equations
\cite{KV2018}.
In that work, it was shown using numerics, 
physical considerations, and some Fourier analysis
that DHB is important for a variety
of reaction-diffusion equations
in which there is slow passage 
through super-critical Hopf bifurcations.
The reaction-diffusion examples 
in which DHB has been found 
\cite{KV2018} include
the Complex Ginzburg-Landau equation,
the Brusselator model of the Belousov-Zhabotinsky reaction,
the Hodgkin-Huxley PDE,
the FitzHugh-Nagumo PDE, and
a spatially-extended pituitary lactotroph cell model.

There has been rigorous analysis
in \cite{ADVW2020}
of spatio-temporal canards and delayed bifurcations
in a class of infinite-dimensional systems on bounded domains,
which includes slow passage through Hopf bifurcation,
slow passage through a Turing bifurcation,
and some bifurcations in delay-differential equations.
In that article, 
under the assumption
that a spectral gap exists in the fast sub-system,
a center manifold analysis is performed 
using the infinite-dimensional invariant manifold theory
of Haragus and Iooss \cite{HI2010}.
See also the references in \cite{ADVW2020,VBK2020}
for more on spatio-temporal canards.

In this article,
we use the methods of stationary phase 
and steepest descents to analyze the DHB
created when the bifurcation parameter $\mu$
slowly increases through 
a supercritical Hopf bifurcation (at zero)
in the Complex Ginzburg Landau equation 
on the real line,
\begin{equation} \label{eq:gen-CGL}
\begin{split}
A_t &= (\mu + i \omega_0) A - (1+i\alpha)\vert A\vert^2 A 
+ \eps^\beta I_a(x) + \eps^\gamma d A_{xx}
\\
\mu_t &=  \eps.
\end{split}
\end{equation}
Here, $x$ is real,
$t \ge 0$,
$A = A(x,t)$
is complex-valued,
and $0< \eps \ll 1$
is a small parameter.
The linear growth rate
$\mu = \mu(t)$ is real 
for the main phenomena we study;
however, for the mathematical analysis,
it will be advantageous to consider
complex values of $\mu$ 
in a horizontal strip 
with mid-line on the real axis
and of sufficient height.
The system parameters satisfy
$\omega_0>0$ and $\mathcal{O}(1)$ independent of $\eps$,
$\alpha$ is real,
$\beta>0$,
$\gamma \ge 0$,
$d$ may be complex-valued 
($d = d_R+i d_I$) with $d_R>0$,
and they
are independent of $\eps$.
For real values of $x$,
the source term $I_a(x)$,
which breaks the symmetry $ A \to A e^{i\theta}$
for any real $\theta$ of the CGL equation,
is typically taken to be bounded
and positive, 
with uniformly bounded derivatives.
The initial data 
at $\mu(0)=\mu_0 < 0$ is 
$A(x,0)=A_0(x)$,
and typically taken to be bounded and continuous 
for all real $x$.
Also, it will be useful 
to distinguish between initial data
given at $\mu_0\le -\omega_0$
and data 
given at $-\omega_0 < \mu_0 < 0$.

The PDE \eqref{eq:gen-CGL}
has an attracting Quasi-Steady State (QSS)
for all $\mu < -\delta$,
where $\delta>0$, small, and $\mathcal{O}(1)$,
which solutions approach at an exponential rate.
Similarly,
it has a repelling QSS
for all $\mu > \delta$,
from which solutions diverge at an exponential rate.
For example, 
in the base case 
of $\beta=\frac{1}{2}$
and $\gamma=1$, 
the attracting QSS (for $\mu<-\delta$)
and the repelling QSS (for $\mu>\delta$)
are given by
\begin{equation}\label{eq:QSS-basecase}
A_{\rm QSS} (x,\mu)
= -\sqrt{\eps} \frac{I_a(x)}{\mu+i\omega_0} 
+ \eps^{\frac{3}{2}} \left( \frac{I_a(x)+d(\mu+i\omega_0){I_a}''(x)}
              {(\mu+i \omega_0)^3}
        -\frac{(1+i\alpha)I_a^3(x)}{(\mu+i\omega_0)^2(\mu^2+\omega_0^2)}
        \right) 
+ \mathcal{O}(\eps^{\frac{5}{2}}).
\end{equation}
Here, the $\mathcal{O}(\eps^\frac{5}{2})$
terms depend on $x$ and $\mu$.
The QSS may also be derived
for other $\beta$ and $\gamma$.

\begin{figure}[htbp]
   \centering
   \includegraphics[width=5in]{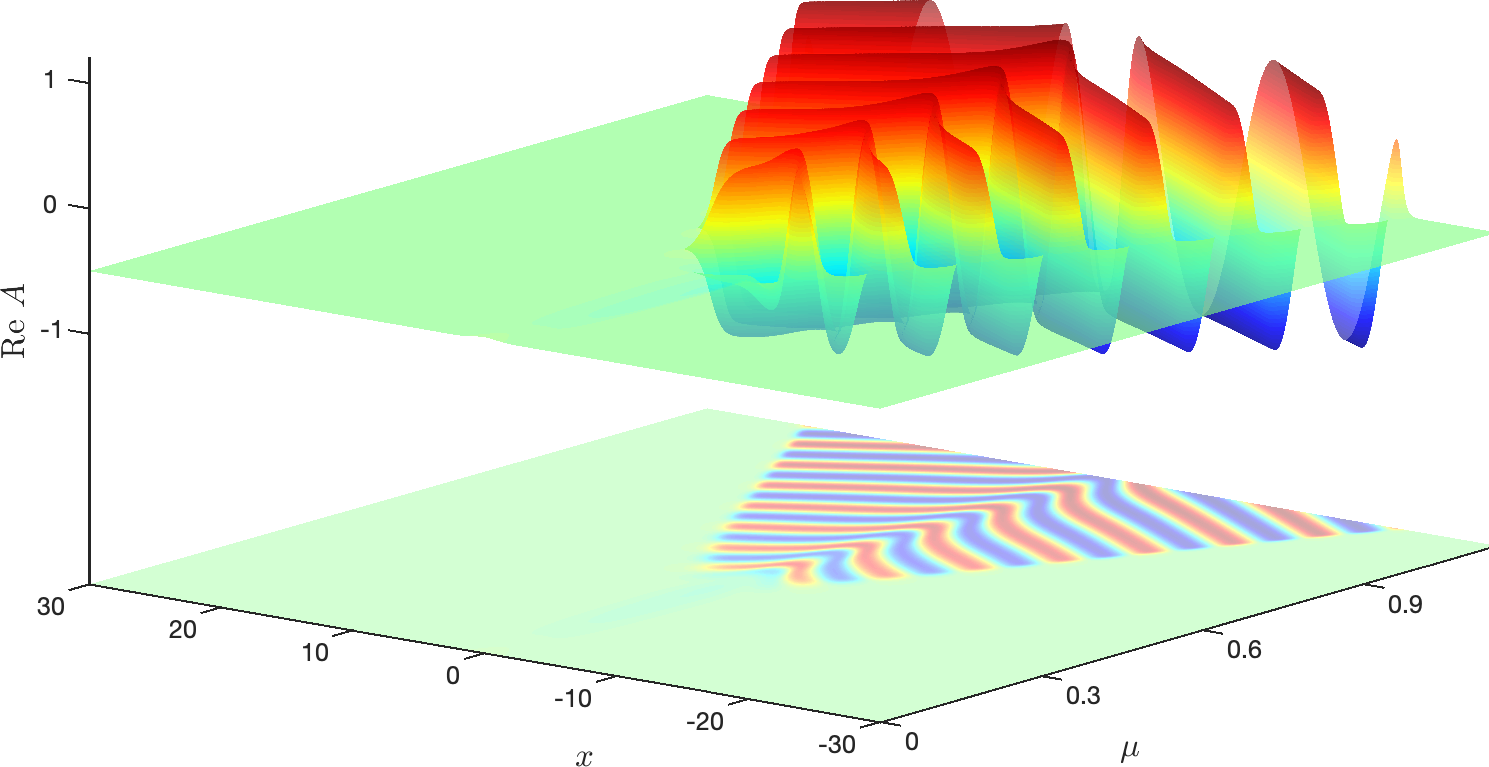}
\caption{
DHB in 
the PDE \eqref{eq:gen-CGL}
in the base case of $\beta=\frac{1}{2}$
and $\gamma=1$,
with a Gaussian source term $I_a(x)= e^{-\frac{x^2}{4}}$
on the domain $[-\ell,\ell]$ ($\ell=30$).
The solution 
is near the attracting QSS while $\mu \in [\mu_0,0)$
(with $\mu_0=-1$, not shown).
Then, for $\mu>0$,
{\it i.e.,} after the instantaneous Hopf bifurcation,
the solution stays near the repelling QSS (green state)
at least until $\mu=\omega_0=0.5$
at all points $x$.
The duration of the DHB 
({\it i.e.,} the length of time the solution
stays near the repelling QSS)
is spatially dependent.
At $x=0$,
$A(x,\mu)$ first leaves a neighbourhood 
of the repelling QSS at time $\mu=\omega_0$,
and the large-amplitude oscillations 
first set in there.
Then, as one steps outward from the center,
with $\vert x \vert > 0$,
the duration of the DHB 
grows beyond $\mu=\omega_0$.
The hard onset of the large-amplitude oscillations
is governed to leading order
by the space-time buffer curve studied in this article.
Here, $\eps = 0.01, \omega_0 = \tfrac{1}{2}, d_R = 1, d_I = 0,$
and $\alpha = 0$.
	}
	   \label{fig:DHB_3D}
	\end{figure}

In \cite{KV2018},
DHB is observed
for solutions which
are in a fixed 
$\mathcal{O}(1)$ neighbourhood
of the attracting QSS
for any $\mu_0$
sufficiently negative.
These solutions all
continue to approach the attracting QSS
until $\mu=0$, 
where the instantaneous Hopf bifurcation occurs.
However, rather than immediately tracking 
the stable (post-DHB) oscillatory state 
as it grows in amplitude,
these solutions 
remain near the repelling QSS for long times into $\mu>0$.
Moreover, the amount of time any such solution $A(x,t)$
spends near the repelling QSS 
can, and generally does, depend on $x$.
See Figure~\ref{fig:DHB_3D}.

Our first goal in this article
is to derive a general formula 
for the space-time buffer curve
of system \eqref{eq:gen-CGL}
with general source terms $I_a(x)$.
The space-time buffer curve
corresponds to the 
$x$-dependent (post-DHB) time
at which the solution cannot 
remain near the repelling QSS any longer 
for each point $x$,
irrespective of how far
before the slowly-varying Hopf point
the solution was attracted 
to the (pre-Hopf) stable QSS.
We directly to leading order
apply the classical methods
of stationary phase and steepest descents
(see for example \cite{BO1999,KC1981,M1984})
to the linear CGL equation,
obtained by linearising \eqref{eq:gen-CGL} 
about $A=0$.
The coefficient
$(\mu + i \omega_0)$ 
on the linear homogeneous term vanishes 
at $\mu = - i \omega_0$
in the complex $\mu$ plane.
This is a saddle point of the complex phase 
$-\frac{1}{2}(\mu + i \omega_0)^2$,
since the derivative of the phase vanishes there.
Moreover, 
the lines of stationary phase of the linear PDE
through this saddle point,
along which the real part 
of the phase vanishes,
are given by $\mu_I = \pm \mu_R - \omega_0$,
where $\mu= \mu_R + i \mu_I$.
In the vicinity of the saddle,
analysis along the relevant Stokes and anti-Stokes lines 
shows that all solutions with initial data
at $\mu_0 \le -\omega_0$
stay near the attracting QSS while $\mu$ is negative
and then, also, that they stay near the repelling QSS 
at all points $x$
at least until $\mu$ reaches $+\omega_0$,
to leading order.
More importantly,
application of the methods of stationary phase and steepest descents
yields the formula for the space-time buffer curve,
$(\mu_{\rm stbc}(x),x)$ 
with $\mu_{\rm stbc}(x) \ge \omega_0$ 
for all $x$,
for the general class of source terms considered here.
This space-time buffer curve
represents the $x$-dependent times
at which the particular solution
of the linearised PDE 
ceases to be exponentially small,
and at which the hard onset
of oscillations occurs, to leading order.
Similarly, 
there is a homogeneous exit time curve, $\mu_h(x)$,
along which the homogeneous component 
of the PDE has magnitude one.
Therefore, it is important to determine,
at each point $x$,
which of the two times $\mu_{\rm stbc}(x)$ and $\mu_h(x)$
is smaller, causing the solution of the full
cubic PDE to stop being exponentially small first.
The smaller of the two times marks the duration of the DHB,
to leading order.

After completing this first goal,
we study how the properties of DHB
depend on the outcome of this competition
between $\mu_{\rm stbc}(x)$
and $\mu_h(x)$,
as well as on the properties
of $I_a(x)$.
For solutions with initial data given at any 
$\mu_0 \le  -\omega_0$,
we identify three cases of DHB:
one in which the particular solution
determines the onset of oscillations
at all $x$ in the domain;
another case in which the time of onset
is determined
by the particular solution 
for some intervals on the domain
and by the homogeneous solution 
on the complementary intervals;
and, a third case in which
the homogeneous solution 
determines the exit time at all points.
We examine spatially uni-modal,
spatially periodic,
and smoothed step function source terms.
In these examples,
quantitative agreement
is found 
between the leading order analysis
in all cases of DHB
and the results of direct numerical simulations
of \eqref{eq:gen-CGL}
with zero-flux boundary conditions on $[-\ell,\ell]$ 
and all values of $\ell$ for which simulations were conducted.
Balanced symmetric Strang operator splitting
\cite{SGMS2013}
was used, with centered finite differences
for the Laplacian and fourth-order Runge-Kutta
for the time-stepping.

In addition, we 
show that solutions of \eqref{eq:gen-CGL} can exhibit
spatially-dependent DHB
also when the initial data is given at times
$\mu_0 \in (-\omega_0,-\delta]$
for some small $\delta > 0$.
We label this as Case 4 of DHB.
For these solutions with initial data
given at times 
much closer to the time of the instantaneous
Hopf bifurcation,
the spatial dependence of $\mu_h(x)$
causes the exit time
from a neighborhood of the repelling QSS
to be asymmetric about $\mu=0$ 
with respect to the time of entry
into a neighborhood of the attracting QSS.
This contrasts with the dynamics of DHB
in analytic ODEs,
where the entry-exit function
(also known as the way-in way-out function)
is symmetric for initial conditions
given close to the instantaneous Hopf bifurcation.

After having carried out the analysis
of DHB in the above four cases,
we use a formal analysis to show that solutions
of the full, cubic CGL equation \eqref{eq:gen-CGL}
with initial data given at time $\mu_0\le  -\omega_0$
are also close to the attracting QSS
(now of the full nonlinear PDE)
on $\mu<0$
and remain close to the repelling QSS
after the instantaneous Hopf bifurcation ($\mu=0$)
at least until $\mu=\omega_0$
in Case 1 of DHB.
The nonlinear analysis is carried out using
the integral form of the equation
governing the difference between the solution
of the full nonlinear PDE 
and the particular solution of the linear PDE.
Use of an iterative method 
then establishes the closeness to the repelling QSS 
in the cubic PDE, 
and it reveals how the asymptotic expansion
of the QSS is naturally generated.
The main result is that, to leading order,
solutions of the cubic CGL \eqref{eq:gen-CGL}
stay near the repelling QSS until
the same space-time buffer curve in Case 1 of DHB.
Moreover, we note that,
the situation here 
for DHB in the CGL PDE
is similar in this respect to that 
for DHB in analytic ODEs,
where the linear problem determines 
the buffer point to leading order
in the ODEs,
and the nonlinear terms in the analytic ODE 
(such as the cubic term $-\vert A \vert^2 A$)
only contribute at higher order to DHB.

Finally, we extend  
the main DHB results
for the base case of the PDE \eqref{eq:gen-CGL}
in several directions.
The simplest extension
is to take into account the higher order terms
in the instantaneous Hopf bifurcation curve
for the base case.
To leading order, this curve is given by $\mu=0$ 
in the space-time plane.
The first non-zero correction occurs
at $\mathcal{O}(\eps)$,
and we will study its impact on DHB.
As a second extension,
we study the DHB also in the base case
but now with source terms
which do not satisfy the hypotheses 
imposed on $I_a(x)$ 
in the general analysis, 
namely with an algebraically-growing and a sign-changing function.
Nevertheless, 
for each of these sources,
we also find
good agreement 
between the analytically calculated space-time buffer curve
and the numerically calculated spatially-dependent
times at which the solutions leave a neighborhood of 
the repelling QSS and the oscillations set in
for the nonlinear PDE \eqref{eq:gen-CGL}.

In a third direction, we extend the analysis 
of the base case 
to asymptotically large-amplitude 
($\beta=-\frac{1}{2}$) source terms
in the PDE \eqref{eq:gen-CGL},
while retaining small diffusivity ($\gamma=1$).
Here, the QSS are highly-nontrivial,
and the instantaneous Hopf bifurcation times 
are spatially dependent,
instead of being homogeneous at $\mu=0$,
to leading order.
We explore 
the more complex spatial dependence
of the Hopf bifurcation curve 
and the space-time buffer curve.
An example shows that it is possible to choose
the amplitude and form of the large source term
(and hence of the resulting QSS)
to design even more complex spatio-temporal
onset of oscillations,
giving region-specific control over the onset of oscillations.
In a fourth direction, 
we briefly extend the analysis of the base case 
to an example with
$\mathcal{O}(1)$ diffusivity
and $\mathcal{O}(1)$ amplitude source term,
(i.e., $\gamma=0$ and $\beta=0$
in \eqref{eq:gen-CGL}).
The space-time buffer curve
gets somewhat flattened out
compared to the case
of $\mathcal{O}(\eps)$ diffusivity.

While our primary motivation 
is to carry out this analysis
of DHB in the CGL PDE \eqref{eq:gen-CGL}
and to derive a method that can be used on other
reaction-diffusion (R-D) systems 
known to exhibit DHB \cite{KV2018},
another motivation 
for understanding the phenomenon
of DHB in PDEs
is that some ODE models exhibiting DHB 
are simplified versions 
or conceptual models of more complex phenomena
which involve diffusion and advection.
An example is the Maasch-Saltzman ODE model
of glacial cycles \cite{MS1990},
in which DHB is advanced as a possible mechanism
for the mid-Pleistocene transition
from 40,000 year glacial cycles
to approximately 100,000 year cycles.
See also \cite{EKKV2017}.
Since the Maasch-Saltzman model 
is a useful conceptual ODE model,
one would also like to know 
whether or not the corresponding PDE models, 
such as a more fully developed PDE model 
of the Pleistocene glacial cycles,
can also exhibit DHB.
Otherwise, for these problems,
the phenomenon of DHB 
would only be of more limited interest.
Along with \cite{ADVW2020,KV2018},
this work presents a step in that larger direction,
showing that DHB also occurs
in nonlinear spatially-extended systems.

We observe that the analysis presented
herein builds naturally on the results
known about DHB in analytic ODEs.
In fact, in the case of $d=0$ and
a spatially-homogeneous source term $I_a(x) \equiv I_a$,
system \eqref{eq:gen-CGL} 
reduces to a prototypical form 
of DHB in analytic ODEs.
In this case, we directly recover 
the known DHB results 
for analytic ODEs, 
see \cite{BER1989,HKSW2016,N1987,N1988,S1973}.
The hard onset of oscillations
occurs to leading order at the buffer point 
$\mu=\omega_0$ of the ODEs, and 
it is spatially uniform.
An example is provided 
by the Shishkova equation,
$\eps z_{\mu} = (\mu + i \omega_0) z
- (1+i\alpha)\vert z \vert^2 z 
+ \eps h(\mu)$
for $z(\mu;\eps)$,
with $h(-i \omega_0)\ne 0$
and $0< \eps \ll 1$.
See \cite{S1973}.
This equation
has a family of attracting slow invariant manifolds
for $\mu < -\delta$, 
where $\delta>0$ is small and independent of $\eps$,
and a family of repelling slow invariant manifolds
for $\mu > \delta$.
These families of manifolds
may be extended in to the regions
$\mu>0$ and $\mu < 0$,
respectively,
and they
are exponentially close to each other
on $(-\omega_0,\omega_0)$ to leading order. 
By the theory of DHB in ODEs,
any solution 
which enters a fixed, small neighbourhood
of an attracting slow invariant manifold
at a value of $\mu < -\omega_0$
must exit a neighbourhood 
of the repelling slow invariant manifold
at $\mu = \omega_0$ to leading order,
which is the buffer point.
This delayed loss of stability
occurs at $\mu = \omega_0$ (to leading order)
for these solutions 
independently of how much time 
they have spent spiraling in 
toward the attracting manifold
before $\mu$ reaches $-\omega_0$.
This is because
the term $\eps h(\mu)$ breaks the symmetry
of $z \to -z$ in the Shishkova equation
and because
the intersection 
at $\mu=\omega_0$ 
of the Stokes line through the saddle
(or nilpotent) point at $\mu=-i\omega_0$
with the real axis 
acts as a barrier, or buffer point.
Hence, all solutions 
that have been attracted to that slow manifold
must diverge away 
from the repelling manifold 
along with it,
irrespective of how far back in the past
they approached the attracting manifold.

Finally, we remark that, just as is observed here,
lines of stationary phase and 
lines of steepest ascents and descents 
play central roles
in the asymptotics of solutions
of many linear and nonlinear ODEs and PDEs.
For the general theory of the Stokes phenomenon, see for example
\cite{BO1999,CKA1998,CM2005,D1973,E1983,KC1981,M2006,M1984,O1974}.
We follow naming conventions
in \cite{BO1999,KC1981,M1984}.
Moreover, in many of these equations,
there are multiple components 
which cease to be exponentially small, 
by crossing Stokes lines
for example,
and which transition through modulus one to becoming large.

This article is organised as follows. 
In Section~\ref{sec:lin-CGL},
the linear CGL equation is studied
in the base case
of \eqref{eq:gen-CGL}
with $\beta=1/2$ and $\gamma=1$,
establishing that
solutions with $\mu_0 < -\omega_0$
remain near the attracting and repelling QSSs 
at least until $\mu=+\omega_0$
for all $x$.
In Section~\ref{sec:stbc}, 
the space-time buffer curve $(\mu_{\rm stbc}(x),x)$ is derived,
where $\mu_{\rm stbc}(x) \ge \omega_0$ for
bounded, positive source terms,
showing that DHB occurs in the CGL PDE.
Also, the examples
of different $I_a(x)$
are given 
to show how their extrema
and spatial form (uni-modal, periodic, smoothed step function)
determine the location 
of the space-time buffer curve
and the dynamics of DHB.
In Section~\ref{sec:Ah-transition}, 
we examine the homogeneous solutions
and derive the homogeneous exit time curve,
$\mu_h(x)$.
This establishes the spatially-dependent exit times
caused by the homogeneous components
of the linear solutions,
both when the initial data is
given at a time $\mu_0 \le  -\omega_0$,
and when it is given at a time $\mu_0 \in (-\omega_0,0)$.
In Section~\ref{sec:DHB-fourcases},
the four cases of DHB are classified, with examples.
In Section~\ref{sec:gen-CGL}, 
the analysis of the nonlinear PDE is presented,
establishing
the closeness 
of the solutions to the repelling QSS
of the full cubic CGL equation.
In Section~\ref{sec:Hopf-basecase},
it is shown how the DHB of the solutions
of \eqref{eq:gen-CGL} 
is influenced by the $\mathcal{O}(\eps)$ terms
in the time of the instantaneous Hopf bifurcation.
In Section~\ref{sec:ag+sc},
examples are presented with algebraically-growing
and sign-changing source terms,
to push beyond the analysis of the general system.
Then, the impact of the large-amplitude
source terms ($\beta = -1/2$)
is analyzed in Section~\ref{sec:largeamp-sourceterm}.
Also, the extension of the DHB results and
space-time buffer curve formula
to the case of $\mathcal{O}(1)$ diffusivity
($\gamma=0$) and amplitude source term ($\beta=0$)
is given
in Section~\ref{sec:DO1}.
Conclusions and discussion are presented 
in Section~\ref{sec:CD}.

%------------------------------------------------------------------------------
\section{Linear analysis for solutions with $\mu_0<-\omega_0$}
\label{sec:lin-CGL}
%------------------------------------------------------------------------------

In this section, 
we analyze the linear CGL equation,
obtained by linearising \eqref{eq:gen-CGL} 
about $A=0$,
in the base case 
of moderate-amplitude source terms
($\beta=\frac{1}{2}$) 
and small-amplitude diffusivity 
($\gamma=1$),
\begin{equation} 
\begin{split}
A_t &= (\mu + i \omega_0) A 
+ \sqrt{\eps} I_a(x) + \eps d A_{xx}
\nonumber
\\
\mu_t &=  \eps.
\nonumber
\end{split}
\nonumber
\end{equation}
Equivalently,
the system may be written 
as a scalar equation
for $A(x,\mu)$,
\begin{equation} 
\label{eq:lin-CGL}
\eps A_{\mu} = (\mu + i \omega_0) A 
+ \sqrt{\eps} I_a(x) + \eps d A_{xx}.
\end{equation}
Here, $\mu(t)=\mu_0 + \eps t$,
and we focus 
on data for which 
$\mu_0 \le -\omega_0$.
The other case,
with data given at a time $\mu_0 \in (-\omega_0,0)$,
is analyzed 
in Section~\ref{sec:Ah-transition}.

For $0<\eps \ll 1$,
all solutions with $\mu_0 \le -\omega_0$
rapidly and exponentially approach 
an attracting QSS, 
see \eqref{eq:QSS-basecase} 
(which also contains the terms in the QSS for the full nonlinear equation),
since the real part of the coefficient on the linear term
is negative and stays well bounded 
away from zero for these $\mu<0$.
In this section, we show that 
the solutions with $\mu_0\le -\omega_0$
not only remain close to the attracting QSS
until the time of the instantaneous Hopf bifurcation,
but after the parameter crosses 
the instantaneous Hopf bifurcation
at $\mu=0$ they remain close 
to the repelling QSS as well,
at least until the time $\mu=+\omega_0$
at all points $x$
for the functions $I_a(x)$ we study.
This will be shown using
the classical methods of stationary phase 
and steepest descents,
see \cite{BO1999,KC1981,M1984}.

%------------------------------------------------------------------------------
\subsection{The homogeneous and particular solutions}  \label{sec:lin-CGL-1}
%------------------------------------------------------------------------------

Define the following new dependent variable,
which is based on an integrating factor,
\begin{equation}\label{eq:AB}
B(x,\mu) = A(x,\mu) e^{-\frac{1}{2\eps} (\mu+ i \omega_0)^2}.
\end{equation}
Equation \eqref{eq:lin-CGL}
may then be written as
\begin{equation} \label{eq:B-lin}
\sqrt{\eps} B_{\mu} 
= I_a(x)e^{-\frac{1}{2\eps}(\mu + i \omega_0)^2}
  + \sqrt{\eps} d B_{xx}.
\end{equation} 

By Duhamel's Principle,
the solution consists
of homogeneous and particular components,
\begin{equation} \label{eq:B-Duhamel}
B(x,\mu)  =  B_h (x,\mu) + B_p(x,\mu).
\end{equation}
The homogeneous component satisfies
$(B_h)_\mu = d (B_h)_{xx}$ 
with $B_h(x,\mu_0)=A_0(x)$,
\begin{equation} \label{eq:B_h}
B_h(x,\mu) 
= \frac{e^{-\frac{1}{2\eps}(\mu_0+i \omega_0)^2}}
          {\sqrt{4\pi d (\mu - \mu_0)}}
\int_{\mathbb{R}}  e^{\frac{-(x-y)^2}{4d(\mu-\mu_0)}} A_0(y) dy.
\end{equation}
This homogeneous solution is valid (at least) for all real $\mu > \mu_0$,
recalling that $\mu_0 < -\omega_0$,
and throughout this article,
we will evaluate or estimate it on the ${\rm Re}(\mu)$-axis.
Nevertheless, we note that, 
with the initial data used in the examples,
$B_h(x,\mu)$ is actually analytic
in the complex $\mu$ plane,
excluding the branch point and cut.
Note that, for general initial data,
$A_0(x)$ is Gevrey regular of order $k=\frac{1}{2}$
on an appropriate domain
implies that the homogeneous solution $A_h(x,\mu)$
is analytic,
by standard theory
for homogeneous heat equations.
See for example \cite{R2019},
and recall that a function $f(z)$ 
is Gevrey regular of order $k$ on a set $\vert z \vert < r$,
if there exist positive constants $C_1, C_2$ such that
${\rm max}_{\vert z\vert \le r} \vert \frac{d^k f}{dz^k}(z) \vert
\le C_1 C_2^k (k!)^\frac{1}{2}$.

The particular solution 
satisfies the full linear PDE
\eqref{eq:B-lin},
but with zero initial condition at $\mu_0$, 
\begin{equation} \label{eq:Bp-def}
\begin{split}
B_p(x,\mu) 
&= \frac{1}{\sqrt{\eps}}
\int_{\mu_0}^{\mu} 
g(x,\mu-\tilde{\mu}) 
e^{-\frac{1}{2\eps}(\tilde{\mu} + i \omega_0)^2}
d\tilde{\mu}, \\
g(x,\mu-\tilde{\mu}) &= 
\frac{1}{\sqrt{4\pi d (\mu - \tilde{\mu}) }}
\int_{\mathbb{R}}  e^{\frac{-(x-y)^2}{4d(\mu-\tilde{\mu})}} I_a(y) dy.
\end{split}
\end{equation}
The source terms $I_a(x)$
we study
are such that
$g(x,\mu-{\tilde{\mu}})$
is analytic in a region of the complex plane
which includes the portion of the real axis 
with $\mu > {\tilde{\mu}}$, but which excludes a small neighborhood
of the branch point and cut.
Then, in turn, $B_p(x,\mu)$
is analytic in an appropriate region  
about the segment $(\mu_0,\omega_0)$ 
of the ${\rm Re}(\mu)$-axis.
For more general source terms,
one needs to require that $I_a(x)$ is
Gevrey regular of order $k=\frac{1}{2}$
in a region containing a segment of the ${\rm Re}(x)$-axis,
sufficiently large to guarantee that
$B_p(x,\mu)$ is analytic in $\mu$
for ${\rm Re}(\mu) > {\rm Re}({\tilde{\mu}})$.
This follows 
from standard theory for the analyticity of solutions.
We refer to 
\cite{B2005,LMS1999,R2019}
for the general 
theory of analyticity of solutions
and Gevrey regularity of order $k$
for homogeneous and inhomogeneous heat equations.

In the complex ${\tilde{\mu}}$ plane
(${\tilde{\mu}}= {\tilde{\mu}}_R + i {\tilde{\mu}}_I$),
the phase function
in the integrand of $B_p$ is
\begin{equation}
\label{eq:complexphase}
-\frac{1}{2} ( {\tilde{\mu}} + i \omega_0 )^2 
= \phi + i \psi, \ \ {\rm where}  \ \ 
\phi 
= -\frac{1}{2} \left( {\tilde{\mu}}_R^2 - ({\tilde{\mu}}_I+\omega_0)^2 \right)
\ \ {\rm and} \ \ 
\psi
=  - {\tilde{\mu}}_R ({\tilde{\mu}}_I + \omega_0).
\end{equation}
This phase has 
a saddle point
at ${\tilde{\mu}}=-i\omega_0$,
and the topography induced by this saddle
will play a central role
in the analysis.
The level sets of $\phi$ 
are hyperbolas 
and also known 
as Stokes lines. 
The Stokes lines 
with $\phi=0$
through the saddle 
(which are the asymptotes
of the hyperbolas)
bound the valleys and hills.
They may be parametrised by ${\tilde{\mu}}_R$
via ${\tilde{\mu}}_I = \pm {\tilde{\mu}}_R - \omega_0.$
Also, the level sets of $\psi$
are hyperbolas 
(with asymptotes
given by the axes),
and they are referred to as
anti-Stokes lines.
The geometry
is illustrated in 
Figure~\ref{fig:Contour-Ca}.
See \cite{BO1999,D1973,KC1981,M1984},
for example,
for the general theory of Stokes and anti-Stokes lines.

\bigskip
\noindent
{\bf Remark.}
The integral 
involving $A_0(y)$ in \eqref{eq:B_h}
is the same as the integral
involving $I_a(y)$ in \eqref{eq:Bp-def}(b),
in the case
when $A_0(y)$ is chosen as a multiple of $I_a(y)$
and with $\mu_0$ replacing ${\tilde {\mu}}$.

%------------------------------------------------------------------------------
\subsection{Tracking the particular solutions near the attracting QSS}
\label{sec:lin-CGL-2}
%------------------------------------------------------------------------------

In this section,
we briefly show using steepest descents 
that solutions with $\mu_0 \le -\omega_0$
stay near the attracting QSS
before the instantaneous Hopf bifurcation.
Although the result is of course well known,
it is useful to demonstrate briefly 
how the asymptotic method of steepest descents
naturally finds the asymptotic expansion 
of the attracting QSS.

\begin{figure}[htbp]
   \centering
   \includegraphics[width=5in]{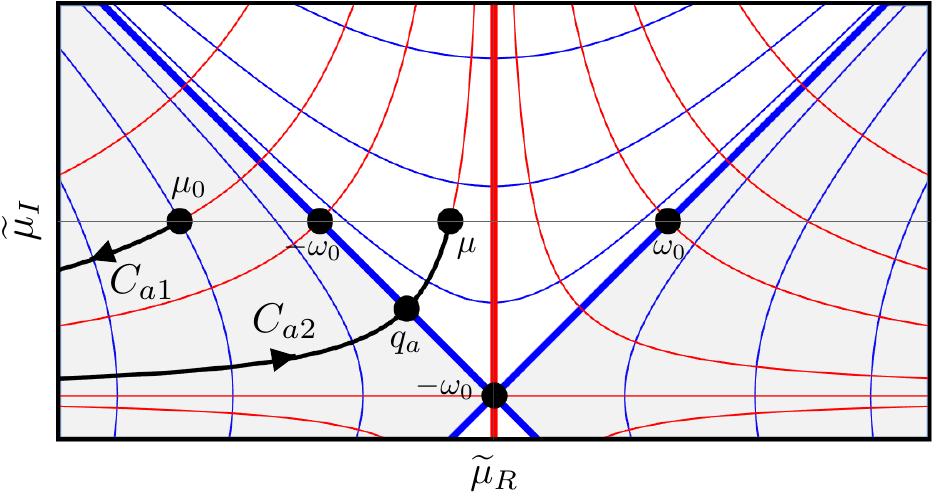}
\caption{
The contour $C_a =C_{a1} \bigcup C_{a2}$
in the complex ${\tilde{\mu}}$ plane.
Here, and in other figures below,
the Stokes lines are the blue curves, and the anti-Stokes lines 
are the red curves. 
The Stokes lines through the saddle (thick blue lines) 
separate the hills (unshaded) and valleys (grey shaded) 
of the phase function. 
	}
	   \label{fig:Contour-Ca}
	\end{figure}

Let $\delta>0$ 
be sufficiently small 
and independent of $\eps$.
We fix an arbitrary value 
of $\mu \in (-\omega_0,-\delta]$,
and track the particular solution $B_p$
from $\mu_0$ to the fixed value $\mu$.
Let $C_a$ denote the contour $ C_{a1} \bigcup C_{a2} $,
where $C_{a1}$
consists of the semi-infinite segment
of the steepest descent curve
$\psi \equiv -\omega_0\mu_0$
from the point $\mu_0$ on the real axis
out to infinity  ($-\infty-i\omega_0$)
toward the horizontal asymptote
${\tilde{\mu}}_I = -\omega_0$,
and $C_{a2}$
consists of the semi-infinite segment
of the steepest ascent curve
$\psi \equiv -\omega_0\mu$
from infinity back up to the fixed value
$\mu$ on the real axis.
Note that $C_{a2}$ crosses 
the Stokes line $\mu_{I} = -\mu_R - \omega_0$
at the point
$q_a=-\sqrt{-\omega_0\mu} + i (\sqrt{-\omega_0\mu}-\omega_0)$.
See Figure~\ref{fig:Contour-Ca}.

	We track $B_p$ along $C_a$ 
	from $\mu_0$ 
	to the fixed value $\mu$ 
	on $(-\omega_0,-\delta]$.
	By \eqref{eq:Bp-def},
	the solution is
	\begin{equation}
	\label{eq:Bp-at-mu-Ca}
	\begin{split}
	B_p(x,\mu) 
	& = \mathcal{I}_{a1} + \mathcal{I}_{a2}, \quad {\rm where} \\
	\mathcal{I}_{ai} &=
	\frac{1}{\sqrt{\eps}}
	\int_{C_{ai}}
	g(x,\mu - {\tilde{\mu}}) 
	e^{-\frac{1}{2\eps} ({\tilde{\mu}} + i \omega_0)^2}
	d{\tilde{\mu}}, \ \ i=1,2, 
	\qquad \mu \in (-\omega_0,-\delta].
	\end{split}
	\end{equation}
	Directly from $\phi$,
	the real part 
	of the complex phase \eqref{eq:complexphase},
	and the analyticity of $g$,
	one finds 
	\begin{equation}\label{eq:Ia1}
	\vert \mathcal{I}_{a1} \vert
	\le C e^{-\frac{1}{2\eps}(\mu_0^2-\omega_0^2)},
	\quad {\rm for} \ {\rm  some} \  C>0.
	\end{equation}

	The main work then
	is to derive the result for $\mathcal{I}_{a2}$,
	which we do using two different parametrisations of $C_2$,
	explicitly here using ${\tilde{\mu}}_R$,
	and implicitly in Appendix~\ref{sec:App-A},
	\begin{equation}\label{eq:Ca2-muR}
	C_{a2}:  \quad
	{\tilde {\mu}} 
	= {\tilde{\mu}}_R 
	  + i {\tilde{\mu}}_I({\tilde{\mu}}_R),
	\quad {\rm with} \ \ 
	{\tilde{\mu}}_I({\tilde{\mu}}_R) 
	= -\omega_0 \left( 1 - \frac{\mu}{{\tilde{\mu}}_R} \right).
	\end{equation}
	Along $C_{a2}$,
	${\tilde{\mu}}_R$ increases
	from $-\infty$ to $\mu$,
	and ${\tilde{\mu}}_I$ increases
	from $-\omega_0$ to zero.
	Also, we observe that,
	by \eqref{eq:Ca2-muR},
	$\frac{d{\tilde{\mu}}}{d{\tilde{\mu}}_R} 
	= 1 - \frac{i\omega_0\mu}{{\tilde{\mu}}_R^2}$.
	Hence,
	\begin{equation}
	\mathcal{I}_ {a2}
	= \frac{1}{\sqrt{\eps}}
	\int_{-\infty}^{\mu}
	g\left( x,\mu-{\tilde{\mu}}_R 
		  +i \omega_0 \left(1 - \frac{\mu}{{\tilde{\mu}}_R}\right) \right)
	e^{-\frac{1}{2\eps}
	    \left({\tilde{\mu}}_R+i\frac{\omega_0\mu}{{\tilde{\mu}}_R}\right)^2}
	\left( 1 - \frac{i \omega_0 \mu}{{\tilde{\mu}}_R^2} \right)
	d{\tilde{\mu}}_R.
	\nonumber
	\end{equation}
	Now, it is useful to define
	\begin{equation}\label{eq:ibp-uv}
	u({\tilde{\mu}}_R) 
	= -\eps
	    \frac{g\left( x,\mu-{\tilde{\mu}}_R 
		  +i\omega_0 \left(1 - \frac{\mu}{{\tilde{\mu}}_R}\right) \right)}
	       {\left( {\tilde{\mu}}_R + \frac{i\omega_0\mu}{{\tilde{\mu}}_R} \right)}
\qquad
	v({\tilde{\mu}}_R) 
	= e^{-\frac{1}{2\eps}
	    \left({\tilde{\mu}}_R+i\frac{\omega_0\mu}{{\tilde{\mu}}_R}\right)^2}.
	\end{equation}
	Hence,
	using integration by parts,
	formula \eqref{eq:Bp-def}(b) for $g$,
	and $\lim_{\chi \to 0^+} g(x,\chi) = I_a(x)$,
	we find
	\begin{equation}
	\begin{split}
	\mathcal{I}_{a2}
	&= \frac{1}{\sqrt{\eps}}
	\int_{-\infty}^{\mu} 
	u \frac{dv}{d{\tilde{\mu}}_R} d{\tilde{\mu}}_R \\
	&= 
	-\sqrt{\eps}
	\frac{I_a(x)e^{-\frac{1}{2\eps}(\mu+i\omega_0)^2}}
	     {(\mu+i\omega_0)} 
	-\frac{1}{\sqrt{\eps}}
	\int_{-\infty}^{\mu} 
	v \frac{du}{d{\tilde{\mu}}_R} d{\tilde{\mu}}_R 
	\end{split}
	\end{equation}
	Proceeding to higher order,
	using $g_{\mu} = d g_{xx}$ and 
	$\lim_{\chi \to 0^+} g_{xx} (x,\chi) 
	= {I_a}''(x)$, 
	we find
	\begin{equation} \label{eq:Ia2-final}
	\mathcal{I}_{a2}
	=\left\{
	-\sqrt{\eps} \frac{I_a(x)}{\mu+i\omega_0}
	+\eps^{\frac{3}{2}} 
	   \frac{\left[ I_a(x) + d {I_a}''(x) (\mu+i\omega_0) \right]}
		{(\mu+i\omega_0)^3} 
	+\mathcal{O}\left(\frac{\eps^{\frac{5}{2}}}{(\mu+i\omega_0)^5}\right) \right\}
	e^{-\frac{1}{2\eps}(\mu+i\omega_0)^2}.
	\end{equation}
	Summing \eqref{eq:Ia1} and \eqref{eq:Ia2-final}
	and recalling $\mu_0 \le  -\omega_0 < \mu \le -\delta$,
	we have
	\begin{equation}
	\label{eq:Ia1+Ia2}
	B_p(x,\mu)=\mathcal{I}_{a1} + \mathcal{I}_ {a2}
	= \left\{ 
	-\frac{\sqrt{\eps} I_a(x)}{\mu + i \omega_0}
	+   \frac{
	\eps^{\frac{3}{2}}
	\left[ I_a(x) + d {I_a}''(x) (\mu+i\omega_0) \right]}
		{(\mu+i\omega_0)^3}
	+ \mathcal{O}\left(\frac{\eps^{\frac{5}{2}}}
				{(\mu+i\omega_0)^5} \right) 
	\right\}
	e^{-\frac{1}{2\eps}(\mu+i\omega_0)^2} .
	\end{equation}

	Finally, we translate the formula
	back to the $A$ equation 
	using \eqref{eq:AB},
	\begin{equation}\label{eq:Ap-QSS-a}
	A_p(x,\mu) 
	= -\sqrt{\eps} \frac{I_a(x)}{\mu + i \omega_0}
	  +{\eps}^\frac{3}{2} 
	    \frac{\left[I_a(x) + d {I_a}''(x)(\mu + i\omega_0)\right]}
		 {(\mu + i \omega_0)^3}
	+ \mathcal{O}\left(\frac{\eps^{\frac{5}{2}}}{(\mu+i\omega_0)^{5}}\right),
	\quad -\omega_0 < \mu \le  -\delta.
	\end{equation}
	The first and second terms here
	are exactly the first and second order terms
	in the asymptotic expansion
	of the attracting QSS for the linear CGL;
	cf. \eqref{eq:QSS-basecase},
	where the expansion is given for the cubic CGL.
	(Note that the cubic term
	$-\vert A\vert^2 A$ in \eqref{eq:gen-CGL}
	gives rise to an additional
	term at $\mathcal{O}(\eps^{\frac{3}{2}})$
	given by $-\frac{(1+i\alpha)I_a^3(x)}{(\mu+i\omega_0)^2(\mu^2+\omega_0^2)}$
	in the asymptotics of the QSS, see \eqref{eq:QSS-basecase}.)
	The remainder, which is uniform in $x$,
        contains $\mathcal{I}_{a1}$ and
	the higher order terms
	in the asymptotic expansion of the attracting QSS,
	and one may continue 
	using integration by parts
	on $\mathcal{I}_{a2}$
	to derive them.
	The remainder terms also include
an exponentially small term
coming from the attraction
of the initial data to the QSS.
Therefore,
we have shown
that, with $\mu_0 < - \omega_0$,
the particular solution $A_p$
is close to the attracting QSS
for all $\mu \in (-\omega_0,-\delta]$.

%---------------------------------------------------------
\subsection{Tracking the particular solutions near the repelling QSS}
\label{sec:lin-CGL-3}
%---------------------------------------------------------

In this section,
we track the solutions
with initial data given at $\mu_0 \le -\omega_0$
beyond the instantaneous Hopf bifurcation point at $\mu = 0$
into the regime where $\mu>0$.
We show that for any $\mu \in [\delta,\omega_0]$
these solutions are close to the repelling QSS
at all points $x$.
We use 
the method of
stationary phase,
as well as steepest descents,
taking advantage of the saddle point
at ${\tilde{\mu}}=-i \omega_0$
in \eqref{eq:complexphase}.

\begin{figure}[htbp]
   \centering
   \includegraphics[width=5in]{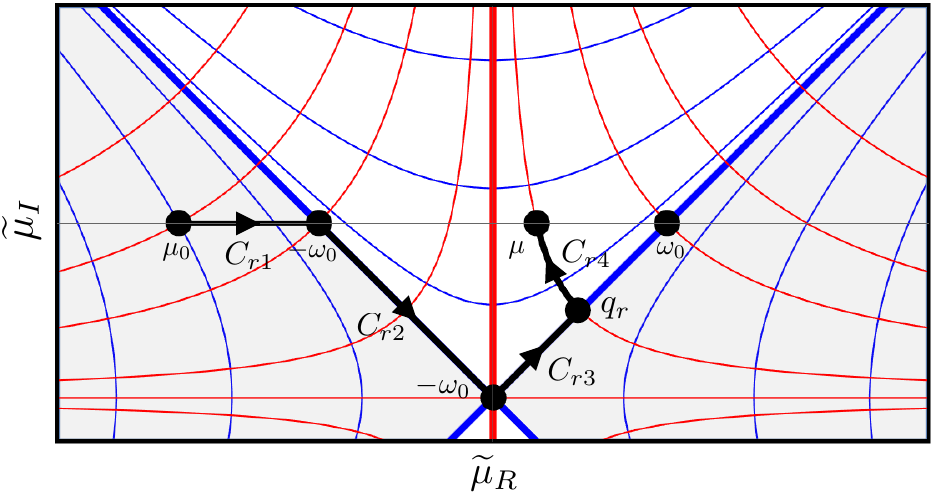}
   \caption{
The contour $C_r = C_{r1} \bigcup C_{r2} \bigcup C_{r3} \bigcup C_{r4}$
in the complex ${\tilde{\mu}}$ plane.
}
   \label{fig:Contour-Cr}
\end{figure}

We fix an arbitrary value 
of $\mu \in [\delta,\omega_0]$.
We consider the contour
$C_r = C_{r1} \bigcup C_{r2} \bigcup C_{r3} \bigcup C_{r4}$,
where $C_{r1} = [\mu_0,-\omega_0]$;
$C_{r2}$ is the segment of the Stokes line
${\tilde{\mu}}_{I} = - {\tilde{\mu}}_R - \omega_0$
from $-\omega_0$
down to the saddle point at $-i\omega_0$;
$C_{r3}$ consists of the segment
of the Stokes line 
${\tilde{\mu}}_{I} = {\tilde{\mu}}_R - \omega_0$
from the saddle point at $-i\omega_0$
up to the point 
$q_r=\sqrt{\omega_0\mu} + i (\sqrt{\omega_0\mu}-\omega_0)$,
for this fixed value of $\mu$;
and, $C_{r4}$
consists of the segment of the steepest ascent
curve $\psi = -\omega_0\mu$
from $q_r$ up to the point $\mu$.
See Figure~\ref{fig:Contour-Cr}.

We take any solution $B_p$
with $\mu_0 \le -\omega_0$
on or near 
the attracting QSS,
and we track it along $C_r$ 
to the fixed value $\mu$.
At that point,
the solution is
\begin{equation}
\label{eq:Bp-at-qr}
\begin{split}
B_p(x,\mu) 
& = \mathcal{I}_{r1} + \mathcal{I}_{r2} + \mathcal{I}_{r3} + \mathcal{I}_{r4}, \quad {\rm where} \\
\mathcal{I}_{ri} &=
\frac{1}{\sqrt{\eps}}
\int_{C_{ri}}
g(x,\mu - {\tilde{\mu}}) 
e^{-\frac{1}{2\eps} ({\tilde{\mu}} + i \omega_0)^2}
d{\tilde{\mu}}, \ \ i=1,2,3,4,
\quad \delta \le \mu\le\omega_0.
\end{split}
\end{equation}
The integral $\mathcal{I}_{r1}$
along the segment
$C_{r1} = [\mu_0,-\omega_0]$
may be evaluated in the same manner
as used for the integral
\eqref{eq:Ia2-final}.
However, here $\mu \in [\delta,\omega_0]$
and the Taylor expansion 
is about ${\tilde{\mu}}=-\omega_0$,
\begin{equation} \label{eq:Ir1-asymptotics}
\mathcal{I}_{r1} 
=\left\{
\frac{\sqrt{\eps}}{2\omega_0} (1+i) g(x,\mu+\omega_0)
+\mathcal{O}\left(\eps^{\frac{3}{2}}\right) \right\}
e^{\frac{i \omega_0^2}{\eps}}.
\end{equation}

Next, 
we parametrise $C_{r2}$ 
by ${\tilde{\mu}}_R$,
with ${\tilde{\mu}}_I ({\tilde{\mu}}_R)
= - ({\tilde{\mu}}_R + \omega_0)$,
and ${\tilde{\mu}}_R: -\omega_0\to 0$.
Hence,
$-\frac{1}{2}(\tilde{\mu} + i \omega_0)^2= i{\tilde{\mu}}_R^2$
for all ${\tilde{\mu}}$ on $C_{r2}$.
It is purely imaginary,
corresponding to the fact 
that $C_{r2}$ 
lies on a Stokes line $\phi=0$.
Hence, for each ${\tilde{\mu}}$ on $C_{r2}$, 
the integrand 
is of the form 
to which the method of stationary phase applies,
namely $g \cdot e^{\frac{i}{\eps}h({\tilde{\mu}}_R)}$
with $h({\tilde{\mu}}_R)={{\tilde{\mu}}_R}^2$.
Moreover, the end point 
$\tilde{\mu}=-i\omega_0$ 
of $C_{r2}$
(${\tilde{\mu}}_R=0$)
is a point of stationary phase,
since
$h'(0)=0$
and $h''(0)=2 \ne 0$,
and it is the only such point along $C_{r2}$.
[For the general method 
in which an end point is a saddle (or turning) point,
see for example Section 4.1 of \cite{M1984},
especially formula (4.14).]

Applying the method of stationary phase,
we insert the parametrisation
of $C_{r2}$,
use ${\tilde{\nu}}=-{\tilde{\mu}}_R$,
Taylor expand about 
${\tilde{\nu}}=0$ 
({\it i.e.,} ${\tilde{\mu}}=-i\omega_0$),
and observe that the dominant contribution
asymptotically comes 
from the point of stationary phase 
at the saddle, 
\begin{equation}\label{eq:Ir2-asymptotics}
\begin{split}
\mathcal{I}_{r2}
&= \frac{1}{\sqrt{\eps}} 
\int_{-\omega_0}^0 
g(x,\mu - [ {\tilde{\mu}}_R - i({\tilde{\mu}}_R + \omega_0) ])
e^{\frac{i}{\eps} {\tilde{\mu}}_R^2}
(1-i) d{\tilde{\mu}}_R \\
&= \frac{1}{\sqrt{\eps}} 
(1-i) 
\int_0^{\omega_0} 
g(x,\mu + [ {\tilde{\nu}} + i(-{\tilde{\nu}} + \omega_0) ])
e^{\frac{i}{\eps} {\tilde{\nu}}^2}
d{\tilde{\nu}} \\
&= \frac{1}{\sqrt{\eps}} 
(1-i) 
g(x,\mu+i\omega_0) 
\int_0^{\omega_0} 
(1 + \mathcal{O}({\tilde{\nu}})) 
e^{\frac{i}{\eps} {\tilde{\nu}}^2}
d{\tilde{\nu}} \\
&= \sqrt{\frac{\pi}{2 h''(0)}} 
e^{\frac{i\pi}{4}} (1-i)
g(x,\mu+i\omega_0) 
+ \mathcal{O}(\sqrt{\eps}) \\
&= \sqrt{\frac{\pi}{2}} g(x,\mu+i\omega_0) 
+ \mathcal{O}(\sqrt{\eps}),
\quad
{\rm for} \ {\rm any} \ \mu \in
[\delta,\omega_0].
\end{split}
\end{equation}
This leading order term in $\mathcal{I}_{r2}$
will turn out to be half
of the leading order term 
in the total integral for $B_p(x,\mu)$ 
for each $\mu \in [\delta,\omega_0]$.

Next, we show that $\mathcal{I}_{r3}$
gives the other half of the leading order term
in $B_p$.
By the definition of $C_{r3}$, 
for any $\mu \in [\delta,\omega_0]$,
we have
$\mathcal{I}_{r3}
= \frac{1}{\sqrt{\eps}}
  \int_{-i\omega_0}^{q_r} 
  g(x,\mu - {\tilde{\mu}})
  e^{-\frac{1}{2\eps} ({\tilde{\mu}} + i \omega_0)^2} d{\tilde{\mu}}.$
We also use ${\tilde{\mu}}_R$ 
to parametrise the segment $C_{r3}$
as ${\tilde{\mu}}_I = {\tilde{\mu}}_R - \omega_0$,
now with ${\tilde{\mu}}_R: 0 \to \sqrt{\omega_0\mu}$.
Hence,
$-\frac{1}{2}(\tilde{\mu} + i \omega_0)^2= -i{\tilde{\mu}}_R^2$
along $C_{r3}$; and,
for each ${\tilde{\mu}}$ on $C_{r3}$, 
the integral is also of the form 
to which the method of stationary phase applies,
namely 
$g \cdot e^{\frac{i}{\eps}{\tilde h}(\tilde{\mu}_R)}$,
with ${\tilde h}({\tilde{\mu}_R})=-{\tilde{\mu}}_R^2$.
Moreover, the initial point 
$\tilde{\mu}=-i\omega_0$ of $C_{r3}$
(${\tilde{\mu}}_R=0$)
is a point of stationary phase,
since
${\tilde h}'(0)=0$ 
and ${\tilde h}''(0)=-2 \ne 0$,
and it is the only such point along $C_{r3}$.
We find
\begin{equation}\label{eq:Ir3-asymptotics}
\begin{split}
\mathcal{I}_{r3}
&= \frac{1}{\sqrt{\eps}} 
\int_0^{\sqrt{\omega_0\mu}} 
g(x,\mu - [ {\tilde{\mu}}_R + i({\tilde{\mu}}_R - \omega_0) ])
e^{-\frac{i}{\eps} {\tilde{\mu}}_R^2}
(1+i) d{\tilde{\mu}}_R \\
&=\frac{1}{\sqrt{\eps}} g(x,\mu+i\omega_0) 
\int_0^{\sqrt{\omega_0\mu}} e^{-\frac{i}{\eps} {\tilde{\mu}}_R^2}
(1 + \mathcal{O}({\tilde{\mu}}_R)) (1+i) d{\tilde{\mu}}_R  \\
&= \sqrt{\frac{\pi}{2 \vert {\tilde h}''(0)\vert}} g(x,\mu+i\omega_0) 
e^{-\frac{i\pi}{4}}(1+i)
+ \mathcal{O}(\sqrt{\eps}) \\
&= \sqrt{\frac{\pi}{2}} g(x,\mu+i\omega_0) 
+ \mathcal{O}(\sqrt{\eps}),
\quad
{\rm for}  \ {\rm any} \ \mu \in
[\delta,\omega_0].
\end{split}
\end{equation}

Finally, we calculate $\mathcal{I}_{r4}$.
Implicitly parametrise $C_{r4}$
using $\sigma$,
\begin{equation}
C_{r4}: \quad
-\frac{1}{2} ({\tilde{\mu}} + i \omega_0)^2 
= -\frac{1}{2} (\mu + i \omega_0)^2  + \sigma.
\nonumber
\end{equation}
The parameter $\sigma$ starts from 
$-\frac{1}{2}(\omega_0^2 - \mu^2)$ 
at the point $q_r$
and increases monotonically along $C_{r4}$ 
to zero at the point $\mu$.
The explicit representation is
\begin{equation}
{\tilde{\mu}}(\sigma) =
-i\omega_0 + \left[ (\mu + i \omega_0)^2 - 2\sigma \right]^{\frac{1}{2}}.
\nonumber
\end{equation}
The integration along $C_{r4}$
follows in a manner similar
to that along $C_{a2}$ in Appendix \ref{sec:App-A},
except that here one starts at $q_r$
and also here $\mu>0$,
\begin{equation}
\mathcal{I}_{r4}
= \sqrt{\eps} e^{-\frac{1}{2\eps}(\mu+i\omega_0)^2}
\int_{-\frac{1}{2}(\omega_0^2 - \mu^2)}^{0}
g(x,\mu+i\omega_0- [(\mu+i\omega_0)^2 - 2\sigma]^{\frac{1}{2}})
e^{\frac{\sigma}{\eps}}
\left[(\mu + i \omega_0)^2 - 2\sigma
\right]^{-\frac{1}{2}} d\sigma.
\nonumber
\end{equation}
Then, with a similar Taylor expansion,
one finds
\begin{equation}\label{eq:Ir4-result}
\mathcal{I}_{r4}
= \left[
-\frac{\sqrt{\eps} I_a(x)}
      {\mu+i\omega_0} 
+ \eps^{\frac{3}{2}} 
     \left(\frac{I_a(x)+d(\mu+i\omega_0){I_a}''(x)}
                {(\mu+i \omega_0)^3}
     \right) 
+ \mathcal{O}\left(\frac{\eps^{\frac{5}{2}}}
                        {(\mu+i\omega_0)^5}\right) 
\right]
e^{-\frac{1}{2\eps}(\mu + i \omega_0)^2}.
\end{equation}

Summing 
\eqref{eq:Ir1-asymptotics},
\eqref{eq:Ir2-asymptotics},
\eqref{eq:Ir3-asymptotics},
and \eqref{eq:Ir4-result},
we have
\begin{equation}
\label{eq:Bp-mu-sec23}
\begin{split}
B_p(x,\mu)
&=\mathcal{I}_{r1} + \mathcal{I}_{r2} + \mathcal{I}_{r3} + \mathcal{I}_{r4} \\
&= \left[ - \frac{\sqrt{\eps} I_a(x)}{\mu + i \omega_0}
 + \eps^{\frac{3}{2}} \left(\frac{I_a(x)+d(\mu+i\omega_0){I_a}''(x)}
              {(\mu+i \omega_0)^3}
        \right) \right] 
e^{-\frac{1}{2\eps}(\mu+i\omega_0)^2} \\
&+ \mathcal{O}\left(\frac{\eps^{\frac{5}{2}}
e^{-\frac{1}{2\eps}(\mu+i\omega_0)^2}}{(\mu+i\omega_0)^5}\right)
+ \sqrt{2\pi} g(x,\mu+i\omega_0)
 + \mathcal{O}(\sqrt{\eps}),
\quad {\rm for} \ \delta \le \mu \le \omega_0.
\end{split}
\end{equation}
Finally, we translate the formula
back to the $A$ equation 
using \eqref{eq:AB},
\begin{equation}
\label{eq:Ap-sec23}
\begin{split}
A_p(x,\mu) 
&= - \sqrt{\eps} \frac{I_a(x)}{\mu + i \omega_0}
+ \eps^{\frac{3}{2}} \left( \frac{I_a(x)+d(\mu+i\omega_0){I_a}''(x)}
              {(\mu+i \omega_0)^3}
        \right)
+ \mathcal{O}\left(\frac{\eps^{\frac{5}{2}}}{(\mu+i\omega_0)^5}\right) \\
&+ \left( \sqrt{2\pi} g(x,\mu+i\omega_0)
+\mathcal{O}(\sqrt{\eps}) \right)
e^{\frac{1}{2\eps}(\mu+i\omega_0)^2},
\quad {\rm for} \ \delta\le\mu\le\omega_0.
\end{split}
\end{equation}
The first and second terms
are precisely the leading order terms
in the expansion 
of the repelling QSS for the linear CGL
(cf. \eqref{eq:QSS-basecase},
where the QSSs 
are given for the cubic CGL equation).
The third term contains the higher order terms
in the asymptotic expansion of the repelling QSS,
and continued integration by parts
will yield them.

The fourth term 
is exponentially small for 
$\mu \in [\delta,\omega_0 -K\eps^r)$,
for some $K>0$ and any $0<r<1$.
It is a classic Stokes type term.
This term is not 
in the expansion \eqref{eq:Ap-QSS-a} 
of the attracting QSS (on $\mu < -\delta$)
to all orders
or in the expansion
of the repelling QSS (on $\mu>\delta$)
to all orders.
Rather, it is beyond all orders,
$\mathcal{O}\left(e^{-\frac{\omega_0^2}{2\eps}}\right)$,
arising naturally from tracking solutions
on (and near) the attracting QSS
along a contour over the saddle point 
in the complex $\mu$ plane
and into the regime of ${\rm Re} (\mu)>0$.
It is a measure of the exponentially small distance
between the attracting and repelling QSS at $\mu=0$.
In Section~\ref{sec:stbc},
we will determine when it becomes
$\mathcal{O}(1)$ (and then exponentially large).

Overall, therefore,
formulas \eqref{eq:Ap-QSS-a}
and \eqref{eq:Ap-sec23} 
give the asymptotics of solutions $A_p(x,\mu)$
for all $\mu \in (-\omega_0,-\delta]$
and all $\mu \in [\delta,\omega_0]$,
respectively.
They show that,
for all $x$,
the solutions of the linear CGL equation
with Gevrey regular data $A_0(x)$
given at any $\mu_0 \le -\omega_0$
are near the attracting QSS 
until the Hopf bifurcation;
and then, once ${\tilde{\mu}}$ has become positive,
they remain near the repelling QSS 
at least until ${\tilde{\mu}} = \omega_0$
to leading order.
This completes the analysis of this subsection.
Note that the solutions are
also close to the QSS on $(-\delta,\delta)$,
as shown in Appendix~\ref{sec:App-B}.

%------------------------------------------------------------------------------
\section{The space-time buffer curve}
\label{sec:stbc}
%------------------------------------------------------------------------------

In this section,
we track the 
particular solution 
$B_p(x,\mu)$ (and hence also $A_p(x,\mu)$ via \eqref{eq:AB})
from the initial time $\mu_0$,
satisfying $\mu_0 \le -\omega_0$,
to an $x$-dependent, 
maximal value of $\mu$, beyond $\mu=\omega_0$.
For each $x$,
this maximal value, labeled $\mu_{\rm stbc}(x)$,
denotes the space-dependent value of $\mu$
at which $\vert A_p(x,\mu) \vert = 1$,
{\it i.e.,} at which the real part
of the space-time dependent
phase of $A_p$ first vanishes.
To leading order,
the space-dependent time 
$A_p(x,\mu)$ 
is exponentially small for $\mu \in (\mu_0, \mu_{\rm stbc}(x))$
and then transitions to being
exponentially large for $\mu > \mu_{\rm stbc}(x)$.
Hence, at each $x$,
$\mu_{\rm stbc}(x)$ is the maximum of $\mu$ 
for which solutions 
with initial data at $\mu_0 \le -\omega_0$
can remain near the repelling QSS.
We label the union of $\mu_{\rm stbc}(x)$
over all $x$
as the space-time buffer curve.

%------------------------------------------------------------------------------
\subsection{Derivation of the space-time buffer curve, $\mu_{\rm stbc}(x)$}
\label{sec:stbc-1}
%------------------------------------------------------------------------------

From the result of Section~\ref{sec:lin-CGL-3},
we know that the curve 
lies to the right of $\mu=\omega_0$,
since solutions with $\mu_0 \le -\omega_0$
remain close to the repelling QSS
at least until $\mu=\omega_0$
to leading order
at all points $x$.
To track the solutions past this value,
we again deform the contour 
in the complex $\mu$ plane,
taking advantage of the saddle point
at ${\tilde{\mu}}=-i \omega_0$
in \eqref{eq:complexphase}.
Several of the calculations needed here
follow directly from those performed above.

\begin{figure}[h!]
   \centering
   \includegraphics[width=5in]{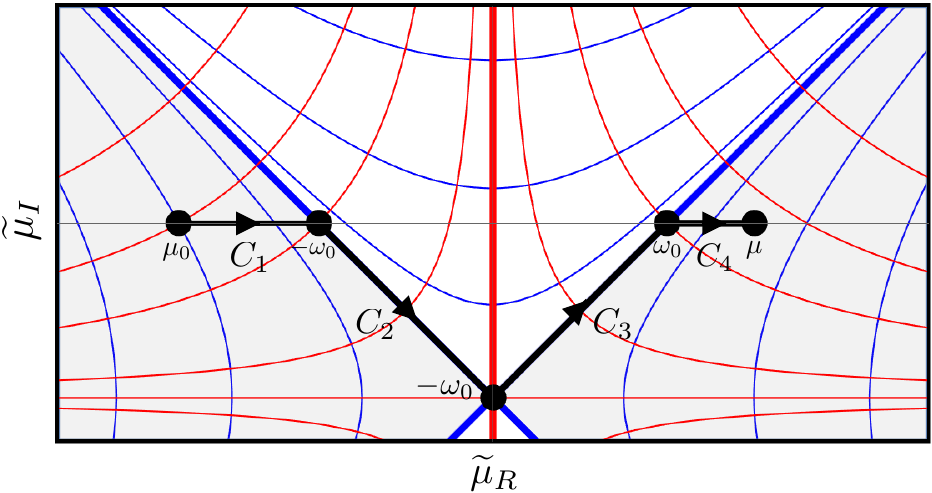}
   \caption{
The contour $C = C_1 \cup C_2 \cup C_3 \cup C_4$
in the complex ${\tilde{\mu}}$ plane.
}
   \label{fig:Contours}
\end{figure}

Let the contour 
\begin{equation}
C = C_1 \bigcup C_2 \bigcup C_3 \bigcup C_4
\nonumber
\end{equation}
consist of the following four segments:
$C_1=[\mu_0,-\omega_0]$,
along the negative real axis;
$C_2$ is the segment of the Stokes line
${\tilde{\mu}}_I = - {\tilde{\mu}}_R - \omega_0$ from $-\omega_0$
down to the saddle at $-i\omega_0$; 
$C_3$ is the segment of the other Stokes line
${\tilde{\mu}}_I =  {\tilde{\mu}}_R - \omega_0$ 
from the saddle at $-i\omega_0$
up to $\omega_0$; and,
$C_4=[\omega_0, \mu_{\rm stbc}(x))$,
along the positive real axis.
See Figure~\ref{fig:Contours}.
Note that
$C_1 = C_{r1}$
and $C_2 = C_{r2}$.

From \eqref{eq:Bp-def}(a)
and the composition of the contour $C$,
one finds
for any arbitrary time $\mu$ on $C_4$,
\begin{equation}\label{eq:Bp-decomp}
\begin{split}
B_p(x,\mu) 
&= \mathcal{I}_1 + \mathcal{I}_2 + \mathcal{I}_3 + \mathcal{I}_4, 
\qquad {\rm where} \ \\
\mathcal{I}_j &= 
  \frac{1}{\sqrt{\eps}}
  \int_{C_j} 
  g(x,\mu - {\tilde{\mu}})
  e^{-\frac{1}{2\eps} ({\tilde{\mu}} + i \omega_0)^2} d{\tilde{\mu}}, \qquad
j=1, 2, 3, \ \\
\mathcal{I}_4 &= 
  \frac{1}{\sqrt{\eps}}
  \int_{\omega_0}^{\mu} 
  g(x,\mu - {\tilde{\mu}})
  e^{-\frac{1}{2\eps} ({\tilde{\mu}} + i \omega_0)^2} d{\tilde{\mu}},
\quad
\mu \in C_4.
\end{split}
\end{equation}
For each $j$,
the value of $\mu$ 
in the integrand along $C_j$
is the same fixed value on $C_4$.

Here, 
since $C_1=C_{r1}$,
one finds by
recalling \eqref{eq:Ir1-asymptotics},
\begin{equation}\label{eq:I1-final}
\mathcal{I}_1= \mathcal{I}_{r1}
= \sqrt{\eps} 
      e^{\frac{i\omega_0^2}{\eps}}
      \left( \frac{1+i}{2\omega_0} \right) g(x,\mu+\omega_0)
      + \mathcal{O}(\eps^{\frac{3}{2}}).
\end{equation}
Next, 
since $C_2=C_{r,2}$,
one finds by
recalling \eqref{eq:Ir2-asymptotics},
\begin{equation}\label{eq:I2-asymptotics}
\mathcal{I}_2
=\mathcal{I}_{r2} 
=\sqrt{\frac{\pi}{2}} g(x,\mu+i\omega_0) + \mathcal{O}(\sqrt{\eps}).
\end{equation}
Then, from the definition of $C_3$, 
we have
\begin{equation}
\mathcal{I}_3 
= \frac{1}{\sqrt{\eps}}
  \int_{-i\omega_0}^{\omega_0} 
  g(x,\mu - {\tilde{\mu}})
  e^{-\frac{1}{2\eps} ({\tilde{\mu}} + i \omega_0)^2} d{\tilde{\mu}} .
\nonumber
\end{equation}
This integral may be evaluated
in the same manner using stationary phase
as that in $\mathcal{I}_{r3}$,
except here one integrates 
all the way up to $\omega_0$
(instead of stopping at $q_r$),
\begin{equation}\label{eq:I3-asymptotics}
\begin{split}
\mathcal{I}_3 
&= \frac{1}{\sqrt{\eps}} 
\int_0^{\omega_0} 
g(x,\mu - [ {\tilde{\mu}}_R + i({\tilde{\mu}}_R - \omega_0) ])
e^{-\frac{i}{\eps} {\tilde{\mu}}_R^2}
(1+i) d{\tilde{\mu}}_R \\
&= \sqrt{\frac{\pi}{2}} g(x,\mu+i\omega_0) 
+ \mathcal{O}(\sqrt{\eps}).
\end{split}
\end{equation}
The contribution
along the segment $C_4$ is
$\mathcal{I}_4
= \frac{1}{\sqrt{\eps}}
  \int_{\omega_0}^\mu
  g(x,\mu - {\tilde{\mu}})
  e^{-\frac{1}{2\eps} ({\tilde{\mu}} + i \omega_0)^2} d{\tilde{\mu}}.$
The contour integral may be estimated
in a manner similar to that used 
along $C_1$,
\begin{equation}\label{eq:I4-final}
\vert \mathcal{I}_4 \vert \le \mathcal{O}(\sqrt{\eps}).
\end{equation}

Then, substituting the results 
for $\mathcal{I}_1, 
\mathcal{I}_2, 
\mathcal{I}_3,$
and $\mathcal{I}_4$ 
(see \eqref{eq:I1-final} --
\eqref{eq:I4-final})
into \eqref{eq:Bp-decomp},
we find
\begin{equation}\label{eq:Bp-stationaryphase}
B_p (x,\mu) 
=  \sqrt{2\pi} g(x,\mu+i\omega_0) 
+\mathcal{O}(\sqrt{\eps}),\quad
{\rm for} \ {\rm any} \ \mu \in C_4.
\end{equation}
This represents the particular solution
valid for all $\mu$ on $C_4$.

Finally, using the change of variables
\eqref{eq:AB},
we see from \eqref{eq:Bp-stationaryphase} 
that the particular solution 
of the linear CGL equation 
\eqref{eq:lin-CGL} 
with $\mu_0 \le  -\omega_0$ is
\begin{equation}\label{eq:Ap-basecase}
A_p(x,\mu)
=  \left( \sqrt{2 \pi} g(x,\mu+i\omega_0)
+\mathcal{O}(\sqrt{\eps})\right)
e^{\frac{1}{2\eps} (\mu + i \omega_0)^2},
\qquad {\rm for} \ \ \mu \in C_4.
\end{equation}
Therefore, one finds that
$\vert A_p(x,\mu)\vert = 1$ to leading order
along the curve $\mu_{\rm stbc}(x)$
defined by
\begin{equation}\label{eq:stbc}
\left\{ (\mu_{\rm stbc}(x),x) \vert  \quad
{\rm Re}\left( \ln(\sqrt{2\pi}g(x,\mu_{\rm stbc}(x)+i\omega_0)) 
   + \frac{1}{2\eps}(\mu_{\rm stbc}(x) + i \omega_0)^2 \right)
= 0
\right\} .
\end{equation}
This is the space-time buffer curve,
where the real part
of the argument 
of the space-time-dependent phase
in the exponential function
in \eqref{eq:Ap-basecase} vanishes,
to leading order.
It marks the transition
between $A_p$ being exponentially small 
for $\mu \in (\mu_0 , \mu_{\rm stbc}(x))$,
to it being exponentially large
for $\mu > \mu_{\rm stbc}(x)$.
Moreover, to leading order,
the implicit form of the analytical formula is
\begin{equation}\label{eq:stbc-2}
(\mu_{\rm stbc}(x))^2 
= \omega_0^2 -\eps \ln (2\pi) - 2\eps \ln \left| g(x,{\mu_{\rm stbc}}(x)+i\omega_0) \right|.
\end{equation}
In summary,
formula \eqref{eq:stbc}
defines the space-time buffer curve,
and \eqref{eq:stbc-2}
gives the leading order asymptotics
for solutions of \eqref{eq:lin-CGL}
with any $\mu_0 \le -\omega_0$
and the source terms $I_a(x)$ considered here.

\bigskip
\noindent
{\bf Remark.}
In the limit 
$\vert d \vert \to 0$,
the PDE \eqref{eq:gen-CGL} reduces
to a one-parameter family of ODEs in time, 
in which $x$ is the parameter
through $I_a(x)$.
Here, we briefly show that, in this limit,
the space-time buffer curve 
of the PDE \eqref{eq:stbc}
reduces to the buffer point of the $x$-dependent ODE.
In fact, in this limit,
the fundamental solution
of the heat equation 
approaches a delta function,
and $A_h(x,\mu) 
\to 
A_0(x) e^{\frac{1}{2\eps}(\mu+i\omega_0)^2 - (\mu_0+i\omega_0)^2}.$
Then, \eqref{eq:Bp-def} and \eqref{eq:Ap-basecase}
imply $g(x,\mu-{\tilde{\mu}}) \to I_a(x)$,
for all $\mathcal{O}(1)$ values 
$\mu-{\tilde{\mu}}$,
and
$A_p(x,\mu) 
= \left( \sqrt{2\pi} I_a(x) + \mathcal{O}(\sqrt{\eps}) \right)
e^{\frac{1}{2\eps} (\mu+ i \omega_0)^2},
\ {\rm as} \ \vert d \vert \to 0.$
Hence, at each $x$,
$A_p(x,\mu)$
is the same 
as the solution of the corresponding $x$-dependent Shishkova ODE.
Therefore,
at each point,
the space-time buffer curve reduces
to $\omega_0$, which is
the buffer point of the $x$-dependent ODE.
See \cite{HKSW2016,N1987,S1973}
for the general
theory of DHB and buffer points in analytic ODEs.

%------------------------------------------------------------------------------
\subsection{Examples of the space-time buffer curve, $\mu_{\rm stbc}(x)$}
\label{sec:stbc-2}
%------------------------------------------------------------------------------

To study the space-time buffer curve \eqref{eq:stbc},
we give three examples,
involving different types of source terms:
uni-modal (Gaussian), 
spatially-periodic, 
and smoothed step function.
The first example
consists of 
the PDE \eqref{eq:gen-CGL}
with a Gaussian source term, 
\begin{equation}\label{eq:Gaussian}
I_G(x) = e^{\frac{-x^2}{4\sigma}}, \quad \sigma>0.
\end{equation}
Gaussian source terms
are simple models for spatially-localized sources
in R-D equations,
such as the amplitude of a light-source
in chemical pattern formation,
the refractive index of waveguides in nonlinear optics,
and spatially-localized electrical currents
applied to arrays of nerve cells.

\begin{figure}[h]
   \centering
   \includegraphics[width=5in]{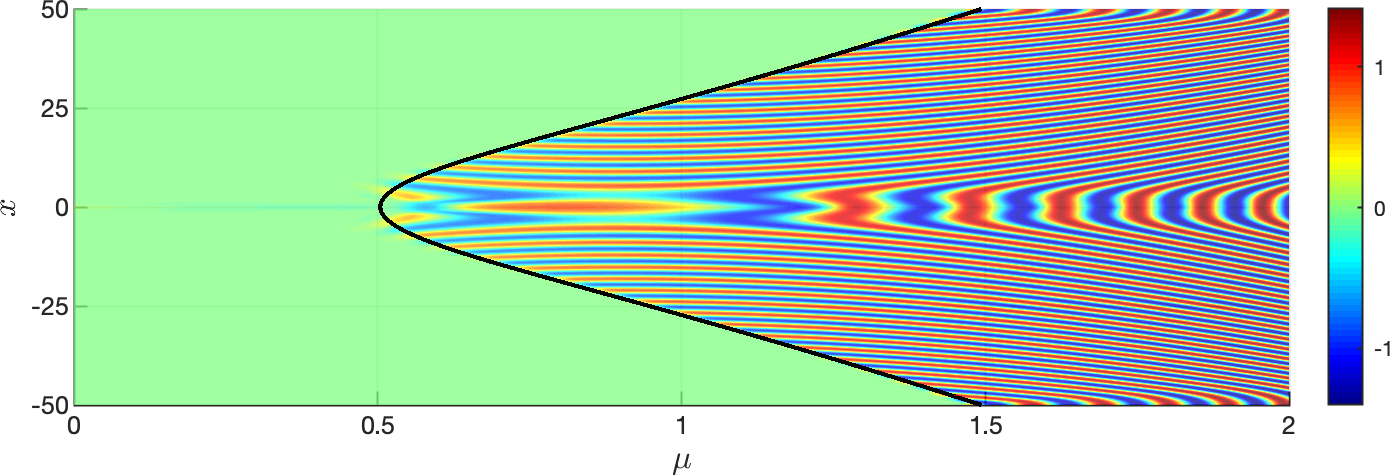}
   \caption{
${\rm Re}(A(x,\mu))$ 
obtained from 
\eqref{eq:gen-CGL} 
with source term
$I_G(x)$.
In the green region,
the solution 
lies close to the repelling QSS.
Superimposed is the black space-time buffer curve 
obtained by solving \eqref{eq:stbc}
with $g$ given by \eqref{eq:g-Gaussian}.
The oscillations commence
just before the space-time buffer curve,
where they have small amplitude.
Beyond it, the oscillations have large amplitude,
since the amplitude of the stable limit cycles is large by this time.
The maxima and minima of the oscillations propagate 
to the center of the interval,
toward $x=0$, where they disappear.
(The phase velocities and group velocities
of the periodic waves are calculated in Section~\ref{sec:Post-DHB},
where the nature of the defect is studied.)
The instantaneous Hopf bifurcation occurs at $\mu=0$ to leading order
at all $x$.
Here, $\eps=0.01, \omega_0=\frac{1}{2}, d_R = 3, d_I = 1$,
$\alpha=0$, and $\sigma=\frac{1}{4}$.
The initial data at $\mu_0= -1.5$ 
is $A_0(x) = -\sqrt{\eps} \frac{I_G(x)}{\mu_0+i \omega_0}$. 
Examples with other initial data are presented below.
}
   \label{fig:Gaussian-stbc}
\end{figure}

For $\mu_0 \le {\tilde{\mu}} < \mu$,
we find from \eqref{eq:Bp-def}(b) that
\begin{equation} \label{eq:g-Gaussian}
g(x,\mu-\tilde{\mu})
=  \sqrt{ \frac{\sigma}{d(\mu- \tilde{\mu}) + \sigma}  }
   \ \ 
    e^{\frac{-x^2}{4(d(\mu-\tilde{\mu}) + \sigma)}}.
\end{equation}
The function $g$
is analytic 
along the contour $C$ 
and in a neighbourhood of it,
except at the branch point 
and along the cut.
Application 
of the general formula 
\eqref{eq:Bp-stationaryphase}
for $\mu$ on $C_4$
yields
\begin{equation}
B_p (x,\mu) 
=  \sqrt{\frac{2 \pi \sigma}{d(\mu + i \omega_0) + \sigma}}
  \ {\rm exp} \left[ \frac{-x^2}{4(d(\mu + i \omega_0) + \sigma)} \right]
+\mathcal{O}(\sqrt{\eps}).
\nonumber
\end{equation}
Translating back using \eqref{eq:AB},
we find 
for all $\mu$ on $C_4$,
\begin{equation}\label{eq:Ap-Gaussian}
A_p(x,\mu)
=  \left( \sqrt{\frac{2 \pi \sigma}{d(\mu + i \omega_0) + \sigma}}
    \ 
  {\rm exp} \left[ \frac{-x^2}{4(d(\mu + i \omega_0) + \sigma)} \right]
+\mathcal{O}(\sqrt{\eps})\right)
e^{\frac{1}{2\eps} (\mu + i \omega_0)^2}.
\end{equation}
Therefore, 
by condition \eqref{eq:stbc},
$\mu_{\rm stbc}(x)$ is given implicitly
to leading order by
\begin{equation}\label{eq:stbc-Gaussian}
\mu^2 = \omega_0^2 
+ \frac{\eps x^2 (d_R \mu - d_I \omega_0 + \sigma)}
       {2\left( (d_R \mu - d_I \omega_0 + \sigma)^2 
                + (d_R\omega_0 + d_I\mu)^2 \right) }
-\eps \ln (2\pi \sigma)
+\frac{\eps}{2} \ln( (d_R\mu - d_I\omega_0 + \sigma)^2 
          + (d_R\omega_0 + d_I\mu)^2),
\end{equation}
provided
\[
{\rm Re} \left(
\frac{d(\mu+i\omega_0) + \sigma}{d(\mu+i\omega_0)}
\right) \ge 0.
\]
Condition \eqref{eq:stbc}
implicitly defines
the space-time buffer curve 
$(\mu_{\rm stbc}(x),x)$ 
along which $\vert A_p \vert = 1$
for Gaussian sources in \eqref{eq:gen-CGL}.

Figure~\ref{fig:Gaussian-stbc}
illustrates this result.
For $0<\mu<\mu_{\rm stbc}(x)$,
the solution is exponentially close 
to the repelling QSS.
Then, at each $x$,
the solution diverges from the repelling QSS there,
and the post-DHB
oscillations set in,
as soon as 
$\mu$ is $\mu_{\rm stbc}(x)$ to leading order
for that $x$.
Overall, 
\eqref{eq:stbc-Gaussian}
and the results presented in Figure~\ref{fig:Gaussian-stbc}
show that $x=0$ is the minimum of $\mu_{\rm stbc}(x)$,
and the solution first diverges from
the repelling QSS there.
Then, as $\vert x \vert$ increases,
$\mu_{\rm stbc}(x)$ increases, quadratically near the tip.
Hence, the duration of the DHB
({\it i.e.}, the time when the solution
leaves a neighborhood of the repelling QSS)
grows quadratically with $\vert x\vert$
near the tip. 

\begin{figure}[h]
   \centering
   \includegraphics[width=5in]{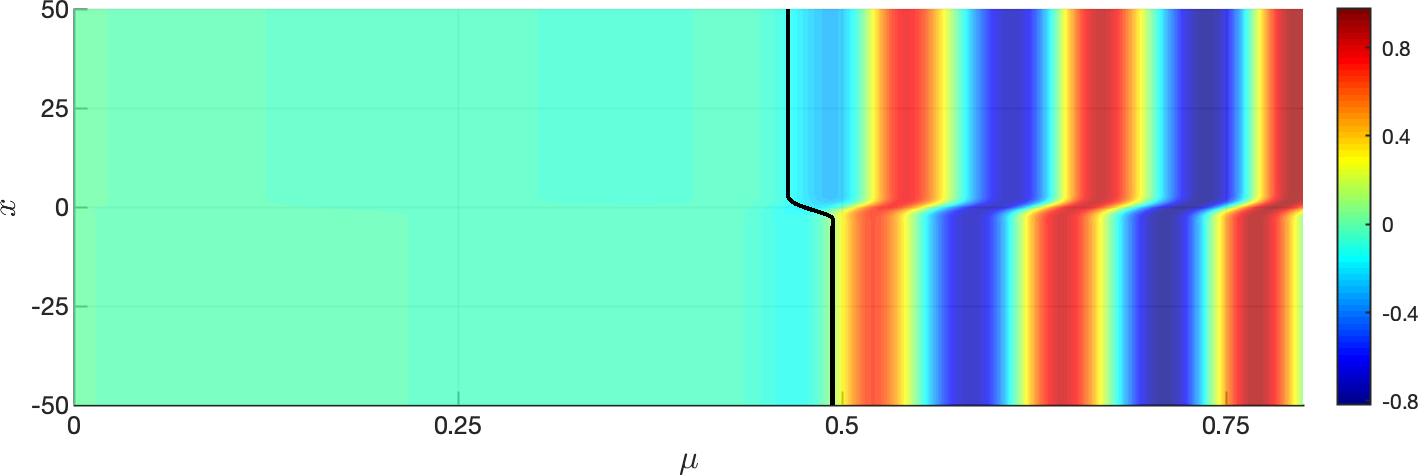}
   \caption{
${\rm Re}(A(x,\mu))$ obtained from
\eqref{eq:gen-CGL}
with the error function source term
$I_{\rm erf}(x)$, 
$I_{\rm ave} = 0.5$,
and $I_{\rm e} = 0.125$.
The black space-time buffer curve 
is super-imposed,
showing that it gives, to leading order,
the time of onset of the oscillations 
at all points $x$ 
in the domain.
Here, the initial data at $\mu_0= -1$ is 
$A_0(x) = -\sqrt{\eps} \frac{I_{\rm erf}(x)}{\mu_0+i \omega_0}$, and
the parameter values are
$\eps=0.01,$
$\omega_0=0.5,$
$\alpha=0$,
$d_R=1$,
and $d_I=0$.
}
   \label{fig:erf-stbc}
\end{figure}

The second example of the space-time buffer curve 
is given by a smoothed step function, 
\begin{equation}\label{eq:erf}
I_{\rm erf}(x) = I_{\rm ave} + I_{\rm e} {\rm erf}(x),
\end{equation}
with $I_{\rm ave} > I_{\rm e} > 0$.
(The error function is
${\rm erf}(x) = \frac{2}{\sqrt{\pi}} \int_0^x e^{-t^2} dt$,
and ${\rm erf}(-x) = - {\rm erf}(x)$.)
This is a simple form
for R-D systems 
in which there is (approximately) piecewise constant input,
with some portion of the domain
receiving one level 
($I_{\rm ave} + I_{\rm e}$)
and a complementary part
receiving a different level 
($I_{\rm ave} - I_{\rm e}$),
with a smooth transition in between.
By \eqref{eq:Bp-def}(b), 
one finds
\begin{equation}\label{eq:g-erf}
g(x,\mu-{\tilde{\mu}})
= I_{\rm ave} + I_{\rm e} {\rm erf} \left(
      \frac{x}{\sqrt{1+4d(\mu-{\tilde{\mu}})}} \right),
\end{equation}
see for example \cite{NG1969}, 
provided that the argument of the error function
lies in $(-\frac{\pi}{4},\frac{\pi}{4})$.
Hence,
by \eqref{eq:stbc},
the space-time buffer curve for \eqref{eq:gen-CGL}
with $I_{\rm erf}(x)$
is to leading order
\begin{equation}
\label{eq:stbc-erf}
(\mu_{\rm stbc}(x))^2 = \omega_0^2 
- 2 \eps \ln \left( I_{\rm ave} + I_{\rm e} {\rm erf}
\left( \frac{x}{\sqrt{1 + 4d(\mu_{\rm stbc}(x)+i\omega_0)}} \right) \right)
-\eps \ln(2\pi).
\end{equation}
See Figure~\ref{fig:erf-stbc}.
Small-amplitude oscillations set in
just before $\mu_{\rm stbc}(x)$.
Then, at each point $x$,
the amplitude of the oscillations becomes large
as soon as $\mu$ reaches $\mu_{\rm stbc}(x)$,
to leading order.

\begin{figure}[h]
   \centering
   \includegraphics[width=5in]{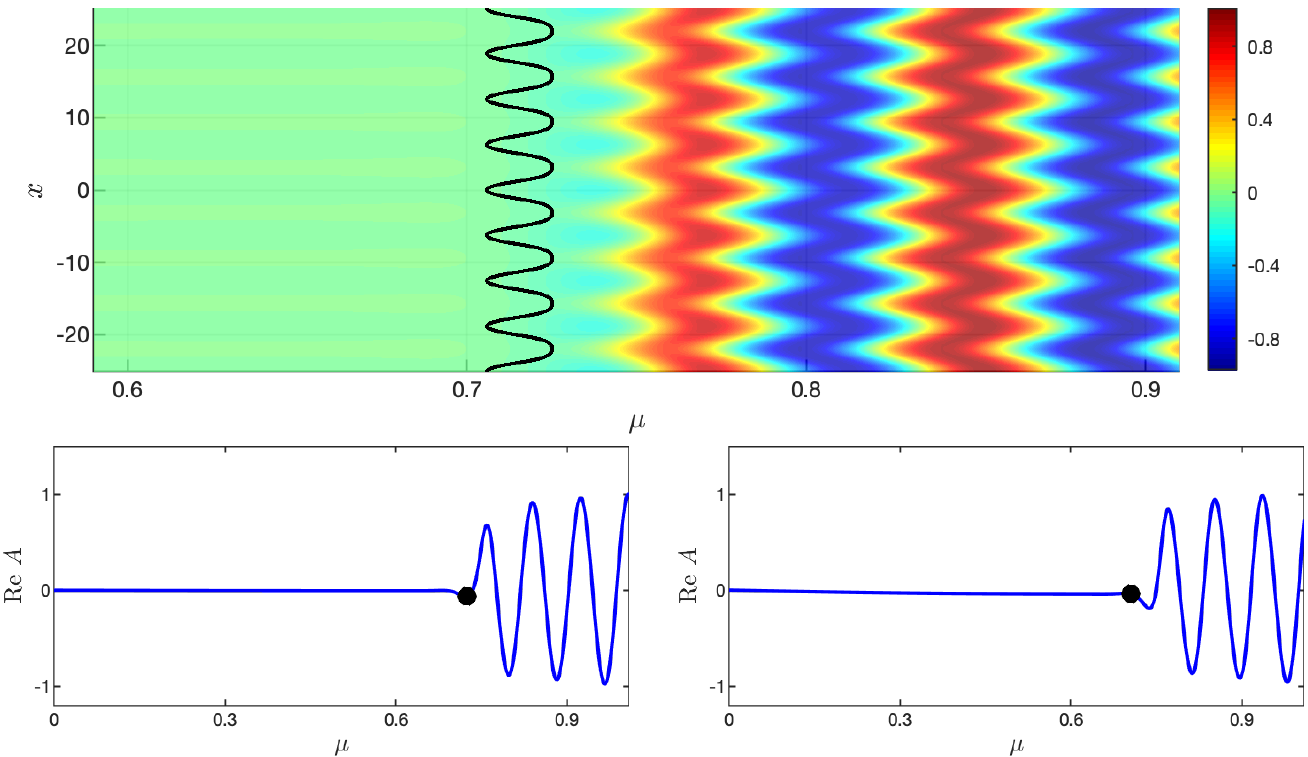}
   \caption{
${\rm Re}(A)$ for the solution of \eqref{eq:gen-CGL}
with the spatially-periodic source term,
$I_{\rm per}(x)$,
$p_1 = \frac{1}{3}$,
and $p_2 = \frac{1}{4}$.
The space-time buffer curve obtained by setting
$\vert A_p(x,\mu) \vert =\sqrt{\eps}$
is super-imposed (black curve).
Time series at $x=5\pi$ (a local minimum
of $I_{\rm per}(x)$)
and $x=6\pi$ (a local maximum)
are also shown.
The initial data at $\mu_0= -1$
is $A_0(x) = -\sqrt{\eps} \frac{I_{\rm per}(x)}{\mu_0+i \omega_0}$.
Here, $\eps=0.01$,
$\omega_0=0.75$,
$\alpha=0$,
$d_R=1$,
and $d_I=0$.
}
   \label{fig:periodic}
\end{figure}

The third example 
of the space-time buffer curve
consists of a spatially periodic source term,
\begin{equation}\label{eq:per}
I_{\rm per}(x) = p_1 + p_2 \cos(x),
\end{equation}
with $p_1 > p_2 > 0$,
$\mathcal{O}(1)$ independent of $\eps$.
Spatially-periodic forcing
arises in various pattern formation problems,
see for instance \cite{DBZE2001,HHM2015}.
From the general definition 
\eqref{eq:Bp-def}(b)
of $g$,
one finds
\begin{equation}\label{eq:g-per}
g(x,\mu-{\tilde{\mu}})
= \frac{1}{\sqrt{4\pi d (\mu - {\tilde{\mu}})}}
\int_R e^{- \frac{(x-y)^2}{4d (\mu - {\tilde{\mu}})}} 
(p_1 + p_2 \cos(y)) dy
= p_1 + p_2 e^{- d(\mu - {\tilde{\mu}})}\cos(x).
\end{equation}
Now, with this elementary form
of $g$,
the integral 
\eqref{eq:Bp-def}(a) for $B_p(x,\mu)$
may be evaluated
in closed form in this example. 
Specifically, recalling
\eqref{eq:AB},
we find
that the particular solution is 
\begin{equation}\label{eq:Ap-per}
\begin{split}
A_p(x,\mu)
&=\sqrt{\frac{\pi}{2}} p_1
\left[ {\rm erf} \left( \frac{\mu+i\omega_0}{\sqrt{2\eps}} \right)
-{\rm erf} \left( \frac{\mu_0+i\omega_0}{\sqrt{2\eps}} \right) \right]
e^{\frac{(\mu+i\omega_0)^2}{2\eps}} \\
&+\sqrt{\frac{\pi}{2}} p_2 \cos(x) 
\left[ {\rm erf} \left( \frac{\mu+i\omega_0 -\eps d}{\sqrt{2\eps}} \right)
-{\rm erf} \left( \frac{\mu_0+i\omega_0 - \eps d}{\sqrt{2\eps}} \right) \right]
e^{\frac{(\mu+i\omega_0-\eps d)^2}{2\eps}}.
\end{split}
\end{equation}
The space-time buffer curve,
derived from this exact solution,
is shown in Figure~\ref{fig:periodic}.
There is good agreement 
with the onset of oscillations
observed numerically 
in the cubic CGL \eqref{eq:gen-CGL}.

The space-time buffer curve in this case is determined by the
condition $|A_p(x,\mu)| = \sqrt{\eps}$ instead of the usual
$|A_p(x,\mu)| = 1$ criterion. This choice was made because whilst the
nonlinear terms in the cubic CGL equation do not affect the onset of
the oscillations, they do influence the spatial phase of the
oscillations. In this case of a periodic source term, the cubic
nonlinearities induce a phase shift in space causing the buffer curve
derived from $|A_p(x,\mu)| = 1$ to be $\pi$ units out of phase in the
$x$-direction with the numerically observed onset. By setting the
space-time buffer criterion to be $|A_p(x,\mu)| = \sqrt{\eps}$, the
contribution from the nonlinearities to the phase shift is still
small, and hence there is better agreement between the space-time
buffer curve prediction and the onset
of large-amplitude oscillations 
in the numerically-calculated solutions of \eqref{eq:gen-CGL}.

\bigskip
\noindent
{\bf Remark.}
For this example with a periodic source term,
we have derived the space-time buffer curve above
from the exact, closed form
expression for the particular solution
$A_p$, \eqref{eq:Ap-per}.
One may also find the leading order asymptotics
using \eqref{eq:stbc}.
In fact,
recalling \eqref{eq:Ap-basecase},
we see that
the leading order term in the particular solution 
for each $\mu$ on $C_4$ is
$A_p(x,\mu)
= (\sqrt{2\pi} g(x,\mu+i\omega_0) + \mathcal{O}(\sqrt{\eps}))
  e^{\frac{1}{2\eps} (\mu + i \omega_0)^2},$
with $g$ given by \eqref{eq:g-per}.
Hence, to leading order,
the space-time buffer curve 
is given implicitly by
$(\mu_{\rm stbc}(x))^2 = \omega_0^2 - \eps \ln(2\pi) 
+ 2\eps \ln \left| p_1+p_2 e^{-d(\mu_{\rm stbc}(x)+i\omega_0)}\cos(x) \right|.$

%------------------------------------------------------------------------------
\section{The homogeneous exit time curve, $\mu_h(x)$}
\label{sec:Ah-transition}
%------------------------------------------------------------------------------

In this section,
we study the homogeneous 
component, $A_h(x,\mu)$, of the solution
of the linear CGL PDE \eqref{eq:lin-CGL}.
The main result 
is the homogeneous exit time curve,
along which
$\vert A_h(x,\mu)\vert = 1$,
{\it i.e.,} where $A_h$ transitions
from being exponentially small to large.
We label it as $\mu_{h}(x)$.

Recalling the change of variables \eqref{eq:AB},
we see that
formula \eqref{eq:B_h} gives
\begin{equation}
\label{eq:Ah-basecase}
A_h(x,\mu) 
= \frac{ e^{\frac{1}{2\eps}
              \left( (\mu+i\omega_0)^2-(\mu_0+i \omega_0)^2 \right) } }
       {\sqrt{4\pi d (\mu - \mu_0)}}
\int_R  e^{\frac{-(x-y)^2}{4d(\mu-\mu_0)}} A_0(y) dy,
\qquad {\rm for} \ \ \mu>\mu_0.
\end{equation}
The integral in \eqref{eq:Ah-basecase} can be evaluated 
for many different types of bounded initial data $A_0(x)$.
Moreover, the function $A_h(x,\mu)$
is bounded for $\mu \in (\mu_0,-\mu_0)$
and,
for our examples,
also analytic in a region of the complex plane about this interval,
excluding the branch point and cut.

\begin{figure}[h]
   \centering
   \includegraphics[width=5in]{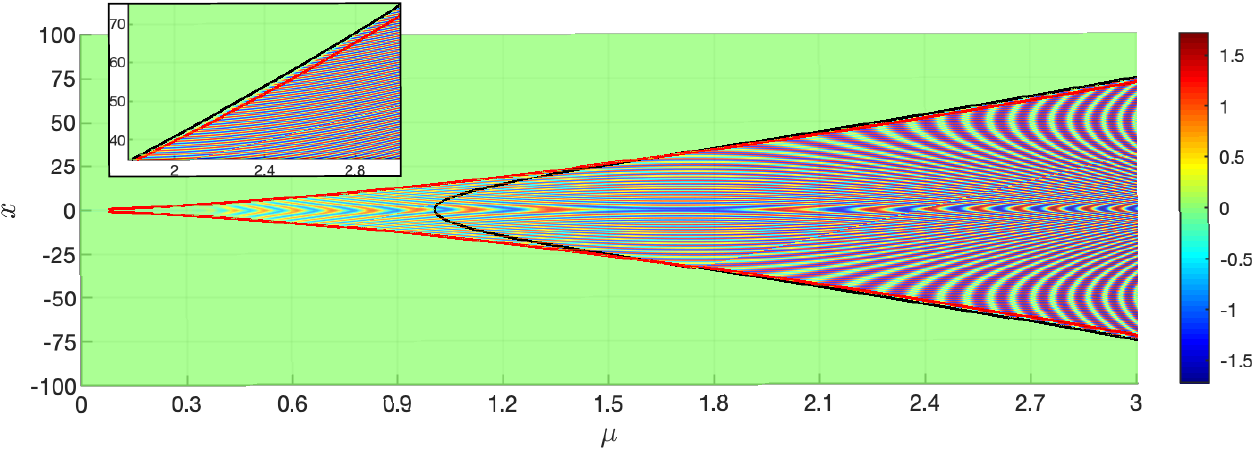}
   \caption{
${\rm Re}(A)$ 
of the solution
obtained from the direct numerical simulation
of \eqref{eq:gen-CGL}
with $\omega_0=1$
and source term 
$I_G(x)=e^{-\frac{x^2}{4\sigma}}$, $\sigma=\frac{1}{4}$.
The initial data
$A_0(x) = -\sqrt{\eps} \frac{I_G(x)}{\mu_0+i \omega_0}$ 
is given at $\mu_0=-0.2$.
For each $x$,
the solution stays near the repelling QSS 
(green region)
for a long time
past the instantaneous Hopf bifurcation at $\mu=0$.
Then, the exit time 
from a neighbourhood of the repelling QSS 
is $x$-dependent.
In the center,
it is given by 
\eqref{eq:hom-exittime} (red curve),
where $\mu_{h}(x)
< \mu_{\rm stbc}(x)$.
In contrast, for $x$ outside this interval,
$\mu_{\rm stbc}(x)$ comes first,
and the space-time buffer curve 
\eqref{eq:stbc}
(black curve)
determines the exit time.
The inset shows a magnification for $\mu > 1.8$.
The convex hull given by
$(x, \mu = {\rm min}(\mu_h(x),\mu_{\rm stbc}(x)))$
agrees to leading order with the onset of oscillations 
for all $x$.
Here, $\eps=0.01, 
\omega_0=1,
\alpha=0.6,
d_R = 1, d_I = 0.$
}
   \label{fig:exittime-determined-by-Ah}
\end{figure}

To give a first illustration,
we choose 
$A_0(x) = -\frac{\sqrt{\eps}I_a(x)}{\mu_0+i\omega_0}$,
which is the leading order term 
in the attracting QSS,
and we use the Gaussian source $I_G(x)=e^{-\frac{x^2}{4\sigma}}$.
(Examples with more general initial data
and with other source terms will be given in Section~\ref{sec:DHB-fourcases}.)
The integral in \eqref{eq:Ah-basecase} yields
\begin{equation} \label{eq:Ah-G}
A_h(x,\mu)
= - \left( \frac{\sqrt{\eps}}{\mu_0+i\omega_0} \right)
e^{\frac{1}{2\eps}\left[ (\mu+i\omega_0)^2 - (\mu_0+i\omega_0)^2 \right]}
\sqrt{ \frac{\sigma}{d(\mu-\mu_0) + \sigma}}
e^{-\frac{x^2}{4(d(\mu-\mu_0) + \sigma)}}.
\end{equation}
Hence, the argument of the total exponential in $A_h$
depends on both $x$ and $\mu$.
Setting $\vert A_h(x,\mu)\vert = 1$,
we find that for $\mu>0$
the exit time $\mu_h(x)$ 
is given implicitly by
\begin{equation}
\label{eq:hom-exittime}
(\mu_h(x))^2
= \mu_0^2 
+ \frac{ \eps x^2 \left[ d_R(\mu_h(x)-\mu_0) + \sigma \right]}
   {2\left[ (d_R(\mu_h(x)-\mu_0)+ \sigma)^2 + d_I^2(\mu_h(x)-\mu_0)^2\right]},
\end{equation}
where to leading order 
$\mu_h(x) = - \mu_0$.
This curve is the homogeneous exit time curve,
$\mu_{h}(x)$.
The logarithmic terms at $\mathcal{O}(\eps)$
are not reported here, but
may be calculated as in Section~\ref{sec:stbc-2}.

Figure~\ref{fig:exittime-determined-by-Ah}
reveals the role
played by the (red) curve 
\eqref{eq:hom-exittime}
in determining the exit time 
(and the time of onset of
oscillations) for solutions
of the cubic CGL with Gaussian source term.
There is a central interval about $x=0$
on which 
$\mu_{h}(x)
< \mu_{\rm stbc}(x)$.
On this interval,
the homogeneous solution $A_h$
with the given initial data
switches from being exponentially small 
to being exponentially large
before the particular solution $A_p$ does so for this source.
Hence, it
determines the exit time 
(and the time of onset of the oscillations) there.
See the red curve 
in Figure~\ref{fig:exittime-determined-by-Ah}.
In contrast,
for $x$ outside this interval,
the situation is reversed,
with $\mu_{\rm stbc}(x)$ occurring first.
Hence, at all points $x$ outside this interval,
the exit time (and onset time for the oscillations)
is determined
by the space-time buffer curve 
$\mu_{\rm stbc}(x)$ given by \eqref{eq:stbc}.
See the black curve 
in Figure~\ref{fig:exittime-determined-by-Ah}
and the inset.

%------------------------------------------------------------------------------
\section{The main cases of DHB: 
one case for each different type of outcome in the competition between 
which of $A_p$ and $A_h$ ceases to be exponentially small first}
\label{sec:DHB-fourcases}
%------------------------------------------------------------------------------

From the analysis 
in Sections~\ref{sec:stbc} 
and \ref{sec:Ah-transition},
we see that at each point $x$ 
there is a competition
between which of $\mu_{\rm stbc}(x)$
and $\mu_{h}(x)$ comes first,
{\it i.e.,} between which component, 
$A_p(x,\mu)$ or $A_h(x,\mu)$,
transitions first from being exponentially small to exponentially large.
Moreover, the formulas
\eqref{eq:stbc-2}
and \eqref{eq:hom-exittime} for
$\mu_{\rm stbc}(x)$
and $\mu_{h}(x)$
show that these times
depend on key parameters,
$\omega_0$ and $\eps$,
the initial data $A_0(x)$ and time $\mu_0$,
as well as on the form
of $I_a(x)$.

In this section,
we analyze both of these functions
and determine various possible outcomes of the competition.
Each different type of outcome leads
to a distinct type of delayed Hopf bifurcation (DHB).
We begin 
in Subsection~\ref{sec:DHB-cases1-3}
with cases of DHB
that arise for solutions of \eqref{eq:gen-CGL}
with initial data 
given at any $\mu_0 \le -\omega_0$.
Then, 
in Subsection~\ref{sec:DHB-case4},
we present a main case of DHB
that arises for solutions of \eqref{eq:gen-CGL}
with initial data 
given at any $-\omega_0 < \mu_0 \le -\delta$,
where $\delta>0$ is a small, $\mathcal{O}(1)$ constant.
Also, we illustrate all of these cases of DHB
using the different types of source terms introduced 
in Section~\ref{sec:stbc}:
Gaussian, spatially-periodic, and smoothed step function.

%------------------------------------------------------------------------------
\subsection{Cases 1-3 of DHB 
for solutions with initial data given at $\mu_0 \le  -\omega_0$}
\label{sec:DHB-cases1-3}
%------------------------------------------------------------------------------

For solutions of \eqref{eq:gen-CGL} 
with initial data given at $\mu_0 \le -\omega_0$,
the competition
between $A_p$ and $A_h$ can have three
possible outcomes
depending on which ceases to be exponentially small first.
These correspond
to the following three cases of DHB:

\bigskip
\noindent
{\bf Case 1 of DHB. 
$\mu_{\rm stbc}(x) < \mu_h(x)$
for all $x \in [-\ell,\ell]$.}
In this case, the parameters $\omega_0$ and $\eps$,
the initial data $A_0(x)$ and time $\mu_0$, and
the source term $I_a(x)$
are such that
$A_p(x,\mu)$ ceases to be exponentially small first,
before $A_h(x,\mu)$ does,
for all $x$,
{\it i.e.,}
$A_h(x,\mu_{\rm stbc}(x))$ is exponentially small
at all points $x$.
Hence, the full solution $A$ is
exponentially close to the repelling QSS
until $\mu=\mu_{\rm stbc}(x)$ to leading order,
and the duration of the DHB and the time of onset of the oscillations,
$\mu_{\rm stbc}(x)$,
is determined completely by $A_p$ on the entire domain.

Case 1 of DHB is illustrated
in Figures~\ref{fig:Gaussian-stbc}, 
\ref{fig:erf-stbc}, 
and \ref{fig:periodic}, 
for the Gaussian, spatially-periodic,
and smoothed step function source terms, respectively.
For the Gaussian source term (with the Gaussian initial data),
one finds that $\mu_{\rm stbc}(x) < \mu_h(x)$ 
for all $x \in [-\ell,\ell]$. 
This is consistent with formulas
\eqref{eq:Ap-Gaussian} and \eqref{eq:stbc-Gaussian}
derived above for $\mu_{\rm stbc}(x)$
and with formula \eqref{eq:hom-exittime} for $\mu_h(x)$.
See Figure~\ref{fig:Gaussian-stbc}. 

Next, for the error function source term,
$I_{\rm erf}(x)$
and the initial data used above,
one finds 
\begin{equation}
A_h(x,\mu)
= \frac{-\sqrt{\eps}}{(\mu_0 + i \omega_0)}
e^{-\frac{1}{2\eps}(\mu_0+i\omega_0)^2}
\left[ I_{\rm ave} 
+ I_{\rm e} {\rm erf}\left( \frac{x}{\sqrt{1+4d(\mu-\mu_0)}} \right) \right].
\end{equation}
Now, the homogeneous exit time curve
$\mu_h(x)$ is  
obtained directly by setting $\vert A_h(x,\mu)\vert=1$
with the exact solution. 
Hence, recalling that $\mu_{\rm stbc}(x)$
is given by setting
$\vert A_p(x,\mu)\vert = 1$
with $g$ given exactly by \eqref{eq:g-erf},
we see that $\mu_{\rm stbc}(x) < \mu_h(x)$
for all $x \in [-\ell,\ell]$.
There is
again good quantitative agreement
between the leading order space-time buffer curve
and the numerically observed onset
of the large-amplitude oscillations,
as shown
in Figure~\ref{fig:erf-stbc}.

The third example of Case 1 of DHB
is given by the simulation 
with the spatially-periodic source term,
$I_{\rm per}(x)$,
and initial data 
$A_0(x)=-\sqrt{\eps} \frac{I_{\rm per}(x)}{\mu_0+i\omega_0}$.
The homogeneous solution is
$A_h(x,\mu)= \sqrt{2} e^{-d(\mu-\mu_0)}
{\rm exp}\left[
\frac{1}{2\eps}\left((\mu+i\omega_0)^2 - (\mu_0+i\omega_0)^2\right)
\right]
\cos(x).$ 
Hence,
the homogeneous exit time curve 
is $\mu_h(x) = -\mu_0$
to leading order,
which is derived from this exact solution.
There are also $\mathcal{O}(\eps)$ corrections,
which are periodic in space.
Then, from \eqref{eq:Ap-per},
we see that $\mu_{\rm stbc}(x) < \mu_h(x)$
for all points on the domain.
See Figure~\ref{fig:periodic}.

\bigskip
\noindent
{\bf Case 2 of DHB. 
$\mu_{\rm stbc}(x) < \mu_h(x)$ 
for some intervals of points $x$ on $[-\ell,\ell]$,
and $\mu_h(x) < \mu_{\rm stbc}(x)$ on the complementary intervals.}
This case arises
when the parameters $\omega_0$ and $\eps$,
the initial data $A_0(x)$ and time $\mu_0$, and
the source term $I_a(x)$
are such that
$\mu_{\rm stbc}(x) < \mu_h(x)$
at some, but not all, points $x$,
and $\mu_h(x) < \mu_{\rm stbc}(x)$ 
on the complementary intervals,
even though $\omega_0 < -\mu_0$.
Here, $A_p$
first causes the solution
to diverge from the repelling QSS
at points $x$ where
$\mu_{\rm stbc}(x) < \mu_h(x)$,
before the homogeneous component $A_h$
can do so.
On the complementary intervals,
where $\mu_h(x) < \mu_{\rm stbc}(x)$,
the homogeneous solution 
stops being exponentially small first,
and hence $A_h(x,\mu)$
determines the DHB.

\begin{figure}[h]
   \centering
\includegraphics[width=5in]{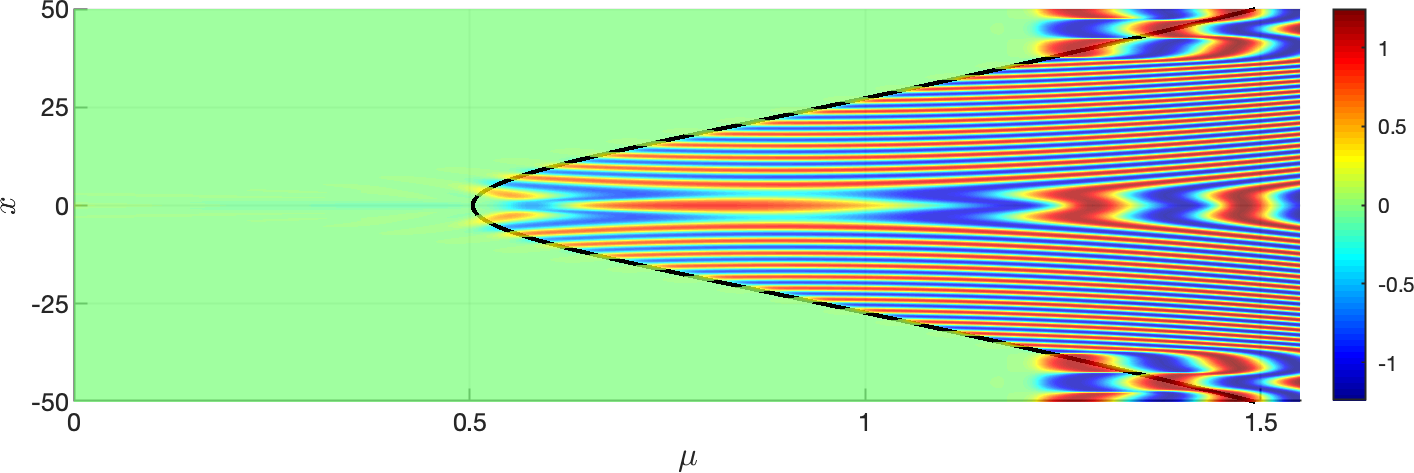}
   \caption{
Illustration of DHB in Case 2
showing ${\rm Re} (A(x,\mu))$ 
obtained from \eqref{eq:gen-CGL} 
with $I_G(x)$. 
The initial data is given at time $\mu_0=-1.2$.
For $\vert x \vert \lesssim  37.1$,
$\mu_{\rm stbc}(x) < -\mu_0 = +1.2$.
Hence, the space-time buffer curve 
given by the particular solution determines the exit time
to leading order for these points $x$.
In contrast, for $\vert x \vert \gtrsim 37.1$,
$\mu_{\rm stbc}(x) > -\mu_0 = +1.2$,
so that for these $x$
the exit time is determined to leading order by $A_h(x,\mu)$.
Here, $\eps=0.01, \omega_0=\frac{1}{2}, d_R = 3, d_I = 1$,
$\alpha=0.6$, and 
$A_0(x) = \cos\left( \frac{10\pi x}{\ell}\right)$
with $\ell=50$. 
Similar results
are obtained for other initial data.
}
   \label{fig:gaussian-DHBcase2}
\end{figure}

An example of DHB in Case 2
is presented in Figure~\ref{fig:gaussian-DHBcase2}.
Here, the source term is
$I_G(x)=e^{-\frac{x^2}{4\sigma}}$
with $\sigma=\frac{1}{4}$,
and the initial data at $\mu_0=-1.2$ is 
$A_0(x)=\cos\left(\frac{n \pi x}{\ell}\right)$,
with $n=10$ and $\ell=50$.
For this initial condition,
we find 
$A_h(x,\mu) = \sqrt{2} e^{-\frac{d n^2 \pi^2}{\ell^2}(\mu-\mu_0)}
{\rm exp}\left[\frac{1}{2\eps}\left( (\mu+i\omega_0)^2 
                              - (\mu_0+i\omega_0)^2\right) \right]
\cos\left( \frac{n\pi x}{\ell} \right)$.
Hence, we find
$\mu_h(x)$ directly from this exact solution.

We observe that $\mu_{\rm stbc}(x)= - \mu_0$
at $\vert  x \vert  \approx 37.1.$
Here, to leading order,
the competition is a tie.
For $\vert x \vert \lesssim 37.1$,
$\mu_{\rm stbc}(x) < \mu_h(x)$,
so that the space-time buffer curve 
\eqref{eq:stbc-Gaussian}
predicts when
the solutions diverge 
from the repelling QSS and begin to oscillate.
In contrast,
for all $\vert x \vert \gtrsim 37.1$,
$\mu_h(x) < \mu_{\rm stbc}(x) $,
{\it i.e.,} $A_h$ first ceases to be
exponentially small,
while $A_p$
remains exponentially small.
Hence, 
for $\vert x \vert \gtrsim 37.1$,
the solution diverges from
the repelling QSS 
as $\mu$ reaches $\mu_h(x) \sim -\mu_0$, 
to leading order.
See Figure~\ref{fig:gaussian-DHBcase2}.

On the outer intervals,
the source term $I_G(x)$ is essentially zero
(below round-off in the simulations).
Hence, 
the cubic CGL PDE
is effectively symmetric under 
$A \to Ae^{i\theta}$ ($\theta$ real),
including under $A \to -A$ ($\theta=\pi$) here.
Hence, to leading order,
it has a symmetric way-in way-out function,
so the homogeneous exit time
is $-\mu_0$, to leading order
at these $x$.
Here, the amount of credit built up 
as $\mu$ increases from $-\mu_0$ to zero
and the solution
spirals toward
the attracting QSS
exponentially 
is exactly spent as $\mu$ increases
from zero to $-\mu_0$,
and the solution
spirals away exponentially
from the repelling QSS.
Simulations with other values of $\mu_0$
in $(-1.5,-0.5)$ 
and with other initial data
show similar results for $\mu_{\rm stbc}(x) $ and $\mu_h(x)$,
with the onset being determined
by $\mu_{\rm stbc}(x)$
in the central portion of the domain
and by $\mu_h(x)$
in the outer portions.

\begin{figure}
   \centering
    \includegraphics[width=6.5in]{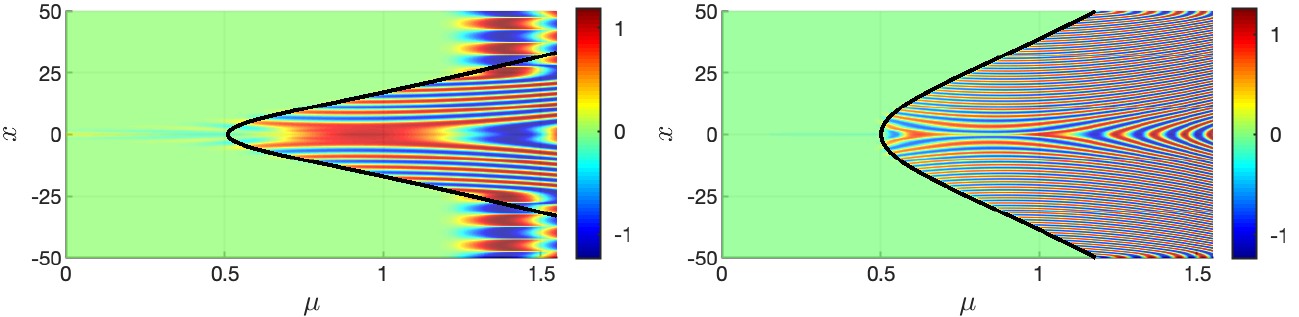}
   \caption{
On a finite domain, the parameter $\eps$ 
changes the width of space-time buffer curve.
Here, we fix the values of all of the parameters (except $\eps$) 
to be the same 
as in Figure~\ref{fig:gaussian-DHBcase2}.
Then, in comparison 
with Figure~\ref{fig:gaussian-DHBcase2} 
where $\eps=0.01$,
we observe that
(a) for $\eps=0.025$ 
the space-time buffer curve is narrower,
since at any point $x$ the duration of the DHB
$\mu_{\rm stbc}(x)$
increases as $\eps$ gets larger, by \eqref{eq:stbc-2}.
Conversely, (b) for $\eps=0.005$,
the space-time buffer curve is wider,
since at any point $x$ the duration of the DHB
$\mu_{\rm stbc}(x)$
decreases as $\eps$ gets smaller. 
In fact, for this smaller value of $\eps$,
$\mu_{\rm stbc}(x) < \mu_h(x)$ for all $x$ on the domain,
and the solution has shifted over to being in Case 1 of DHB
with this domain.
}
   \label{fig:DHB-role-epsilon}
\end{figure}

The parameter $\eps$ plays an important role in determining
the width of $\mu_{\rm stbc}(x)$, and hence whether
a solution on a finite domain exhibits Case 2 or Case 1 of DHB.
This is illustrated in Figure~\ref{fig:DHB-role-epsilon}.

Finally, for DHB in Case 2 with a Gaussian source term,
we observe that
there is a difference between 
the spatio-temporal dynamics
of the large-amplitude oscillations in $A(x,\mu)$
which are observed
in the central portion of the domain
after the space-time buffer curve $\mu_{\rm stbc}(x)$ is crossed
and those which arise 
in the outer portions of the interval $[-\ell,\ell]$,
after the homogeneous exit time curve $\mu_h(x)$ is crossed.
In the central portion 
(where $A_p$ first becomes exponentially large,
which is on $\vert x \vert \lesssim 37.1$ 
in Figure~\ref{fig:gaussian-DHBcase2}),
the large-amplitude oscillations 
propagate spatially,
initially to $x=0$ and then
outward, away from $x=0$
for most $\mu$
($\mu \gtrsim 0.8$ in Figure~\ref{fig:gaussian-DHBcase2}).
In contrast, outside the central portion
(where $A_h$ first becomes exponentially large
{\it i.e.,} where $\mu_h(x) < \mu_{\rm stbc}(x) $),
the oscillations
do not propagate spatially.
Moreover, at the interfaces
({\it e.g.}, at $\vert x \vert \sim 37.1$ in the figure),
the outward propagating pulses get absorbed
by the regime in which the oscillations do not propagate.
The spatio-temporal dynamics 
of the post-DHB oscillations
is discussed briefly in~\ref{sec:Post-DHB}.

\begin{figure}[h]
   \centering
\includegraphics[width=5in]{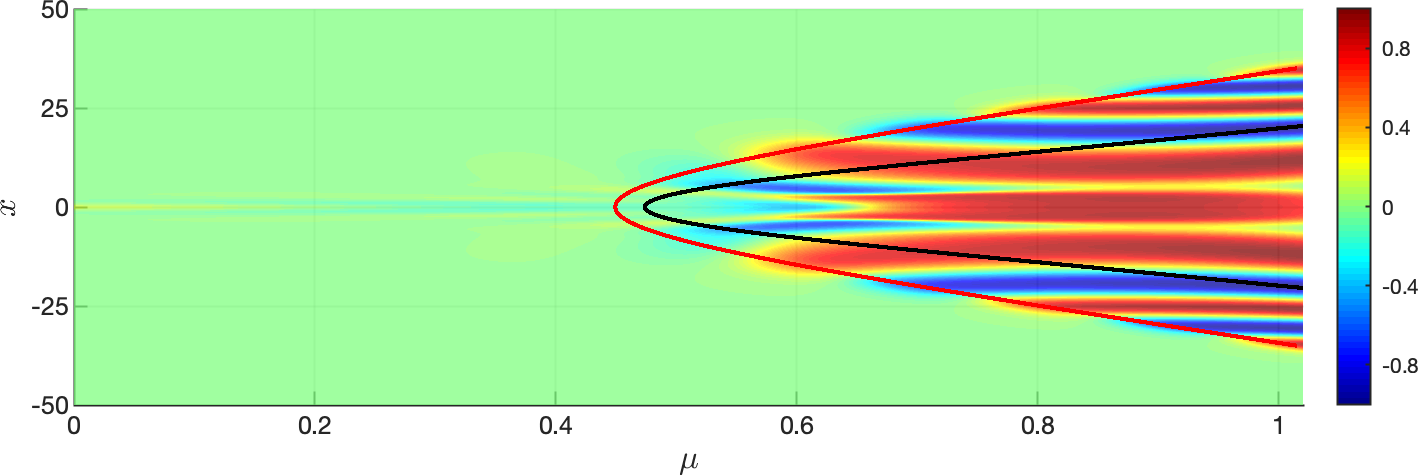}
   \caption{
${\rm Re}(A(x,\mu))$ 
illustrating DHB in Case 3
for \eqref{eq:gen-CGL} 
with $I_G(x)$ ($\sigma=\frac{1}{4}$). 
The curves $\mu_h(x)$ (red) and
$\mu_{\rm stbc}(x)$ (black)
are super-imposed.
For all $x$ in $[-\ell,\ell]$ ($\ell=50$),
$\mu_h(x) < \mu_{\rm stbc}(x)$.
Hence, $A_h$ is the first to stop being exponentially small and grow,
before $A_p$ does, for each $x$.
The oscillations just to the left of the red curve
are small-amplitude, and they are large-amplitude
as soon as $\mu$ reaches $\mu_h(x)$.
The parameters are
$\eps=0.02, 
\omega_0=\frac{1}{2}, 
\alpha=0.6, 
d_R = 3, d_I = 1$.
The initial data,
given at $\mu_0=-0.55$,
is $A_0(x)=-20 e^{-\frac{x^2}{40}}$.}
   \label{fig:case3}
\end{figure}

\bigskip
\noindent
{\bf Case 3 of DHB. 
$\mu_h(x) < \mu_{\rm stbc}(x)$  
for all $x \in [-\ell,\ell]$.}
In this case, the parameters $\omega_0$ and $\eps$,
the initial data $A_0(x)$ and time $\mu_0$, and
the source term $I_a(x)$
are all such that
the homogeneous component $A_h$ 
stops being exponentially small first
at all points $x$. 
It causes the solution $A(x,\mu)$
to diverge from the repelling QSS
at the time $\mu_h(x)$,
since at each point 
$A_p(x,\mu_h(x))$ is exponentially small.
Hence, the DHB 
is determined completely by $A_h$.
An example is given in Figure~\ref{fig:case3}.

%------------------------------------------------------------------------------
\subsection{Case 4 of DHB 
for initial data given at any $-\omega_0 < \mu_0 \le -\delta$}
\label{sec:DHB-case4}
%------------------------------------------------------------------------------

Case 4 of DHB arises for solutions
$A(x,\mu)$ of \eqref{eq:gen-CGL}
with initial data given at $\mu_0 \in (-\omega_0,-\delta]$,
where $\delta>0$ is again any small, $\mathcal{O}(1)$ constant.
With this initial time,
one has $-\mu_0<\omega_0$.
Hence, the left tip of $\mu_h(x)$
(where $\vert A_h \vert=1$,
as calculated from \eqref{eq:Ah-basecase})
comes before the left tip of the space-time buffer curve
$\mu_{\rm stbc}(x)$
(where $\vert A_p \vert=1$,
as calculated from \eqref{eq:stbc}).
In this case,
the source terms $I_a(x)$,
parameters $\omega_0$ and $\eps$,
and initial data $A_0(x)$
are such that, for some intervals of $x$,
the homogeneous component $A_h$ 
stops being exponentially small before $\mu$ reaches $\omega_0$,
{\it i.e.,} before $A_p$ can.
Furthermore, it does so in a manner that $\mu_h(x)$
has non-trivial spatial dependence.
We illustrate this with two examples.

\begin{figure}[h]
   \centering
\includegraphics[width=5in]{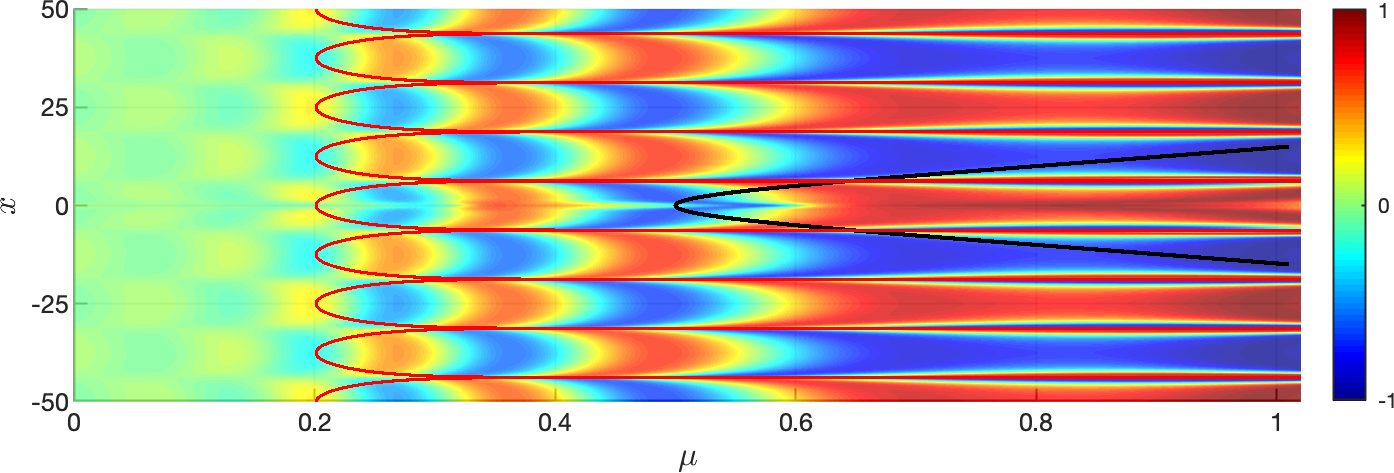}
   \caption{
Illustration of DHB in Case 4
with $I_G(x)$
and initial data $A_0(x)=\cos(n\pi x/\ell)$ 
with $n=4$ and $\ell=50$,
given at $\mu_0=-0.2$.
For all $x$,
the exit time 
for solutions of the CGL PDE \eqref{eq:gen-CGL},
is $\mu_h(x)$ (red curve),
as determined by $A_h$.
We observe that
that $\mu_h(x) \sim -\mu_0=0.2$,
and the $\mathcal{O}(\eps)$ corrections
are spatially periodic.
The temporal oscillations 
which set in after $\mu_h(x)$ are stationary in space;
there are no propagating pulses.
The parameter values 
are $\eps=0.01, \omega_0=\frac{1}{2}, d_R = 1, d_I = 0$,
and $\alpha=0.6$. 
The space-time buffer curve (black curve)
is super-imposed to illustrate that $\mu_{\rm stbc}(x)$
plays no role in the onset of the oscillations.}
   \label{fig:gaussian-DHBcase4}
\end{figure}

The first example of Case 4 of DHB 
is given in Figure~\ref{fig:exittime-determined-by-Ah}.
Here, the solution $A$ of the cubic PDE \eqref{eq:gen-CGL}
stops being exponentially small,
and the large-amplitude oscillations set in,
at $\mu={\rm min}(\mu_h(x),\mu_{\rm stbc}(x))$
at each point.
On the central portion of 
the interval, $\mu_h(x) < \mu_{\rm stbc}(x)$,
{\it i.e.,} the homogeneous exit time curve (red)
lies to the left of the space-time buffer curve (black).
Hence, this solution fits in Case 4 of DHB,
since the oscillations begin to set in 
at $\mu=-\mu_0$ at $x=0$, where $\vert A_h(x,\mu) \vert$ first
stops being exponentially small and grows to one,
well before $A_p$ can transition.
Then, outside this central portion,
$\mu_{\rm stbc}(x) < \mu_h(x)$,
{\it i.e.,} 
the space-time buffer curve (black) occurs before
the homogeneous exit time curve (red).
See also the inset 
in Figure~\ref{fig:exittime-determined-by-Ah}.

The second example of Case 4 of DHB is illustrated in
Figure~\ref{fig:gaussian-DHBcase4},
for \eqref{eq:gen-CGL}
with $I_G(x)$.
For the given choice of initial data
at $\mu_0=-0.2$, 
which lies inside the interval $(-\omega_0,0)$
with $\omega_0=0.5$,
we observe that $\mu_h(x) < \mu_{\rm stbc}(x)$
and $\mu_h(x) \sim -\mu_0$
for all $x \in [-\ell,\ell]$.
That is,
$A_h(x,\mu)$ transitions from being exponentially small
to large at $\mu_h(x) \sim -\mu_0$, 
at which time $A_p$ is still exponentially small.
Therefore, in this case,
the time at which the solution exits from a neighbourhood
of the repelling QSS is $-\mu_0$ to leading order
on $[-\ell,\ell]$.
The $\mathcal{O}(\eps)$ corrections to $\mu_h(x) \sim -\mu_0$
are spatially periodic,
and by zooming in on the homogeneous exit time curve
one can see these wiggles, as well.
Also, at each point $x$,
the large-amplitude oscillations in $A(x,t)$
which set in after $\mu$ reaches $\mu_h(x)$ 
do not propagate spatially in this case,
with a Gaussian source term.

\bigskip
\noindent
{\bf Remark.}
The spatial dependence of the DHB in Case 4
for solutions
with initial data at $\mu_0 \in (-\omega_0,\delta]$
is a non-trivial extension to PDEs
of what is known 
for analytic ODEs with initial conditions given at the same time.
Consider for example,
the Shishkova ODE (aka Stuart-Landau ODE
with slowly-varying bifurcation parameter).
As mentioned above, it corresponds to 
setting $d=0$ in \eqref{eq:gen-CGL} 
and replacing $I_a(x)$
with an analytic function $h(\mu)$ 
satisfying $h(-i\omega_0)\ne 0$, see
\cite{HKSW2016,N1987,S1973}).
For solutions with initial conditions
given at $\mu_0$,
where $\mu_0$ is any value in $(-\omega_0,-\delta]$,
the exit time from a neighbourhood 
of a repelling slow manifold is $-\mu_0>0$
to leading order.
This is because, for each $\mu_0 \in (-\omega_0,0)$,
there is a Stokes line 
in the complex plane that connects it 
to the point $-\mu_0$ on the positive real axis,
without any saddle point or turning point in between.
Recall Figure~\ref{fig:Contours}.
Hence, for these solutions,
the exit time of $-\mu_0$
occurs before the buffer point
created by the particular solution,
and one says that
the entry-exit function (aka way-in way-out function)
of this analytic ODE 
is symmetric to leading order
for any solution 
with $\mu_0 \in (-\omega_0,-\delta]$.
For the PDE,
the same elliptic contour is used,
however the exit time is generally spatially dependent.

\bigskip

To conclude this section with the examples of DHB,
we observe that there is good agreement between the theory
derived for the linear PDE \eqref{eq:lin-CGL}
and the results in all of the numerical simulations
of the nonlinear PDE \eqref{eq:gen-CGL} which we carried out. 
This indicates 
that, for the nonlinear PDE, the cubic terms in $A$ 
are higher order,
and this will be confirmed by the analysis
of the nonlinear terms
in Section~\ref{sec:gen-CGL}.
Moreover, we note that in this respect the phenomenon
of DHB in the CGL PDE is similar to that
for DHB in the Shishkova ODE 
and in other analytic ODEs,
where the cubic and other nonlinear terms are also higher order.
See for example Section 3 of \cite{HKSW2016}.

%------------------------------------------------------------------------------
\section{DHB and the space-time buffer curve for the cubic CGL}   \label{sec:gen-CGL}
%------------------------------------------------------------------------------

In this section,
we build on the results
for the linear CGL equation \eqref{eq:lin-CGL}
established in Section~\ref{sec:stbc}
to study the full nonlinear CGL equation 
\eqref{eq:gen-CGL}
in the base case
in which
$\beta=\frac{1}{2}$
and $\gamma=1$,
\begin{equation} \label{eq:cCGL}
\eps A_{\mu} = (\mu + i \omega_0) A - (1+i\alpha)\vert A\vert^2 A 
+ \eps^{\frac{1}{2}}I_a(x) + \eps d A_{xx},
\end{equation}
with complex-valued $d$ and $d_R>0$,
and with $I_a(x) \ge 0$ for all $x$.
We demonstrate that solutions
with initial data $A(x,\mu_0)=A_0(x)$ given at time $\mu_0<-\omega_0$ in DHB Case 1 stay near the attracting and repelling quasi-stationary states for $\mu\in[\mu_0,\omega_0 - \tilde\delta)$,
for some $\tilde\delta>0$ small but $\mathcal{O}(1)$ with respect to $\eps$. 
Hence, the nonlinear solution
exhibits DHB and
the space-time buffer curve
for this nonlinear equation
is the same to leading order
as the curve \eqref{eq:stbc} 
for the linear CGL equation
\eqref{eq:lin-CGL}.

We use the same dependent variable 
given by \eqref{eq:AB},
$
B(x,\mu) = A(x,\mu) e^{-\frac{1}{2\eps} (\mu+ i \omega_0)^2},
$
to transform the cubic CGL to
\begin{equation} \label{eq:B-cubic}
\eps B_{\mu} 
= - (1+i\alpha) E(\mu) \vert B \vert^2 B  
+ \sqrt{\eps} I_a(x)e^{-\frac{1}{2\eps}(\mu + i \omega_0)^2}
+ \eps d B_{xx},
\end{equation} 
where
\begin{equation}\label{eq:E}
E(\mu) = e^{\frac{1}{\eps}(\mu_R^2 - (\mu_I + \omega_0)^2)}.
\end{equation}
Next, we subtract off the linear particular solution $B_p$, recall \eqref{eq:Bp-mu-sec23}, substituting $B = B_p + b$ into \eqref{eq:B-cubic} to obtain
\begin{align}\label{e:nl00}
\eps b_\mu = \eps d b_{xx} - (1+i\alpha) E(\mu) (B_p + b)|B_p+b|^2,\qquad  \mu\in [\mu_0,-\mu_0].
\end{align}
We suppress the $x$ dependence in the solutions to keep the formulas more manageable. As it is needed throughout this section, we note the general expansion of the nonlinearity $N[b] := b|b|^2$
$$
N[f+g] = (f+g)^2(\overline{f+g}) = f^2\bar f + 2f \bar f g + f^2 \bar g + g^2 \bar f + 2 f g \bar g + g^2 \bar g.
$$
We consider mild solutions of \eqref{e:nl00} using the variation of constants formula
\begin{align}
b(\mu) &= G_d(\mu - \mu_0)*\left( e^{-\frac{1}{2\eps}(\mu_0 + i\omega_0)^2} A_0 \right) -\frac{(1+i\alpha)}{\eps}\int_{\mu_0}^\mu E(\tilde\mu)G_d(\mu - \tilde\mu)*N\left[B_p(\tilde\mu)+b(\tilde\mu)\right]d\tilde\mu,\label{e:mild0}
\end{align}
where $G_d(\mu-\mu_0)$ denotes the Green's function
and $*$ denotes the convolution. 
We let $H\left[v\right](\mu):=B_h(\mu) + \tilde H\left[v\right](\mu)$ denote the right member of this equation with
\begin{align}
\tilde H\left[v\right](\mu)&:= -\frac{(1+i\alpha)}{\eps}\int_{\mu_0}^\mu E(\tilde\mu)G_d(\mu - \tilde\mu)*N\left[v(\tl\mu)\right]d\tilde\mu,\notag\\
B_h(\mu)&:= G_d(\mu - \mu_0)*\left( e^{-\frac{1}{2\eps}(\mu_0 + i\omega_0)^2} A_0 \right).
\end{align}
We shall assume that the initial data $A_0$ is bounded and the inhomogeneity $I_a$ is smooth with uniformly bounded derivatives. That is, we assume there exists a constant $C>0$ with
 \begin{equation}\label{e:iaa}
 |\partial_x^j I_a(x)|\leq 
 \begin{cases}
 C,& \qquad |\omega_0|\geq 1\\
 C|\omega_0|^2,&\qquad |\omega_0|<1
 \end{cases}
 \qquad j\in \mathbb{N}_0, x\in \mathbb{R}.
 \end{equation}
  This is a rather strong assumption which allows us to readily bound remainder terms occurring below uniformly in $x$. We strongly suspect that similar results can be obtained for less restrictive assumptions on $I_a$.

\subsection{Iterative framework and base iterate}
To construct an approximate solution to the mild formulation \eqref{e:mild0}, we use an iterative approach. We set $b_0(\mu):= 0$ and then iteratively define
\begin{equation}
b_{j+1}(\mu):= H\left[B_p+b_j\right](\mu), \qquad j\geq 0.\label{e:bj}
\end{equation}
In this section, we estimate $b_1$. We claim
\begin{equation}\label{eq:b1est}
b_1(\mu) = e^{-\frac{1}{2\eps}(\mu+i\omega_0)^2}\left(-\eps^{3/2} \frac{(1+i\alpha)I_a(\cdot)^3}{(\mu + i\omega_0)^2(\mu^2+\omega_0^2)} +\mc{O}(\eps^{5/2})\right) + \mathcal{O}(e^{\frac{1}{2\eps}(\omega_0^2 - \mu_0^2)}),\quad \mu\in[-\mu_0,\omega_0 - \tl\delta)
\end{equation}
for some $\tl\delta>0$ fixed and small, with error terms uniform in $x$. This gives the leading order terms in \eqref{eq:b1est}. To obtain this estimate, we note the linear term $B_h$ defined in \eqref{eq:B_h} is exponentially small for all $\mu\in [\mu_0,-\mu_0]$ provided $A_0$ is bounded.  Hence, it suffices to estimate the nonlinear term.

The formula for $B_p$ is given by \eqref{eq:Ia1+Ia2} for $\mu \in [\mu_0,-\delta]$, by the formulas in Appendix \ref{sec:App-B} for $\mu \in
(-\delta,\delta)$, and by \eqref{eq:Bp-mu-sec23} for $\mu\in [\delta, \omega_0]$. Overall, for all
$\mu \in [\mu_0,\omega_0]$, we may write the asymptotic expansion for $B_p$ as
$$
B_p = e^{-\frac{1}{2\eps}(\mu+i\omega_0)^2}\left(-\eps^{1/2} \frac{I_a(\cdot)}{\mu+i\omega_0}+\mathcal{O}(\eps^{3/2})\right) + \tl c(\mu)\left( g(\cdot,\mu+i\omega_0)+\mathcal{O}(\eps^{1/2})\right)=: B_{p,1} + B_{p,2}
$$
where $\tl c(\mu)$ is a bounded, monotonic function with $\tl c(\mu)\equiv0$ for $\mu\leq-\delta$, $\tl c(\mu)\rightarrow\sqrt{\pi/2}$ as $\mu\rightarrow0^+$, and $\tl c(\mu)\equiv \sqrt{2\pi}$ for all $\mu\geq\delta$. Here, $B_{p,1}$ is the linear contribution to the QSS and $B_{p,2}$ is defined to be the term that arises in the same solution for $\mu \ge 0$ due to passing through the saddle at $-i\omega_0$.
We remark that, by the assumptions on $I_a$, the error terms are uniform in $x$.

Next, re-write the expansion as
$$
B_p = e^{-\frac{1}{2\eps}(\mu+i\omega_0)^2}\left(-\eps^{1/2} \frac{I_a(\cdot)}{\mu+i\omega_0}+\mathcal{O}(\eps^{3/2}) + e^{\frac{1}{2\eps}(\mu+i\omega_0)^2}(\tl c(\mu) g(\cdot,\mu+i\omega_0)+ \mc{O}(\eps^{1/2}))\right).
$$
Also note that this factorization, which moves $ e^{-\frac{1}{2\eps}(\mu+i\omega_0)^2}$ outside of all terms in $B_p$, illuminates what remains when transitioning back to the $A$ coordinates. Since $\tl c(\mu)g(\cdot,\mu+i\omega_0)$ is $\mathcal{O}(1)$, the corresponding term in $A$-coordinates, $e^{\frac{1}{2\eps}(\mu+i\omega_0)^2}\tl c(\mu)g(\cdot,\mu+i\omega_0)$ is exponentially small for $\mu\in[\mu_0,\omega_0 - \tl\delta]$ for some fixed $\tl\delta >0$ which is small but $\mc{O}(1)$ with respect to $\eps$.

%estimate \tl H[B_p], and derive an estimate for

To estimate the nonlinear term, 
we use this expansion and work separately on
$\mu<0$ and on $\mu>0$, beginning with the former,
\begin{align}
 \tilde H[B_p]&:= -\frac{1+i\alpha}{\eps} \int_{\mu_0}^\mu E(\tl \mu) G_d(\mu - \tl\mu)*\left[ B_p(\tl \mu) |B_p(\tl\mu)|^2  \right] d\tl\mu  \notag\\
 &=-\frac{1+i\alpha}{\eps} \int_{\mu_0}^\mu \left(-\frac{\tl\mu + i\omega_0}{\eps}e^{-\frac{1}{2\eps}(\tl\mu+i\omega_0)^2}\right) \left(\frac{-\eps}{\tl\mu + i\omega_0} G_d(\mu - \tl\mu)*[-\eps^{3/2} \frac{I_a(\cdot)^3}{(\tl\mu + i\omega_0)(\tl\mu^2+\omega_0^2)}+\mathcal{O}(\eps^{5/2}) ]\right)d\tl\mu
 \notag\\
 &= -[e^{-\frac{1}{2\eps}(\tl\mu+i\omega_0)^2}\eps^{3/2} G_d(\mu - \tl\mu)*\frac{(1+i\alpha)I_a(\cdot)^3}{(\tl\mu + i\omega_0)^2(\tl\mu^2+\omega_0^2)} +\mathcal{O}(\eps^{5/2})\Big]_{\tl\mu = \mu_0}^\mu \notag\\
 &\qquad\qquad+\int_{\mu_0}^\mu e^{-\frac{1}{2\eps}(\tl\mu+i\omega_0)^2} \p_{\tl\mu}\left(G_d(\mu - \tl\mu)*[\eps^{3/2} \frac{(1+i\alpha)I_a(\cdot)^3}{(\tl\mu + i\omega_0)^2(\tl\mu^2+\omega_0^2)}+\mathcal{O}(\eps^{5/2}) ]\right)d\tl\mu\notag\\
 &= e^{-\frac{1}{2\eps}(\mu+i\omega_0)^2}\left(-\eps^{3/2} \frac{(1+i\alpha)I_a(\cdot)^3}{(\mu + i\omega_0)^2(\mu^2+\omega_0^2)} +\mc{O}(\eps^{5/2})\right)+\mc{O}(e^{\frac{1}{2\eps}(\omega_0^2 - \mu_0^2)}).\label{eq:bht}
\end{align}
Note that in the second line, we multiplied the integrand by one in a
useful manner and used the approximation \eqref{eq:Ia1+Ia2}. 
In the third line, we have integrated by parts; and, in the fourth line, we have used the fact that the $\tl\mu = \mu_0$ boundary term is exponentially small (in particular $\mc{O}(e^{\frac{1}{2\eps}(\omega_0^2 - \mu_0^2)})$) while the remaining integral is $\mc{O}(\eps^{5/2})$, uniformly in $x$.  This last claim can be obtained by integrating by parts once more and using the fact that the imaginary part of the phase $-(2\eps)^{-1}(\tl\mu + i\omega_0)^2$ is non-stationary for $\tl\mu\in[\mu_0,\mu]$; see for example \cite[\S 6]{BO1999}.

A similar estimate holds for $\mu\in[0,\omega_0 - \delta]$, as the term $B_{p,2}$ only contributes exponentially small effects here relative to $B_{p,1}$, for $\mu\leq\omega_0 - \tl\delta$. To see this in more detail, we estimate a few of the terms in the expansion of $N[B_{p,1}+B_{p,2}]$. For instance, consider the term $B_{p,2}|B_{p,2}|^2$:
\begin{align}
\Bigg\|-\frac{1+i\alpha}{\eps}& \int_{\mu_0}^\mu E(\tl \mu) G_d(\mu - \tl\mu)*\left[ \tilde c(\tilde\mu)^3 g(\cdot,\tilde\mu+i\omega_0)|g(\cdot,\tilde\mu+i\omega_0) |^2 +\mathcal{O}(\eps^{1/2})\right] d\tl\mu \Big\|_{L^\infty}\notag\\
&=\Bigg\| -\frac{1+i\alpha}{\eps} \int_{-\delta}^\mu E(\tl \mu) G_d(\mu - \tl\mu)*\left[ \tilde c(\tilde\mu)^3 g(\cdot,\tilde\mu+i\omega_0)|g(\cdot,\tilde\mu+i\omega_0) |^2+\mathcal{O}(\eps^{1/2}) \right] d\tilde\mu \Bigg\| \notag\\
&\leq |1+i\alpha| \eps^{-1}E(\mu) \tilde c(\mu)^3 \int_{-\delta}^\mu  \| G_d(\mu - \tl \mu)*\left[g(\cdot,\tilde\mu+i\omega_0)|g(\cdot,\tilde\mu+i\omega_0) |^2 +\mathcal{O}(\eps^{1/2})\right]\|_{L^\infty}d\tilde\mu\notag\\
&\leq C E(\mu)(\eps^{-1} +\mathcal{O}(1))\|I_a\|_{L^\infty}^3\label{eq:b2-3}
\end{align}
for some constant $C>0$, possibly dependent on $\omega_0$. We recall that $\tl c(\mu)$ increases monotonically. Note the last line remains exponentially small for $\eps\ll1$, uniformly for $\mu\in [0,\omega_0 - \tilde\delta]$.
In these inequalities we have repeatedly used the estimate on the heat evolution $\| G_d(\mu - \tl \mu)*f\|_{L^\infty} \leq C \|f\|_{L^\infty}$ for $\mu - \tl\mu\geq0$.  Terms which are quadratic in $B_{p,2}$ can also be bound in a similar way by a term of the form $C\eps^{-1}E(\mu)^{1/2}\|I_a\|_{L^\infty}^2$. 

It remains to consider terms which are quadratic in $B_{p,1}$. For example, consider the term $2B_{p,2}|B_{p,1}|^2$, where we can estimate
\begin{align}
\Bigg\|-\frac{1+i\alpha}{\eps}& \int_{\mu_0}^\mu E(\tl \mu) G_d(\mu - \tl\mu)*\left[ 2B_{p,2}|B_{p,1} |^2 \right] d\tl\mu \Big\|_{L^\infty}\notag\\
&=2\Bigg\| -\frac{1+i\alpha}{\eps} \int_{-\delta}^\mu  G_d(\mu - \tl\mu)*\left[ (\tilde c(\tilde\mu) g(\cdot,\tilde\mu+i\omega_0) + \mc{O}(\eps^{1/2}))\cdot\left|\frac{\eps^{1/2} I_a}{\tilde\mu+i\omega_0} + \mc{O}(\eps^{3/2})\right|^2  \right] d\tilde\mu \Bigg\|_{L^\infty} \notag\\
&\leq 2|1+i\alpha|  \tilde c(\mu)\int_{-\delta}^\mu  \Big\| G_d(\mu - \tl \mu)*\left[g(\cdot,\tilde\mu+i\omega_0)\cdot\left|\frac{ I_a}{\tilde\mu+i\omega_0}\right|^2+ \mc{O}(\eps)\right]\Big\|_{L^\infty}d\tilde\mu\notag\\
&\leq C\|I_a\|_{L^\infty}^3 + \mathcal{O}(\eps).\label{eq:b2-1}
\end{align}

Hence, we may conclude the desired estimate \eqref{eq:b1est} by using \eqref{eq:bht} for $\mu<0$, the similar estimate \eqref{eq:b2-3} for $\mu>0$, and the estimates of the terms quadratic in $B_{p,1}$, such as \eqref{eq:b2-1}, as well as by
noting that the terms coming from $B_{p,2}$ are exponentially small relative to $e^{-\frac{1}{2\eps}(\mu+i\omega_0)^2}$ on $\mu\in[-\delta,\omega_0-\tl\delta)$, and thus contained in the $\mc{O}(\eps^{5/2})$ term in \eqref{eq:b1est}.

\subsection{Inductive step}
We claim inductively that
\begin{align}
b_{j+1}(\mu)= b_j(\mu) + e^{-\frac{1}{2\eps}(\mu+i\omega_0)^2}\left( C_{j+1}(\mu,\omega_0)(\eps^{1/2}I_a)^{2(j+1)+1}+\mathcal{O}(\eps^{\frac{2j+5}{2}})  \right),\label{e:bj0}
\end{align}
where $C_j$ is function of $\mu\in \mathbb{R}$ and $\omega_0$ for which $I_a^{2j+1} C_j$ is uniformly bounded in $x$ and for real $\mu$. By \eqref{eq:b1est}, the claim holds for $j=0$ with $C_1(\mu,\omega_0) = -\frac{1+i\alpha}{(\mu+i\omega_0)^2(\mu^2+\omega_0^2)}$. Note, that by \eqref{e:iaa} we have that $|C_1(\mu,\omega_0) I_a^3|\leq C$ for some fixed constant $C>0$. 

We assume formula \eqref{e:bj0} holds for all $0\leq k\leq j$, and prove it holds for $k = j+1$. Expand
\begin{align}
b_{j+1} &= B_h(\mu) +\left(\tilde H[B_p+b_j] - \tilde H[B_p+b_{j-1}]\right) + \tilde H[B_p+b_{j-1}]\notag\\
&= b_j +\left(\tilde H[B_p+b_j] - \tilde H[B_p+b_{j-1}]\right)\notag\\
&= b_j -\frac{(1+i\alpha)}{\eps} \int_{\mu_0}^\mu E(\tl\mu) G_d(\mu - \tl\mu)*(N[B_p+b_j] - N[B_p+b_{j-1}])d\tl\mu
\end{align}
where the difference of nonlinearities above can be expanded as
$$
N[B_p+b_j] - N[B_p+b_{j-1}] = 2|B_p|^2(b_j-b_{j-1}) + B_p^2(\bar b_j - \bar b_{j-1}) + \bar B_p (b_j^2 - b_{j-1}^2) + 2 B_p(|b_j|^2 - |b_{j-1}|^2) +( b_j^2\bar b_j - b_{j-1}^2\bar b_{j-1}).
$$
Also, note that by our inductive hypothesis, and the fact that the homogeneous term $B_h$ is exponentially small, we find
$$
b_j - b_{j-1} = e^{-\frac{1}{2\eps}(\mu + i\omega_0)^2}\left(  C_j(\mu,\omega_0) (\eps^{1/2}I_a)^{2j+1} + \mathcal{O}(\eps^{\frac{2j+3}{2}}) \right).
$$
Next, we notice that the leading order terms in $\eps$ of $N[B_p+b_j] - N[B_p+b_{j-1}]$ are
\begin{align}
2|B_p|^2(b_j-b_{j-1})&= 2E(\mu)^{-1}e^{-\frac{1}{2\eps}(\mu+i\omega_0)^2}\left(  \frac{C_j(\mu,\omega_0)}{\mu^2+\omega_0^2} (\eps^{1/2}I_a)^{2j+3} + \mathcal{O}(\eps^{\frac{2j+5}{2}}) \right), \\
B_p^2(\bar b_j - \bar b_{j-1}) &= E(\mu)^{-1}e^{-\frac{1}{2\eps}(\mu+i\omega_0)^2}\left( \frac{\overline{C_j(\mu,\omega_0)}e^{4i\omega_0\mu}}{(\mu+i\omega_0)^2} (\eps^{1/2}I_a)^{2j+3} + \mathcal{O}(\eps^{\frac{2j+5}{2}}) \right). 
\end{align}
Defining $\tl C_j(\mu,\omega_0) =  \frac{2C_j(\mu,\omega_0)}{\mu^2+\omega_0^2} + \frac{\overline{C_j(\mu,\omega_0)}e^{4i\omega_0\mu}}{(\mu+i\omega_0)^2}$, inserting the expansions for the leading
order terms, and using integration by parts we find
\begin{align}
b_{j+1}(\mu)-b_j(\mu) 
&= -\frac{1+i\alpha}{\eps} \int_{\mu_0}^\mu E(\tl \mu) G_d(\mu - \tl\mu)*\left[ N(B_p+b_j) - N(B_p+b_{j-1})\right]d\tl\mu\notag\\
&= -\frac{1+i\alpha}{\eps} \int_{\mu_0}^\mu e^{-\frac{1}{2\eps}(\tl\mu+i\omega_0)^2} G_d(\mu - \tl\mu)*\left[ \eps^{\frac{2j+3}{2}}  \tl C_j(\tl\mu,\omega_0)  I_a^{2j+3} + \mc{O}(\eps^{\frac{2j+5}{2}})\right]d\tl\mu\notag\\
&= e^{-\frac{1}{2\eps}(\mu+i\omega_0)^2}\frac{(1+i\alpha)\tl C_j(\mu,\omega_0)}{(\mu+i\omega_0)} \eps^{\frac{2j+3}{2}}I_a^{2j+3} + \mathcal{O}(\eps^{\frac{2j+5}{2}}) \notag\\
 &\qquad - \frac{(1+i\alpha)}{\eps}\int_{\mu_0}^\mu e^{-\frac{1}{2\eps}(\tl\mu+i\omega_0)^2} \partial_{\tl \mu} \left( \frac{\eps}{(\tl\mu+i\omega_0)} G_d(\mu - \tl\mu)*\left[ \eps^{\frac{2j+3}{2}}  \tl C_j(\tl\mu,\omega_0)  I_a^{2j+3} + \mc{O}(\eps^{\frac{2j+5}{2}})\right]\right)d\tl\mu\notag\\
 &= e^{-\frac{1}{2\eps}(\mu+i\omega_0)^2}C_{j+1}(\mu,\omega_0) \eps^{\frac{2j+3}{2}}I_a^{2j+3} + \mathcal{O}(\eps^{\frac{2j+5}{2}}) 
\end{align}
where we have defined $C_{j+1}(\mu,\omega_0) =\frac{(1+i\alpha)\tl C_j(\mu,\omega_0)}{(\mu+i\omega_0)} $. Observe that due to the definition of $\tilde C_j$, we have that $|C_{j+1}|\leq C |\omega_0|^{-3(j+1)-1}$. Then, by the assumption \eqref{e:iaa}, and by the boundedness of $C_jI_a^{2j+1}$, the term $C_{j+1} I_a^{2j+3}$ is bounded in $\mu$ and $x$, uniformly in $\omega_0$. Hence, the difference $b_{j+1} - b_j$ is $\mc{O}(\eps^{\frac{2j+3}{2}})$, which becomes small as $j\rightarrow+\infty$.

From this iterative process, we observe that 
each approximation $b_j$ successively reveals 
the $\eps^{\frac{2j+1}{2}}$-order terms 
in the expansion of the nonlinear attracting and repelling QSSs 
for $\mu<0$ and $\mu>0$, respectively. 
Furthermore, this formal iterative method 
makes it clear that solutions 
with bounded data $B_0(x,\mu_0)$ 
remain exponentially close to the QSSs 
for all $\mu\in[\mu_0,\omega_0 - \tl\delta)$. 
As $\mu$ approaches $\omega_0$ from below, 
while terms coming from $B_{h}$ remain exponentially small, 
the $B_{p,2}$ terms coming from tracking the solution 
over the saddle point 
are no longer exponentially small in the original $A$-coordinates. 
For the $x$-dependent value of $\mu \ge \omega_0$ 
given to leading order by $\mu_{\rm stbc}(x)$, 
they induce the delayed Hopf bifurcation. 
Furthermore, we conclude in Case 1 of DHB that,
since the contribution from the nonlinear terms is higher order, 
$A_p(x,\mu)$ 
mediates the spatially dependent bifurcation.
A similar analysis may be done in the other cases of DHB.

%------------------------------------------------------------------------------
\section{The $\mathcal{O}(\eps)$ value of $\mu_{\rm Hopf}(x)$ in the base case}
\label{sec:Hopf-basecase}
%------------------------------------------------------------------------------

In the analysis 
of the base case of the PDE \eqref{eq:gen-CGL}
($\beta=1/2$ and $\gamma=1$)
in Sections~\ref{sec:lin-CGL}-\ref{sec:DHB-fourcases},
we used that 
$\mu_{\rm Hopf}(x)= 0$ for all $x$
to leading order.
In this section, 
we calculate the $\mathcal{O}(\eps)$ term
in the value of $\mu_{\rm Hopf}(x)$,
and we identify the role this asymptotically small correction
plays in determining the spatial dependence 
of the observed onset of oscillations.
The calculations here are performed 
for general sources $I_a(x)$,
and examples are given with
Gaussian and spatially-periodic terms.

Recall from formula \eqref{eq:QSS-basecase} 
that the attracting and repelling QSS 
on $\mu<-\delta$ and $\mu>\delta$,
respectively, are
$\mathcal{O}(\sqrt{\eps})$ to leading order.
Thus, the linearisation about the small-amplitude QSS
is consistent.
We set $A = A_{\rm QSS} + \eps^a \mathcal{A}$,
with $a>1/2$,
and the linearised equation for $\mathcal{A}$ is
\begin{equation}
\label{eq:basecase-lin-calA}
\eps \mathcal{A}_{\mu}
=(\mu + i\omega_0) \mathcal{A}
-(1+i\alpha) \left( 2 \vert A_{\rm QSS} \vert^2 \mathcal{A} + A_{\rm QSS}^2 \bar{\mathcal{A}}\right)
+\eps d \mathcal{A}_{xx}.
\end{equation}
In terms of the real and imaginary parts,
$\mathcal{A} = \mathcal{U} + i \mathcal{V}$ and
$A_{\rm QSS}(x) = u_{\rm Q} + i v_{\rm Q}$,
the linearised equation for $\mathcal{A}$ 
may be expressed as
$$
\eps
\begin{bmatrix}
{\mathcal{U}}_\mu \\
{\mathcal{V}}_\mu
\end{bmatrix}
= M 
\begin{bmatrix}
\mathcal{U} \\
\mathcal{V}
\end{bmatrix}
+ \eps
\begin{bmatrix}
d_R {\mathcal{U}}_{xx} - d_I {\mathcal{V}}_{xx}\\
d_R {\mathcal{V}}_{xx} + d_I {\mathcal{U}}_{xx}
\end{bmatrix}
$$
where
$\mu = \mu_R + i \mu_I$ and
$$
M = \begin{bmatrix}
\mu_R - 3 u_{\rm Q}^2 - v_{\rm Q}^2 
+ 2\alpha u_{\rm Q}v_{\rm Q}
& 
(\mu_I+\omega_0) + \alpha u_{\rm Q}^2 
+3 \alpha v_{\rm Q}^2 -2u_{\rm Q}v_{\rm Q}\\
(\mu_I+\omega_0) - 3\alpha u_{\rm Q}^2 
- \alpha v_{\rm Q}^2 -2u_{\rm Q}v_{\rm Q}
& 
\mu_R - u_{\rm Q}^2 - 3 v_{\rm Q}^2 
- 2\alpha u_{\rm Q}v_{\rm Q}
\end{bmatrix}.
$$
The trace of $M$ is
\begin{equation}\label{eq:traceM}
{\rm tr} (M) 
= 2 \mu_R - 4 (u_{\rm Q}^2 + v_{\rm Q}^2)
= 2 \mu_R - \frac{4 \eps I_G^2(x)}{\mu_R^2 + (\mu_I +\omega_0)^2}
+\mathcal{O}(\eps^2),
\end{equation}
and
$$
{\rm det}(M)
= (\mu_R^2 + (\mu_I+\omega_0)^2)
-4 \mu_R (u_{\rm Q}^2 + v_{\rm Q}^2)(1+\mu_I+\omega_0)
+3(1+\alpha^2)(u_{\rm Q}^2 + v_{\rm Q}^2)^2,
$$
so that ${\rm det}(M) > 0$
for all $\mu_R<0$,
as well as for
a range of values of $\mu_R>0$.

\begin{figure}
   \centering
   \includegraphics[width=5in]{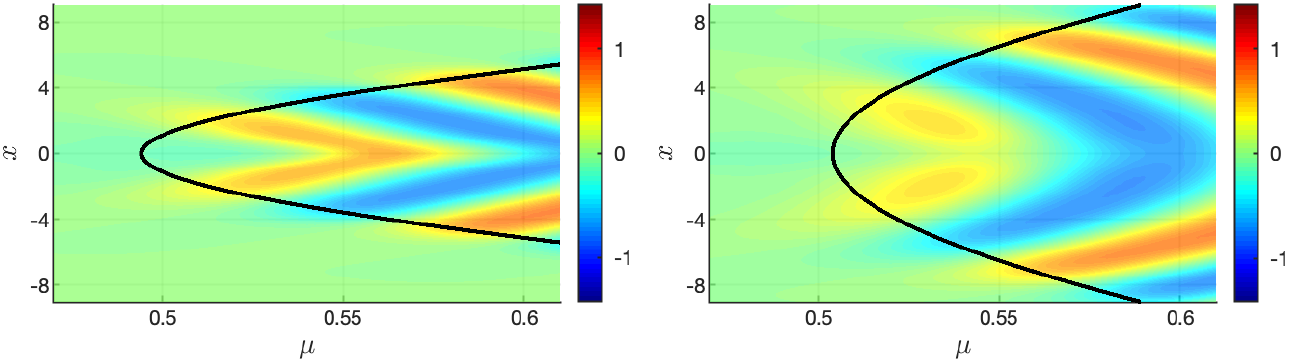}
   \caption{
Onset of oscillations
near the tip of the space-time buffer curve,
of \eqref{eq:lin-CGL}
in a neighbourhood of $x=0$ in the base case
$\beta=1/2$ and $\gamma=1$ of \eqref{eq:gen-CGL}
with Gaussian source ($\sigma=\frac{1}{4}$).
By \eqref{eq:basecase-Hopf},
$\mu_{\rm Hopf}(x)
= \frac{2 \eps (I_G(x))^2}{\omega_0^2} 
+\mathcal{O}(\eps^2)$.
The maximum is at $x=0$,
and $\mu_{\rm Hopf}(x)$ decays rapidly
for $\vert x \vert >0 $.
Hence, near $x=0$,
the delay in the onset of oscillations
at each $x$ is slightly longer,
creating the ``fork in the tongue".
In contrast, 
away from the center of the domain,
the magnitude of $\mu_{\rm Hopf}(x)$ is negligible,
and hence the space-time buffer curve \eqref{eq:stbc-Gaussian}
determines the DHB and the delayed onset of the oscillations there.
Also, this figure illustrates
the effect of the logarithmic terms
at $\mathcal{O}(\eps)$ in \eqref{eq:stbc-2}
and \eqref{eq:stbc-Gaussian}.
Namely, (a) with $d_R = 1$ and $d_I=0$,
the entire buffer curve is shifted leftward
from $\mu=\omega_0=\frac{1}{2}$
by an $\mathcal{O}(\eps)$ amount, 
as is most visible near the tip,
and (b) with $d_R=3$ and $d_I=1$,
the logarithmic terms shift the buffer curve rightward.
The parameters are
$\eps=0.01, \omega_0=\frac{1}{2}$,
and $\alpha=0$. 
The initial data at $\mu_0= -1$ 
is $A_0(x) = -\sqrt{\eps} \frac{I_G(x)}{\mu_0+i \omega_0}$.
}
   \label{fig:Gaussian-tip-of-stbc}
\end{figure}

Therefore,
the Hopf bifurcation curve
for the solutions of \eqref{eq:gen-CGL},
which is obtained
by setting ${\rm tr}(M)=0$,
is given to leading order by
\begin{equation}
\label{eq:basecase-Hopf}
\mu_{\rm Hopf} (x)
= \frac{2 \eps (I_a(x))^2}{\omega_0^2} 
+\mathcal{O}(\eps^2).
\end{equation}
This asymptotic formula holds 
for general 
$I_a(x)$ in the base case of the PDE \eqref{eq:gen-CGL}.

\begin{figure}[h]
   \centering
   \includegraphics[width=5in]{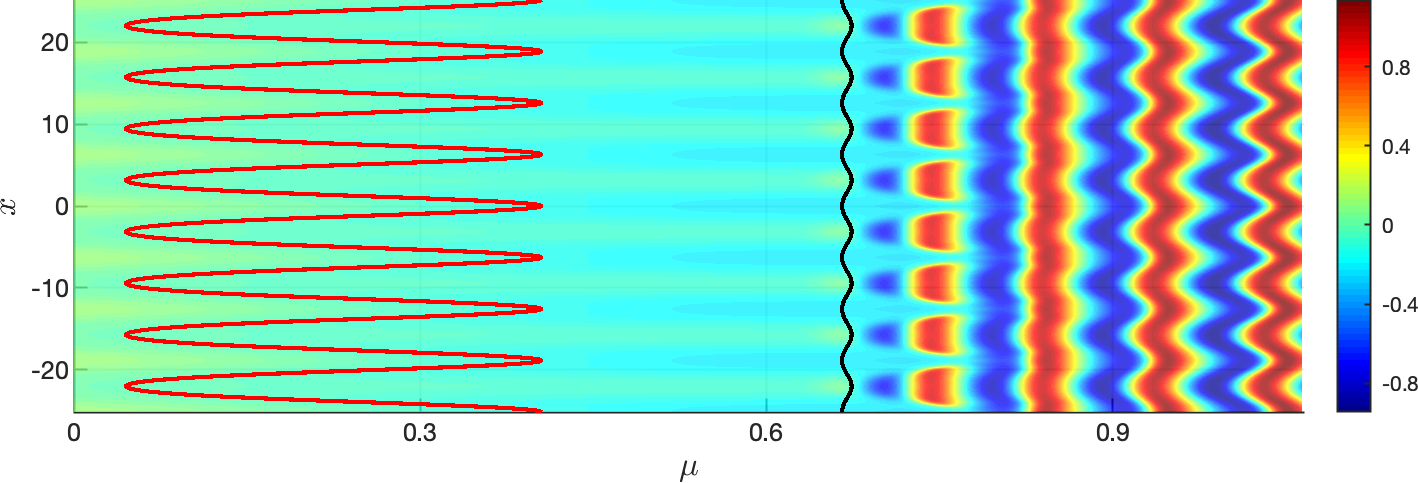} \\ \vspace{5pt}
   \includegraphics[width=2.875in]{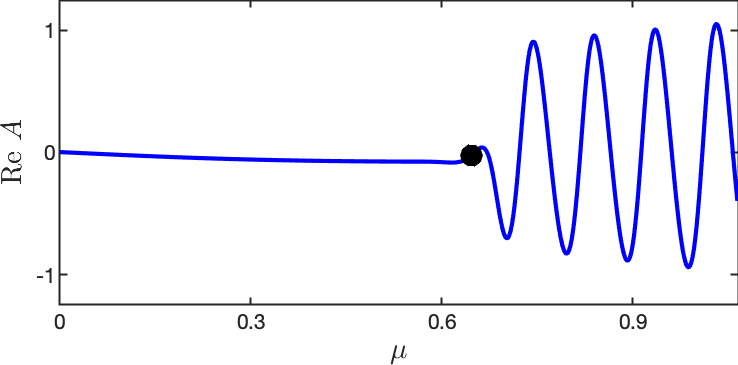} \hspace{0.05in}
   \includegraphics[width=2.875in]{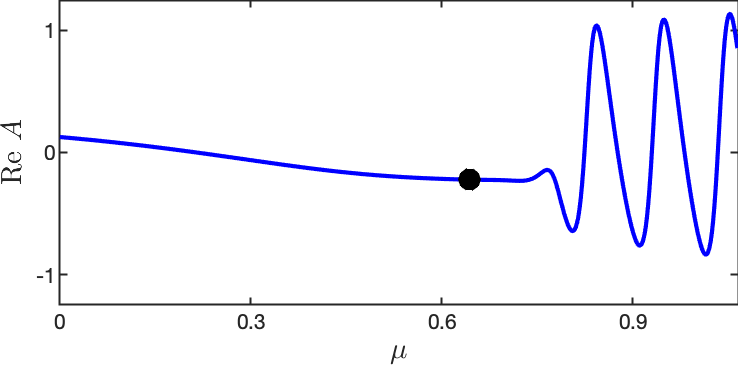}  
   \caption{
${\rm Re}(A(x,\mu))$ 
for an example of Case 1 of DHB
in the PDE \eqref{eq:gen-CGL}
with spatially-periodic source term,
$I_{\rm per}(x)=2+\cos(x)$.
To the left of the exact space-time buffer curve
(black curve),
the solution is near the small-amplitude,
repelling, spatially-periodic QSS
(yellow and green spatially-periodic pattern).
To the right, the oscillations have set in
to leading order,
with the amplitude of the oscillations
increasing rapidly
from small (yellow and
medium blue)
to large (red and dark blue).
The Hopf bifurcation curve
$\mu_{\rm Hopf}(x)$ (red curve)
is super-imposed.
At all points $x$,
this term creates an $\mathcal{O}(\eps)$ increase
in the duration of the DHB and the onset of the oscillations
past the space-time buffer curve.
At the minima, $x=(2k+1)\pi$ for each $k$,
the onset occurs first,
and the maxima,
which occur at $x=2k\pi$ for each $k$,
show a pronounced delay beyond the space-time buffer curve,
since the numerical values of $\mu_{\rm Hopf}(x)$
are approximately 0.4 at the maxima, which is of the same
size as $4\sqrt{\eps}$ here.
The lower frames show time traces at $x = 5\pi$ 
(a local maximum of the space-time buffer curve)
and $x=6\pi$ (a local minimum
of the space-time buffer curve),
illustrating the hard onset
of the oscillations
(the black markers).
The effect of $\mu_{\rm Hopf}(x)$ is particularly visible 
at even integer multiples of $\pi$.
Here, $\eps=0.01, \omega_0=\frac{2}{3}, d_R = 1, d_I = 0$,
and $\alpha=0$.
The initial data at $\mu_0= -1$
is $A_0(x) = -\sqrt{\eps} \frac{I_{\rm per}(x)}{\mu_0+i \omega_0}$.
Similar results are observed for other $A_0(x)$ 
and for other values of $\mu_0< -\omega_0$.
}
   \label{fig:periodic-Hopf}
\end{figure}

The spatial dependence of $\mu_{\rm Hopf}(x)$ 
is illustrated with two different source terms.
First, in Figure~\ref{fig:Gaussian-tip-of-stbc},
we show the results obtained 
with a Gaussian source term, $I_G(x)$.
In a small interval about $x=0$,
the solution of the PDE
remains near the repelling QSS (green region 
just inside the tip of the space-time buffer curve)
for an amount of time equal to 
$\mu_{\rm Hopf}(x)$.
Here, the onset of the oscillations is delayed
slightly longer
than predicted by the space-time buffer curve
by that same amount of time,
and this is manifested by the ``fork in the tongue"
centered at $x=0$.
Then, for $x$ further away from $x=0$,
the amplitude of the Gaussian source term, $I_G(x)$,
is negligibly small. 
Hence, here the function
$\mu_{\rm Hopf}(x)$ is negligibly small, 
and the onset of the oscillations 
coincides with the space-time buffer curve.
In the numerical simulations, 
the small-amplitude oscillations
(light yellow and light blue)
are observed right before the space-time buffer curve,
just as is the case for DHB in analytic ODEs.
Further, one sees that the amplitudes
of the oscillations have become large
(orange, red, dark blue, and purple)
immediately after the space-time buffer curve.

Second, Figure~\ref{fig:periodic-Hopf}
shows the results obtained
with a spatially-periodic source term, $I_{\rm per}(x)$.
The influence
of the $\mathcal{O}(\eps)$
value of $\mu_{\rm Hopf}(x)$,
\eqref{eq:basecase-Hopf},
manifests more here,
since the maximum numerical value of $(I_{\rm per}(x))^2$
over all $x$ is nine, which is of the same
size numerically as $\frac{1}{\sqrt{\eps}}$.

\bigskip
\noindent
{\bf Remark.}
In Figures~\ref{fig:gaussian-DHBcase2},
\ref{fig:case3}, 
and \ref{fig:gaussian-DHBcase4}, 
the $\mathcal{O}(\eps)$ effect
of $\mu_{\rm Hopf}(x)$ is also visible
in the center of the domains,
near the tips of the space-time buffer curve
and the homogeneous exit time curve.

%------------------------------------------------------------------------------
\section{Spatially growing and sign-changing source terms}
\label{sec:ag+sc}
%------------------------------------------------------------------------------

In this section, 
we push somewhat
beyond the theory and examples
for the nonlinear PDE \eqref{eq:gen-CGL} 
in the base case ($\beta=\frac{1}{2}$
and $\gamma=1$),
as presented above
in Sections~\ref{sec:lin-CGL}--\ref{sec:DHB-fourcases}.
There, 
the source terms $I_a(x)$
are taken to be positive with uniformly bounded derivatives
at all points.
Here, we study the PDE \eqref{eq:gen-CGL} in the base case
with an algebraically-growing source 
and with a sign-changing source.

%------------------------------------------------------------------------------
\subsection{An algebraically-growing source term}
\label{sec:ag}
%------------------------------------------------------------------------------

In this section,
we analyze \eqref{eq:gen-CGL} in the base case 
with an algebraically-growing source term
\begin{equation}\label{eq:ag}
I_{\rm ag}(x) = x^2.
\end{equation}

\begin{figure}[h]
   \centering
   \includegraphics[width=5in]{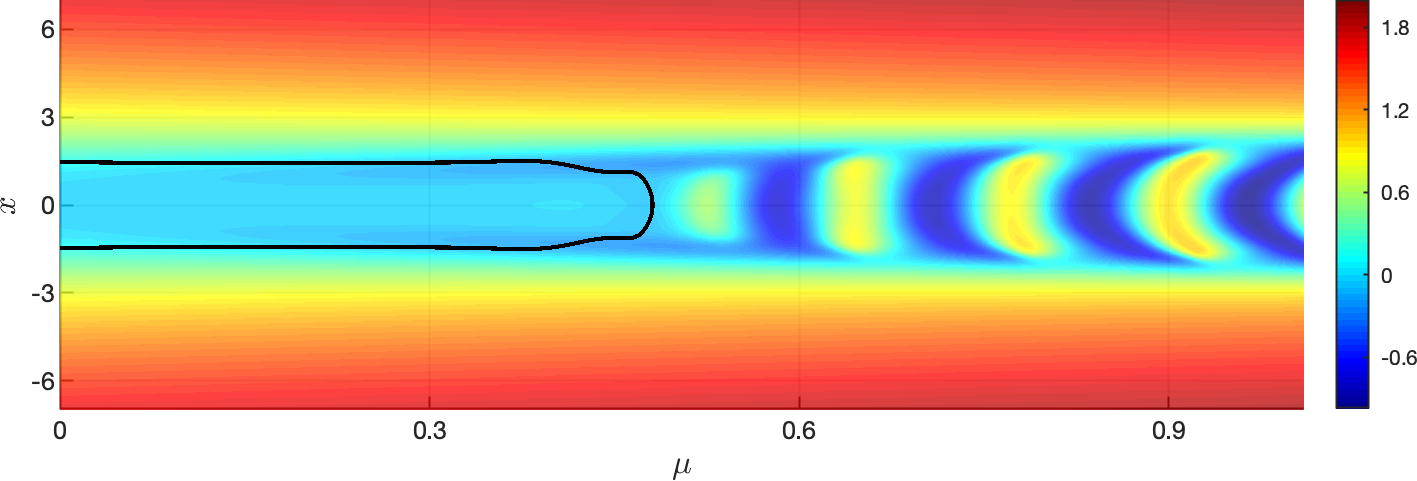}
   \caption{
${\rm Re}(A(x,\mu))$ of the solution
of the cubic CGL \eqref{eq:gen-CGL}
with the algebraically-growing source term,
$I_{\rm ag}(x)$,
along with the exact space-time buffer curve (black curve,
given by the exact solution $A_p(x,\mu)$ in \eqref{eq:Ap-ag}).
DHB and the attendant oscillations are observed
in the central portion of the domain
$( \vert x \vert \lesssim 1.55)$
to the right of the exact space-time buffer curve.
These oscillations propagate outward from $x=0$,
until they reach $\vert x \vert \approx 1.5$.
For $\vert x \vert \approx 1.5$,
$\mu_{\rm stbc}(x) \approx 0$.
Then, outside this central region,
the solution stays near the repelling QSS,
at least until $\mu=10.5$.
Note that the color 
scale differs from that
in the previous figures.
The initial data at $\mu_0= -1$ 
is $A_0(x) = -\sqrt{\eps} \frac{I_{ag}(x)}{\mu_0+i \omega_0}$.
Here, $\omega_0=\frac{1}{2}$, $\alpha=0$,
$d_R=1$,
and $d_I=0$.
}
   \label{fig:alg-growing}
\end{figure}

We start by observing 
that the QSS 
for the cubic PDE \eqref{eq:gen-CGL} with $I_{\rm ag}(x)$
is, to leading order,
\begin{equation}
\label{eq:QSS-ag}
A_{\rm QSS} (x,\mu) = 
\left\{ 
{ \frac{-\sqrt{\eps} x^2}{\mu+i\omega_0} 
\quad {\rm for} \ \ \vert x \vert \ll \eps^{-\frac{1}{4}} } \hfill
\atop
{ \frac{\eps^\frac{1}{6} x^\frac{2}{3} (1-i\alpha)}{(1+\alpha^2)^\frac{2}{3}}
+ \frac{\eps^{-\frac{1}{6}} x^{-\frac{2}{3}}}{3(1+\alpha^2)^\frac{4}{3}}
\left( (3-\alpha^2 -4\alpha i)(\mu+i\omega_0) 
       -2\mu_{\rm R} (1+\alpha^2) \right)
\quad {\rm for} \ \ \vert x \vert \gg \eps^{-\frac{1}{4}} }.
\right.
\end{equation}
Here, we use the fact that 
for $\vert x \vert \ll \eps^{-\frac{1}{4}}$
the QSS is determined to leading order
by balancing the linear term in $A$ 
and the source term in \eqref{eq:gen-CGL},
since the cubic term 
is higher order in this region.
Hence, in this region
with $\vert x \vert \ll \eps^{-\frac{1}{4}}$,
the linearisation
about $A=0$ is valid,
as is the formula 
for the leading order space-time buffer curve.
Also, we note that the higher order terms in the QSS are
$\mathcal{O}(\eps^\frac{3}{2})$ 
and depend on $x$ and $\mu$, 
recall \eqref{eq:QSS-basecase}.

In contrast,
for $\vert x \vert \gg \eps^{-\frac{1}{4}}$,
the QSS in \eqref{eq:QSS-ag} 
is determined by balancing the cubic term 
and the source term in \eqref{eq:gen-CGL}.
Hence, here the QSS has significantly larger amplitude,
and linearisation should be about 
the large-amplitude QSS, 
and no longer about $A=0$.
The higher order terms are
$\mathcal{O}(\eps^{-\frac{1}{2}} x^{-2})$.
Also, to leading order here,
\begin{equation}
\vert A_{\rm QSS}(x,\mu_{\rm R}) \vert 
=  \frac{\eps^\frac{1}{6} x^\frac{2}{3}} {(1+\alpha^2)^\frac{1}{6}}
+ \frac{\eps^{-\frac{1}{6}} x^{-\frac{2}{3}} (\mu_{\rm R} + \omega_0 \alpha)}
       {3(1+\alpha^2)^\frac{5}{6}}.
\end{equation}

\begin{figure}[h]
  \centering
  \includegraphics[width=2.875in]{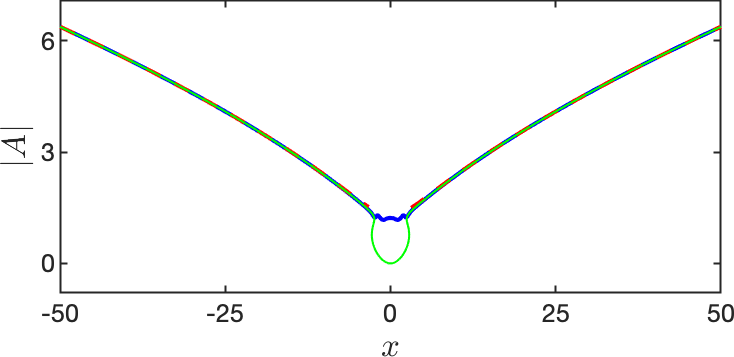} \hspace{0.05in}
  \includegraphics[width=2.875in]{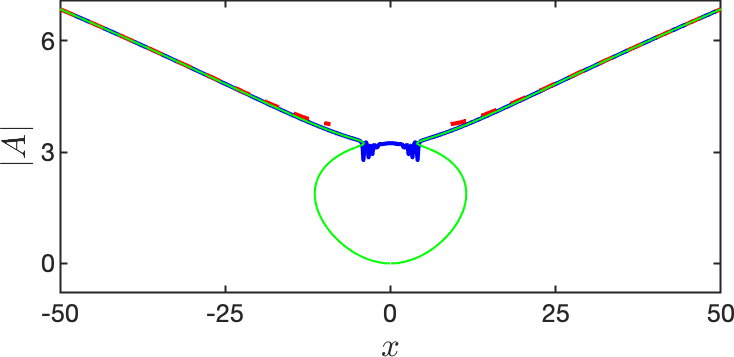}
  \caption{
Spatial profiles 
of $\left| A \right|$ 
(blue curves)
for the same solution
as in Figure~\ref{fig:alg-growing},
at $\mu=1.5$
and $\mu=10.5$.
The oscillations are
seen in the center of the domain,
and the width of this central region 
increases slowly:
$-2 \lesssim x \lesssim 2$ at $\mu=1.5$
and $-3 \lesssim x \lesssim 3$ at $\mu=10.5$.
Outside these central intervals,
the solution is close 
to the repelling QSS,
which is given to leading order by
\eqref{eq:QSS-ag}(b)
(dashed red curve).
Super-imposed is the numerically obtained
solution (green curve)
of the leading order QSS equation,
$r^2 \left( (\mu - r^2)^2 + (\omega_0 - \alpha r^2)^2\right) = \eps x^4,$
valid for all $x$,
as obtained from \eqref{eq:gen-CGL} by neglecting
the time and space derivative terms (which are higher order)
and then squaring the modulus.
Here, $r=\vert A\vert$, $\mu$ is real, and $\alpha=0$.
On $\vert x \vert \gg \eps^{-\frac{1}{4}}$,
the green curve lies essentially on top of the blue and red curves.
For $\mu < \omega_0$,
the green curve does not have any fold points,
and there are no oscillations in the PDE solution,
consistent with the DHB analysis.
Then, the fold points
first appear when $\mu$ reaches $\omega_0$ to leading order,
and they are present (data not shown) for all $\mu$ up
until at least 10.5.
These fold points
delimit the interval
on which the oscillations are observed.
}
   \label{fig:alg-growing-xtraces}
\end{figure}

We now determine the space-time
buffer curve for the region
in which $\vert x \vert \ll \eps^{-\frac{1}{4}}$,
where the QSS has small-amplitude, so that 
the analysis of Section~\ref{sec:stbc} applies.
We require ${\rm Re}(d(\mu-{\tilde{\mu}}))>0$
and $\vert {\rm arg} ( \frac{x}{2d(\mu-{\tilde{\mu}})} ) \vert < \pi$.
Evaluating the integral in definition
\eqref{eq:Bp-def}(b), we find
\begin{equation}\label{eq:g-ag}
g(x,\mu-{\tilde{\mu}})
= x^2 + 2d(\mu-{\tilde{\mu}}).
\end{equation}
With this elementary form
of $g$,
the integral 
\eqref{eq:Bp-def}(a) for $B_p(x,\mu)$
may be evaluated
in closed form in this example,
and hence also $A_p(x,\mu)$ may be found in closed form. 
Specifically, carrying out the integration
in \eqref{eq:Bp-def} and recalling
\eqref{eq:AB},
we find
\begin{equation}\label{eq:Ap-ag}
\begin{split}
A_p(x,\mu)
&= \sqrt{\frac{\pi}{2}} \left(x^2 + 2d(\mu+i \omega_0) \right)
\left[ {\rm erf}\left( \frac{\mu+i\omega_0}{\sqrt{2\eps}} \right)
- {\rm erf}\left( \frac{\mu_0+i\omega_0}{\sqrt{2\eps}} \right) \right]
e^{\frac{1}{2\eps} (\mu + i\omega_0)^2 } \\
+ &{}  2 d \sqrt{\eps} 
\left( 1 - 
e^{\frac{1}{2\eps} (\mu^2 - \mu_0^2 + 2i\omega_0 (\mu - \mu_0))} \right).
\end{split}
\end{equation}
Hence, by taking the real part of the complex-valued, 
space-time-dependent phase of the solution
to be zero, 
we find the exact
space-time buffer curve
for $\vert x \vert \ll \eps^{-\frac{1}{4}}$.
This curve is plotted 
in Figure~\ref{fig:alg-growing},
along with ${\rm Re}(A)$ from the numerical simulation
of the full nonlinear PDE \eqref{eq:gen-CGL} 
with this same source term.
Here, we observe that 
$\eps^{-\frac{1}{4}} \approx 3.16$
for $\eps=0.01$.
Excellent agreement is observed between the onset of the 
oscillations and the exact space-time buffer curve
in the region $\vert x \vert \lesssim 1.55$.

In Figure~\ref{fig:alg-growing-xtraces},
we see that for $\vert x \vert \gg \eps^{-\frac{1}{4}}$,
the solution of the PDE (blue curve)
is near the repelling QSS (red curve)
at least until $\mu=10.5$,
where the QSS is given by \eqref{eq:QSS-ag}(b)
for  $\vert x \vert \gg \eps^{-\frac{1}{4}}$.
Moreover, 
for $\vert x \vert \gg \eps^{-\frac{1}{4}}$,
the Hopf bifurcation (determined
by the linearisation about the non-trivial QSS here,
instead of about $A=0$) is delayed.
See also Section~\ref{sec:largeamp-sourceterm}.

%------------------------------------------------------------------------------
\subsection{A sign-changing source term}
\label{sec:sc}
%------------------------------------------------------------------------------

In this section,
we analyze \eqref{eq:gen-CGL} with
a sign-changing source term
\begin{equation}\label{eq:sc}
I_{\rm sc}(x) = \cos(x).
\end{equation}

\begin{figure}[h]
   \centering
   \includegraphics[width=5in]{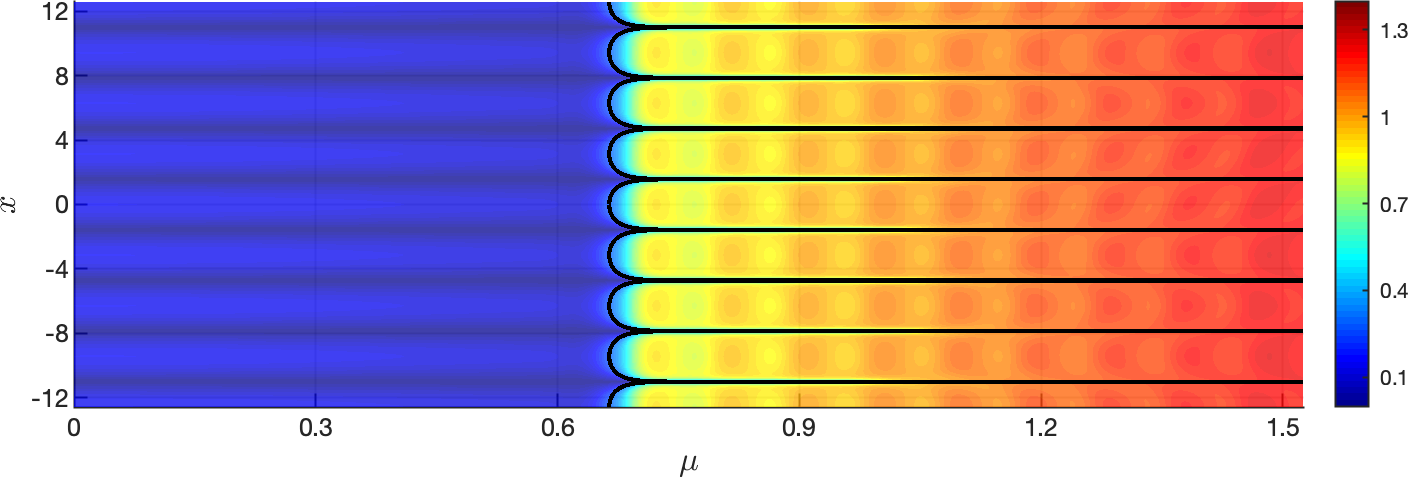}
   \caption{
$\vert A(x,\mu) \vert$ for the solution 
of \eqref{eq:gen-CGL}
with the sign-changing source term,
$I_{\rm sc}(x)=\cos(x)$.
This is an example of Case 2 of DHB.
To the left of the exact space-time buffer curve
(black curve),
the solution is near the small-amplitude,
repelling, spatially-periodic QSS.
About each point
$x = k \pi$,
there is a wide interval
on which $\mu_{\rm stbc}(x) < \mu_h(x)$,
and the hard onset
of the large-amplitude oscillations 
is determined to leading order
by the space-time buffer curve.
On the narrow, complementary intervals
(about $x= \frac{(2n+1)\pi}{2}$),
$\mu_h(x) < \mu_{\rm stbc}(x)$
so that the oscillations 
set in at $\mu_h(x) \sim -\mu_0$,
as determined to leading order
by $A_h$, and there are $\mathcal{O}(\eps)$ amplitude wiggles
about this leading order result,
due to the periodicity of the initial data.
Here, $\eps=0.01, \omega_0=\frac{2}{3}, d_R = 1, d_I = 0$,
and $\alpha=0$. 
The initial data at $\mu_0= -1$ 
is $A_0(x) = -\sqrt{\eps} \frac{I_{\rm sc}(x)}{\mu_0+i \omega_0}$.
}
	   \label{fig:sc}
\end{figure}

One finds
\begin{equation}\label{eq:g-sc}
g(x,\mu-{\tilde{\mu}})
= e^{- d(\mu - {\tilde{\mu}})}\cos(x).
\end{equation}
and
\begin{equation}\label{eq:Ap-sc}
A_p(x,\mu)
= \sqrt{\frac{\pi}{2}} \cos(x) 
\left[ {\rm erf} \left( \frac{\mu+i\omega_0 -\eps d}{\sqrt{2\eps}} \right)
-{\rm erf} \left( \frac{\mu_0+i\omega_0 - \eps d}{\sqrt{2\eps}} \right) \right]
e^{\frac{(\mu+i\omega_0-\eps d)^2}{2\eps}}.
\end{equation}
The resultant space-time buffer curve
(which here is also determined exactly
by setting the real part of the complex-valued,
space-time-dependent phase of $A$ to zero)
is shown in Figure~\ref{fig:sc}.

This is an example of DHB in Case 2.
About each point
$x = k \pi$  (where $\vert \cos(x) \vert = 1$)
there is a wide interval
on which $\mu_{\rm stbc}(x) < \mu_h(x)$,
and the hard onset
of oscillations on these intervals 
is determined to leading order
by the space-time buffer curve.
See Figure~\ref{fig:sc}.
On the complementary intervals,
$\mu_h(x) < \mu_{\rm stbc}(x)$,
so that $A_h$ stops being exponentially small first,
at $\mu_h(x) \sim -\mu_0$, to leading order.
The existence of these narrow intervals
may be understood from the asymptotics of $\mu_{\rm stbc}(x)$.
In particular,
from \eqref{eq:g-sc},
one sees that for this example with a sign-changing source term
there are infinitely many points 
$x = \frac{(2n+1)\pi }{2}$
at which $g(x,\mu+i\omega_0)$ vanishes,
and hence where $\ln \vert g \vert \to - \infty$.
This causes $\mu_{\rm stbc}(x)$
to diverge at these points, recall \eqref{eq:stbc}.
Therefore, for any solution with $\mu_0 \le  -\omega_0$,
there is a (narrow) interval 
about each point $x=\frac{(2n+1)\pi}{2}$
on which $\mu_h(x) < \mu_{\rm stbc}(x)$,
so that $A_h$ 
determines the onset time 
at these points to be
$\mu = \mu_h(x) \sim - \mu_0$.

%------------------------------------------------------------------------------
\section{Asymptotically large source terms}
                  \label{sec:largeamp-sourceterm}
%------------------------------------------------------------------------------

In this section,
we study the CGL equation \eqref{eq:gen-CGL}
with an asymptotically large
$\mathcal{O}(\frac{1}{\sqrt{\eps}})$
source term,
{\it i.e.}, with $\beta = - \frac{1}{2}$ in \eqref{eq:gen-CGL},
\begin{equation} \label{eq:gen-CGL-largeamp-source}
\eps A_{\mu} = (\mu + i \omega_0) A - (1+i\alpha)\vert A\vert^2 A 
+ \frac{1}{\sqrt{\eps}} I_{\tilde a}(x) + \eps d A_{xx},
\end{equation}
while retaining
the $\mathcal{O}(\eps)$ diffusivity, {\it i.e.,} $\gamma=1$,
as in the base case.
The source term, 
which is denoted by $I_{\tilde a}(x)$ in this section,
is taken to be strictly positive.
We find that the Hopf bifurcation
occurs along an $x$-dependent curve
$\mu_{\rm Hopf}(x)$,
and we derive the asymptotics for it,
showing that the large-amplitude source term
causes the bifurcation to occur
well to the right of $\mu=0$.
We quantify how $\mu_{\rm Hopf}(x)$,
together with the space-time buffer curve,
determines the DHB duration
at each $x$,
focusing on regions 
where $I_{\tilde a}(x)=\mathcal{O}(1)$.
We remark that in regions 
where $I_{\tilde a}(x)$ is 
effectively small, {\it i.e.}
of the size of $\mathcal{O}(\eps)$,
then the Hopf term 
only affects the higher order terms.
Overall, the analysis here reveals that, 
by choosing $I_{\tilde a}(x)$ appropriately,
one has region-specific control
over the duration of the DHB,
which can be useful in system design 
for postponing
the onset of undesirable oscillations.

The attracting and repelling QSS 
are given 
on $\mu<-\delta$ and $\mu>\delta$,
respectively, by
\begin{equation} \label{eq:QSS-largeamp-source}
A(x,\mu) 
= \frac{1}{\eps^{\frac{1}{6}}} A_{\rm QSS}(x)
+\eps^{\frac{1}{6}} \mathcal{A}(x,\mu,\eps),
\quad {\rm where} 
\quad
A_{\rm QSS}(x) 
= \frac{1-i\alpha}{(1+\alpha^2)^{\frac{2}{3}}} (I_{\tilde a}(x))^{\frac{1}{3}}, 
\end{equation}
$\mathcal{A} = \mathcal{O}(1)$
for all $x$ and $\mu$ uniformly in $\eps$ for sufficiently small $\eps>0$.
(See also the Remark below.)
The linearised equation for $\mathcal{A}$ is
\begin{equation}
\label{eq:lin-calA}
\eps^{\frac{4}{3}} \mathcal{A}_{\mu}
=\eps^{\frac{1}{3}} (\mu + i\omega_0) \mathcal{A}
-(1+i\alpha) \left( 2 \vert A_{\rm QSS} \vert^2 \mathcal{A} + A_{\rm QSS}^2 \bar{\mathcal{A}}\right)
+\eps^{\frac{4}{3}} d \mathcal{A}_{xx}
+(\mu + i \omega_0) A_{\rm QSS}
+\eps d (A_{\rm QSS})_{xx}.
\end{equation}
In terms of the real and imaginary parts,
$A_{\rm QSS}(x) = u_{\rm Q} + i v_{\rm Q}$ and
$\mathcal{A} = \mathcal{U} + i \mathcal{V}$,
the linearised equation for $\mathcal{A}$ 
may be expressed as
$$
\eps^{\frac{4}{3}}
\begin{bmatrix}
{\mathcal{U}}_\mu \\
{\mathcal{V}}_\mu
\end{bmatrix}
= M 
\begin{bmatrix}
\mathcal{U} \\
\mathcal{V}
\end{bmatrix}
+ \eps^\frac{4}{3} d 
\begin{bmatrix}
{\mathcal{U}}_{xx} \\
{\mathcal{V}}_{xx}
\end{bmatrix}
+ \begin{bmatrix}
\mu_R u_{\rm Q} - (\mu_I+\omega_0) v_{\rm Q}\\
\mu_R v_{\rm Q} + (\mu_I+\omega_0) u_{\rm Q}
\end{bmatrix}
+ \eps \begin{bmatrix}
{\rm Re} (d (A_{\rm QSS})_{xx})  \\
{\rm Im} (d (A_{\rm QSS})_{xx})
\end{bmatrix},
$$
where
$\mu = \mu_R + i \mu_I$ and
$$
M = \begin{bmatrix}
\eps^\frac{1}{3} \mu_R - 3 u_{\rm Q}^2 - v_{\rm Q}^2 + 2\alpha u_{\rm Q}v_{\rm Q}
& 
-\eps^\frac{1}{3} (\mu_I+\omega_0) + \alpha u_{\rm Q}^2 +3 \alpha v_{\rm Q}^2 -2u_{\rm Q} v_{\rm Q}\\
\eps^\frac{1}{3} (\mu_I+\omega_0) - 3\alpha u_{\rm Q}^2 - \alpha v_{\rm Q}^2 
-2u_{\rm Q}v_{\rm Q}
& 
\eps^\frac{1}{3} \mu_R - u_{\rm Q}^2 - 3 v_{\rm Q}^2 - 2\alpha u_{\rm Q}v_{\rm Q}
\end{bmatrix}.
$$
The trace of $M$ is
\begin{equation}\label{eq:traceM-largeampsource}
{\rm tr} (M) 
= 2\eps^{\frac{1}{3}} \mu_R - 4 (u_{\rm Q}^2 + v_{\rm Q}^2),
\end{equation}
and
$$
{\rm det}(M)
= \eps^{\frac{2}{3}} (\mu_R^2 + (\mu_I+\omega_0)^2)
-4 \mu_R \eps^{\frac{1}{3}} (u_{\rm Q}^2 + v_{\rm Q}^2)(1+\mu_I+\omega_0)
+3(1+\alpha^2)(u_{\rm Q}^2 + v_{\rm Q}^2)^2,
$$
so that ${\rm det}(M) > 0$
for all $\mu_R<0$,
as well as for
$\mu_R>0$,
at least until $\mathcal{O}(\eps^{-\frac{1}{3}})$.

\begin{figure}[h]
   \centering
   \includegraphics[width=5in]{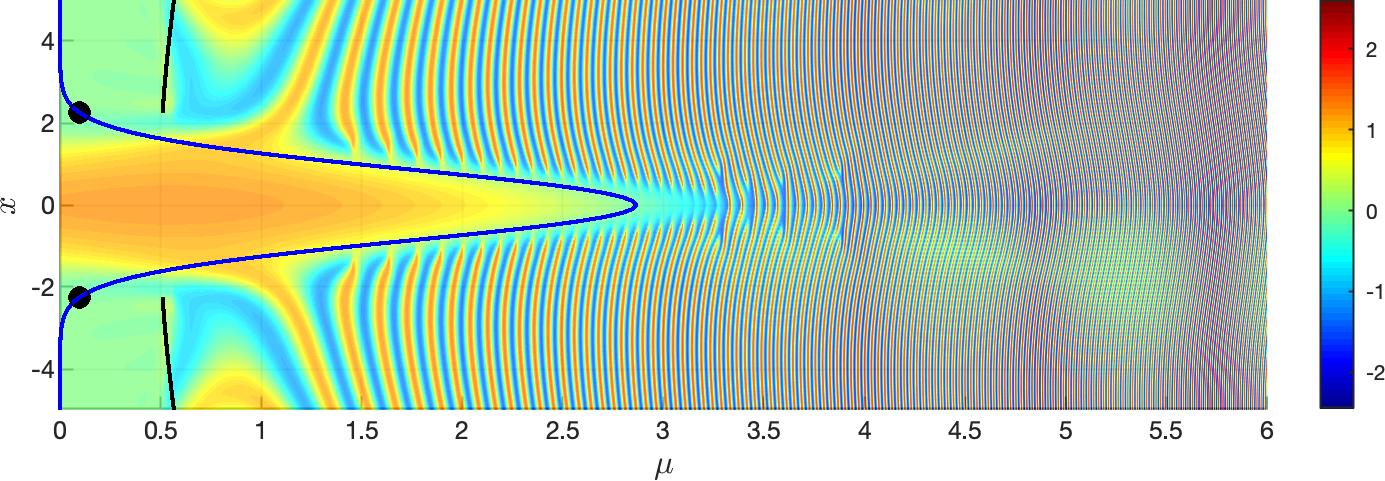}
   \caption{
With the large-amplitude source term
$I_{\tilde a}(x)=\tilde{a} e^{-\frac{x^2}{4\sigma}}$, 
$\tilde{a}=0.2$,
and $\sigma=\frac{1}{4}$,
the instantaneous Hopf bifurcation curve
is $\mu_{\rm Hopf}(x)=(2.8655...) e^{-\frac{2x^2}{3}}$ (blue curve).
For each $x$,
the solution stays near the repelling QSS for some time
(with ${\rm Re}(A)$ being orange-yellow for $\vert x \vert \lesssim 2.24$,
and green for $\vert x \vert \gtrsim 2.24$).
Here, the black dots mark the points
$x = \pm 2.24$,
where $\mu_{\rm Hopf} = \sqrt{\eps}$.
For $\vert x \vert \lesssim 2.24$,
the duration of the delayed onset of the oscillations 
is approximately $\omega_0$ time units
beyond the spatially-dependent curve 
$\mu_{\rm Hopf}(x)$.
Then, for $\vert x \vert \gtrsim 2.24$,
where $\mu_{\rm Hopf} \le \sqrt{\eps}$,
the oscillations set in just before
the spatio-temporal buffer curve (black curve),
as expected since the QSS has small-amplitude here
and the curve is again
given to leading order by
\eqref{eq:stbc-Gaussian} for these $x$.
Simulation performed on 
$[-50,50]$.
For $\vert x \vert > 5$ 
(not shown),
the oscillations 
also commence just before the space-time buffer curve.
The parameters are
$\eps=0.01, \omega_0=\frac{1}{2}, d_R = 1, d_I = 0.5, \alpha=0.6$.
}
   \label{fig:large-amplitude-20}
\end{figure}

\begin{figure}[h]
   \centering
\includegraphics[width=5in]{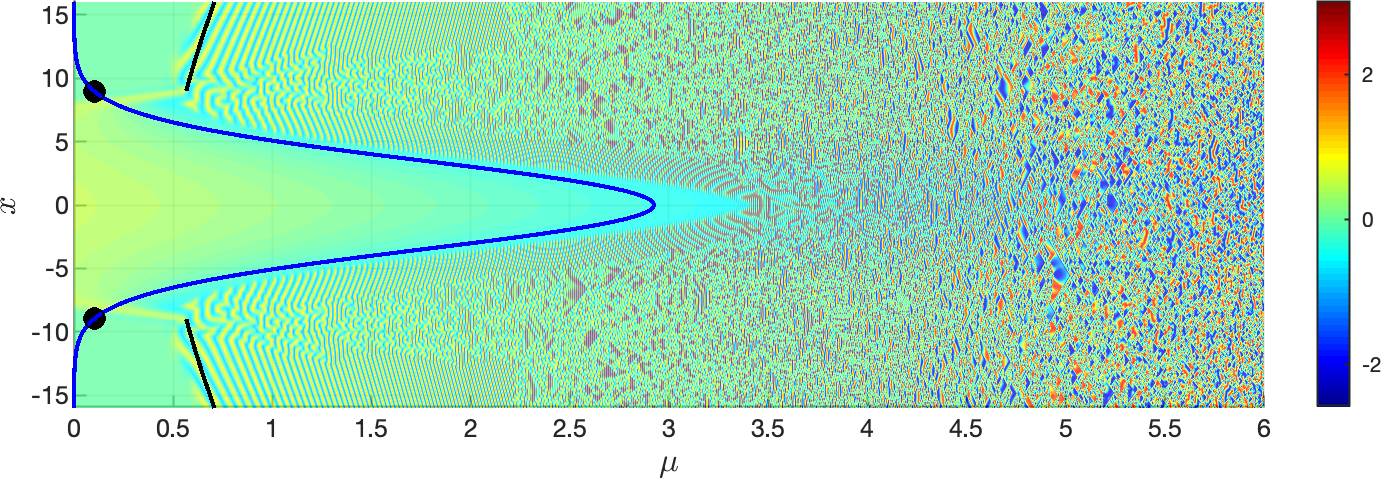}
   \caption{
For this simulation,
$\mu_{\rm Hopf}(x) = (2.9240...) e^{-\frac{x^2}{24}}$
(blue curve).
The parameters are the same 
as in Figure~\ref{fig:large-amplitude-20}, except
$\tilde{a}=0.5$,
$\sigma=4$,
$\alpha=\sqrt{7}$,
and $d_I = 0$.
Compared to the simulation
in Figure~\ref{fig:large-amplitude-20}, 
the magnitude of $\mu_{\rm Hopf}(x)$ is larger
due to the larger source term amplitude $\tilde{a}$,
the half-width $\sigma$ has increased for a larger $\sigma$,
and the space-time dynamics
of the post-onset oscillations
have changed for larger values of $\alpha$.
The space-time buffer curve \eqref{eq:stbc}
(black curve) 
is shown only for those $x$ for which 
$\mu_{\rm Hopf}(x) \ll 1$, {\it i.e.,} where
the linearisation about $A=0$ is valid.
}
   \label{fig:large-amplitude-50}
\end{figure}

\begin{figure}[h]
   \centering
\includegraphics[width=5in]{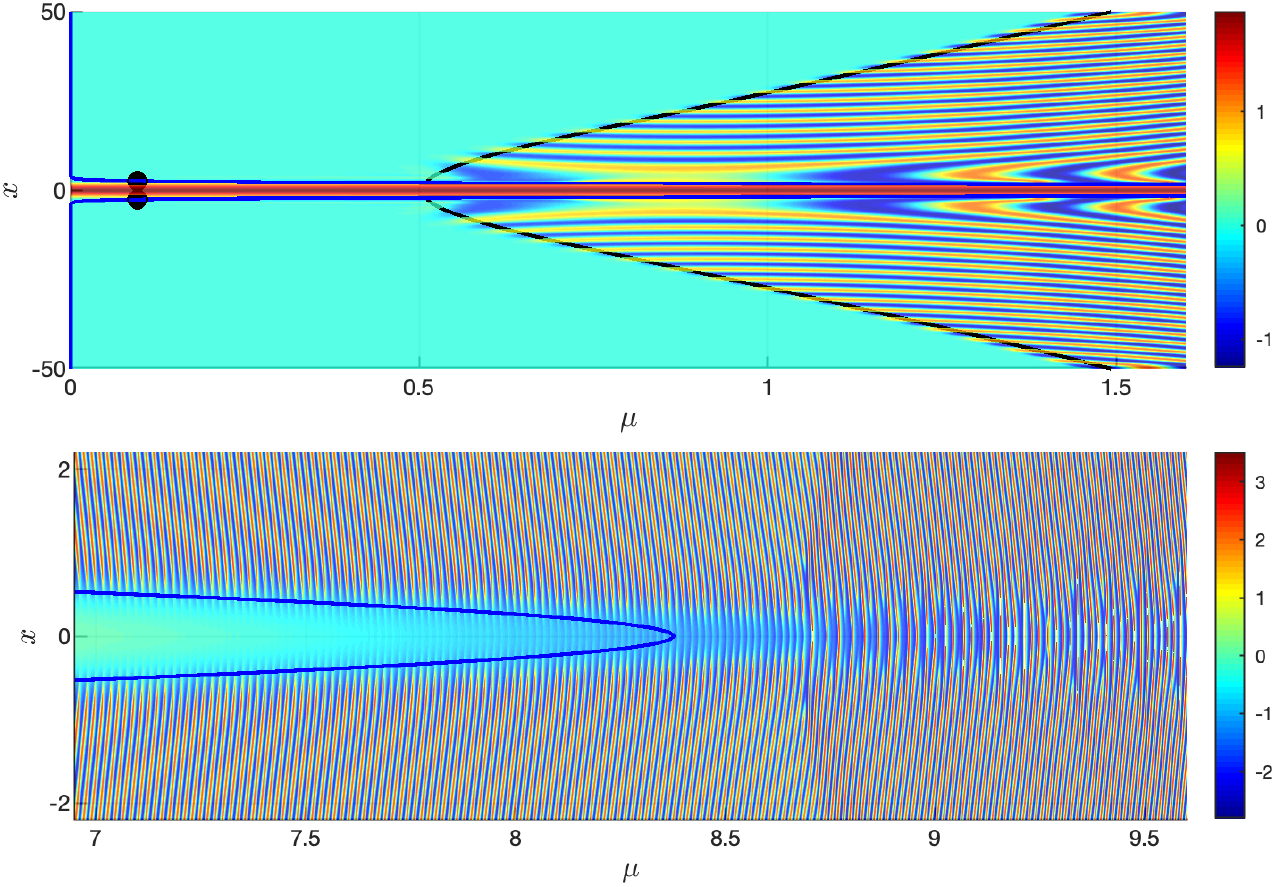}
   \caption{
${\rm Re}(A(x,\mu))$ obtained 
from \eqref{eq:gen-CGL-largeamp-source}
with Gaussian source term.
The QSS has large amplitude 
in an interval about $x=0$.
The space-dependent
Hopf curve
is given by
$\mu_{\rm Hopf}(x) = (8.3788...) e^{-\frac{2x^2}{3}}$ (blue curve).
The points $x$ where 
$\mu_{\rm Hopf}(x)=\sqrt{\eps}$
are marked  by the dots in the upper frame.
The curve $\mu_{\rm Hopf}(x)$
extends to $\mu=8.3788...$,
as is shown in the lower frame.
By the time $\mu$ reaches the tip 
of the Hopf curve, the amplitude of the QSS has decreased
(from red to light green and blue).
Inside the central interval,
the oscillations set in after $\mu$ crosses the blue curve.
In contrast,
outside this interval,
the oscillations set in 
at $\mu_{\rm stbc}(x)$, to leading order,
since $\mu_{\rm Hopf}(x)$ is negligibly small there
(upper frame, 
where the space-time buffer curve (black curve)
is also super-imposed).
The parameters are 
$\tilde{a}=1.0$,
$\sigma=\frac{1}{4}$,
$\eps=0.01$,
$\omega_0=0.5$,
$\alpha=0.6$,
$d_R = 3$,
and $d_I = 1$.
The initial data at $\mu_0=-1.5$
is $A_0(x)=-\sqrt{\eps} \frac{I_G(x)}{\mu_0+i\omega_0}$.
}
   \label{fig:large-amplitude-100}
\end{figure}

Therefore,
for the solutions 
of \eqref{eq:gen-CGL-largeamp-source},
the Hopf bifurcation 
occurs to leading order
at the $x$-dependent time given by
\begin{equation}
\label{eq:sec8-Hopf}
\mu_{\rm Hopf} (x)
= \frac{2}{\eps^{\frac{1}{3}}} (u_{\rm Q}^2 + v_{\rm Q}^2)
= \eps^{-\frac{1}{3}} \frac{2 (I_{\tilde a}(x))^{\frac{2}{3}}}{(1+\alpha^2)^\frac{1}{3}}.
\end{equation}
It is illustrated 
in Figures~\ref{fig:large-amplitude-20},
\ref{fig:large-amplitude-50},
and \ref{fig:large-amplitude-100}.
Also, we have checked (at several points $x$)
that the numerically observed duration of DHB
in the full nonlinear PDE scales as 
$\eps^{-\frac{1}{3}}$ (data not shown).

In Figure~\ref{fig:large-amplitude-20},
we compare the results 
of numerical simulations
of \eqref{eq:gen-CGL-largeamp-source} 
with a Gaussian source term,
$I_{\tilde a}(x)=\tilde{a}e^{-\frac{x^2}{4\sigma}}$,
${\tilde a}=0.2$, and $\sigma=\frac{1}{4}$,
to the analytical results.
The Hopf bifurcation curve
is $\mu_{\rm Hopf} (x) = (2.8655...) e^{-\frac{2x^2}{3}}$,
following \eqref{eq:sec8-Hopf}.
The narrow peak of the repelling QSS 
manifests as the orange and yellow band about $x=0$.
At all points,
the full solution 
stays near the repelling QSS
at least until $\mu$ reaches $\omega_0$.
After that,
the exit time from the neighbourhood
of the repelling QSS 
is $x$-dependent.
In particular, for $\vert x \vert \ge x_* = 2.24...$,
$\mu_{\rm Hopf} (x) \le  0.1=\sqrt{\eps}$.
Hence, the use of the space-time buffer curve
obtained from the linearisation
about $A=0$ (for which $\mu_{\rm Hopf} = 0$ to leading order)
is consistent here,
and we see that, for $\vert x \vert \gtrsim 2.24$,
the time of exit 
from the neighbourhood
of the repelling QSS
(and of the onset of oscillations)
is delayed beyond 
the Hopf bifurcation time
essentially by the amount 
determined by the space-time buffer curve,
given by \eqref{eq:stbc-Gaussian}.

In contrast, 
for $\vert x \vert \lesssim x_*$,
the solution continues to stay near the repelling QSS
for much longer,
and there is a delay 
(beyond
$\mu_{\rm Hopf}(x)$)
of approximately $\omega_0$ in duration
before the oscillations (rapid blue to red transitions)
commence.
For example, at $x=0$, the amplitude 
crosses zero (transition from yellow to green)
near $\mu = 2.8$,
and the oscillations set in near $\mu = 3.3$.

Simulations with other values of 
$\tilde{a}$ (up to and including $\tilde{a}=2$),
$\sigma$,
$\alpha$, 
and $d$ (complex)
show similar results
(over all values simulated)
for the delay in the onset of oscillations.
In the central regions,
the delayed onset occurs 
beyond the Hopf curve, 
$\mu_{\rm Hopf}(x)$,
by an amount approximately equal to $\omega_0$.
Then, outside the central region,
the DHB is given to leading order 
by the space-time buffer curve 
(obtained by linearising the CGL about $A=0$).
See Figure~\ref{fig:large-amplitude-50},
where $\tilde{a}=0.5$,
$\sigma=4$, 
$\alpha=\sqrt{7}$,
$\eps=0.01$,
and the Hopf bifurcation curve
is $\mu_{\rm Hopf}(x) = (2.9240...)e^{-\frac{x^2}{24}}$,
by \eqref{eq:sec8-Hopf}.
See also Figure~\ref{fig:large-amplitude-100},
where $\tilde{a}=1.0$,
$\sigma=\frac{1}{4}$, 
$\eps=0.01$,
$\alpha=0.6$,
$d_R=3$,
$d_I=1$,
and the Hopf bifurcation curve
is $\mu_{\rm Hopf}(x) = (8.3788...)e^{-\frac{x^2}{24}}$,
by \eqref{eq:sec8-Hopf}.

\bigskip
\noindent
{\bf Remark.}
The first correction term $\mathcal{A}$ to the QSS 
in \eqref{eq:QSS-largeamp-source} is
$$
\mathcal{A}(x,\mu,\eps) 
= \frac{1}{3(1+\alpha^2)^{\frac{4}{3}} (I_a(x))^{\frac{1}{3}}}
[ (\mu_R - 3\alpha^2 \mu_R + 4 \alpha (\omega_0 + \mu_I)) 
  + i (3(\omega_0+\mu_I) -\alpha^2 (\omega_0 + \mu_I) - 4\alpha \mu_R) ]
+\mathcal{O}(\eps^{\frac{1}{3}}).
$$

\bigskip
\noindent
{\bf Remark.}
The linearised equation 
\eqref{eq:lin-calA} for $\mathcal{A}$
has $x$-dependent saddle points
(or nilpotent points)
in the complex $\mu$-plane.
These are located 
where ${\rm tr}(M)=0$
and ${\rm det}(M)=0$.
Specifically,
for values of $x$ at which 
$I_{\tilde a}(x)$ is strictly of $\mathcal{O}(1)$,
these saddles occur at
$\mu_R = \mu_{\rm Hopf}(x)$
and $\mu_I= ({\mu_I})_{\pm}$,
where
$({\mu_I})_{+} 
= - \omega_0 + \frac{8 I_{\tilde a}(x)^{\frac{4}{3}}}
                    {\eps^{\frac{2}{3}} (1+\alpha^2)^\frac{2}{3}}
- \frac{3\alpha^2 - 1}{8}
+\mathcal{O}(\eps^{\frac{2}{3}}),$
and $({\mu_I})_{-} 
= - \omega_0 
+ \frac{3\alpha^2 - 1}{8}
+\mathcal{O}(\eps^{\frac{2}{3}})$.
Hence, compared to the case of small-amplitude source terms 
($\beta=\frac{1}{2}$)
for which there is one saddle (at $\mu=-i\omega_0$),
the large-amplitude source term creates a second saddle,
and the maximum and spatial dependence
of $I_{\tilde a}(x)$ determine the saddle locations.
Analysis of the Stokes lines
through these saddles,
especially where they intersect the $\mu_R$-axis,
would enable one to further quantify the DHB in this region,
though the analysis for \eqref{eq:lin-calA} is 
more complex than it is for \eqref{eq:lin-CGL},
where the corresponding matrix $M$ is simpler,
with $M_{11}, M_{22} = \mu_R$,
$M_{12}= - (\mu_I + \omega_0)$, and
$M_{21}= (\mu_I + \omega_0)$.

%------------------------------------------------------------------------------
\section{An example with 
$\mathcal{O}(1)$ diffusivity
and $\mathcal{O}(1)$ source term}
           	  \label{sec:DO1}
%------------------------------------------------------------------------------

In this section,
we extend some of the results of the base case 
of the PDE \eqref{eq:gen-CGL}
to an example of 
$\mathcal{O}(1)$ diffusivity
and $\mathcal{O}(1)$ 
amplitude source term 
in \eqref{eq:gen-CGL}
({\it i.e.}, $\gamma=0$ and $\beta=0$).
The PDE is
\begin{equation} \label{eq:full-DO1}
\eps A_{\mu} = (\mu + i \omega_0) A 
-(1+i\alpha) \vert A \vert^2 A
+ I_a(x) + \hat{d} A_{xx},
\end{equation}
where we now use $\hat{d}$ to denote the diffusivity.
The data 
$A(x,\mu_0)=A_0(x)$
is bounded with sufficiently many continuous derivatives,
with $\mu_0<-\omega_0$
again of primary interest.

%------------------------------------------------------------------------------
\subsection{The $\mathcal{O}(1)$ QSS and $\mu_{\rm Hopf}(x)$ for \eqref{eq:full-DO1} }
\label{sec:DO1-subsec1}
%------------------------------------------------------------------------------

With $\mathcal{O}(1)$ diffusivity $\hat{d}$ 
and source terms $I_a(x)$,
the attracting and repelling QSS
are $\mathcal{O}(1)$,
which contrasts 
with the $\mathcal{O}(\sqrt{\eps})$ 
amplitude of the QSS 
in the base case, recall \eqref{eq:QSS-basecase}. 
For general source terms $I_a(x)$,
one may use variation of constants
to find the QSS of the linearised equation,
followed by an iterative procedure
on the mild form of the PDE to generate
the QSS of the nonlinear PDE.

For example, with the Gaussian 
$I_a(x)=e^{\frac{-x^2}{4\sigma}}$,
the attracting QSS of the linearised version of \eqref{eq:full-DO1}
on $\mu < -\delta$
is found
(using variation of constants)
to be given to leading order by
\begin{equation} \label{eq:dO1-QSS}
\begin{split}
A_{\rm QSS} (x,\mu)
&= 
c_1(x,\mu) 
e^{-\sqrt{\frac{-(\mu+i\omega_0)}{\hat{d}}}x}
+ c_2(x,\mu) 
e^{\sqrt{\frac{-(\mu+i\omega_0)}{\hat{d}}}x},
\quad {\rm with} \\
c_1(x,\mu) 
&=  \frac{1}{2} \sqrt{\frac{\pi\sigma}{-{\hat{d}}(\mu+i\omega_0)}}
e^{-\frac{\sigma (\mu+i\omega_0)}{\hat{d}}}
\left[ 1 + {\rm erf} \left( 
       \frac{x}{2\sqrt{\sigma}} - \sqrt{\frac{-\sigma(\mu+i\omega_0)}{\hat{d}}} 
      \right) \right] \\
c_2(x,\mu) 
&=  \frac{1}{2} \sqrt{\frac{\pi\sigma}{-\hat{d}(\mu+i\omega_0)}}
e^{-\frac{\sigma (\mu+i\omega_0)}{\hat{d}}}
\left[ 1 - {\rm erf} \left( 
       \frac{x}{2\sqrt{\sigma}} + \sqrt{\frac{-\sigma(\mu+i\omega_0)}{\hat{d}}} 
      \right) \right].
\end{split}
\end{equation}
We observe 
$\lim_{x \to -\infty} c_1(x,\mu)=0$
and $\lim_{x \to \infty} c_2(x,\mu)=0$.
A similar formula holds 
for the leading order
repelling QSS on $\mu > \delta$.
Moreover, one can show 
using steepest descents on the integral in $A_p$
(just as for the base case in Section~\ref{sec:lin-CGL}),
that general solutions with $\mu_0< - \omega_0$
stay near the attracting QSS for $\mu < -\delta$
and then near the repelling QSS at least until
$\omega_0$ at each point $x$. 
However, the analysis is more involved,
since the QSSs of the CGL
have $\mathcal{O}(1)$ amplitude.

We now find a formula 
for the $x$-dependent Hopf bifurcation curve $\mu_{\rm Hopf}(x)$.
Let
\begin{equation} \label{eq:QSS-O(1)-source}
A(x,\mu) 
= A_{\rm QSS}(x,\mu)
+\eps^{\frac{1}{3}} \mathcal{A}(x,\mu,\eps),
\end{equation}
where
$\mathcal{A}$ is at most $\mathcal{O}(1)$
uniformly in $\eps$
for all $x$ and $\mu$. 
The linearised equation for $\mathcal{A}$ is
\begin{equation}
\eps \mathcal{A}_{\mu}
=(\mu + i\omega_0) \mathcal{A}
-(1+i\alpha) \left( 2 \vert A_{\rm QSS} \vert^2 \mathcal{A} 
                    + A_{\rm QSS}^2 \bar{\mathcal{A}}\right)
+ {\hat d} \mathcal{A}_{xx}.
\nonumber
\end{equation}
In terms of the real and imaginary parts,
$A_{\rm QSS}(x) = u_{\rm QSS} + i v_{\rm QSS}$
and $\mathcal{A} = \mathcal{U} + i \mathcal{V}$,
the linearised equation for $\mathcal{A}$  may be re-expressed as
\begin{equation}
\label{eq:O(1)-lin-calA}
\eps
\begin{bmatrix}
{\mathcal{U}}_\mu \\
{\mathcal{V}}_\mu
\end{bmatrix}
= M 
\begin{bmatrix}
\mathcal{U} \\
\mathcal{V}
\end{bmatrix}
+ {\hat d}
\begin{bmatrix}
{\mathcal{U}}_{xx} \\
{\mathcal{V}}_{xx}
\end{bmatrix},
\end{equation}
where
$\mu = \mu_R + i \mu_I$ and
$$
M = \begin{bmatrix}
\mu_R - 3 u_{\rm QSS}^2 - v_{\rm QSS}^2 + 2\alpha u_{\rm QSS}v_{\rm QSS}
& 
-(\mu_I+\omega_0) + \alpha u_{\rm QSS}^2 +3 \alpha v_{\rm QSS}^2 -2u_{\rm QSS}v_{\rm QSS}\\
(\mu_I+\omega_0) - 3\alpha u_{\rm QSS}^2 - \alpha v_{\rm QSS}^2 -2u_{\rm QSS}v_{\rm QSS}
& 
\mu_R - u_{\rm QSS}^2 - 3 v_{\rm QSS}^2 - 2\alpha u_{\rm QSS}v_{\rm QSS}
\end{bmatrix}.
$$

Now, the trace of $M$ is
\begin{equation}\label{eq:traceM-O1}
{\rm tr} (M) 
= 2 \mu_R - 4 (u_{\rm QSS}^2 + v_{\rm QSS}^2).
\end{equation}
Hence,
for the solutions of
\eqref{eq:full-DO1},
the Hopf bifurcation 
is given implicitly by
\begin{equation}
\label{eq:sec9-Hopf}
\mu_{\rm Hopf} (x)
= 2 \vert A_{\rm QSS} \vert^2,
\end{equation}
where $A_{\rm QSS} = u_{\rm QSS} + i v_{\rm QSS}$
is evaluated at $(x,\mu_{\rm Hopf}(x))$.
This result for $\mathcal{O}(1)$
diffusivity and amplitude source term
shows that in the regions
where $\mu_{\rm Hopf}(x)$
is $\mathcal{O}(1)$
the QSS changes from being attracting to repelling 
at a value of $\mu$ substantially to the right of zero,
and the DHB needs to be determined from \eqref{eq:O(1)-lin-calA}.
In contrast, for those $\vert x \vert$ 
at which $\mu_{\rm Hopf}(x) \ll 1$, this formula shows that
the changeover in the stability type of the QSS
happens again at $\mu=0$ to leading order,
and the linearisation about $A=0$ 
is again valid.

%------------------------------------------------------------------------------
\subsection{The space-time buffer curve for the linearised version of \eqref{eq:full-DO1}}
\label{sec:DO1-subsec2}
%------------------------------------------------------------------------------

The general solution 
of the PDE obtained by linearising \eqref{eq:full-DO1} about $A=0$
is decomposed 
into homogeneous and particular components,
$A(x,\mu) = A_h(x,\mu) + A_p(x,\mu).$
These two components are derived
using Duhamel's Principle 
(just as in Section~\ref{eq:lin-CGL})
to solve the equation for $B$.
Here, 
\begin{equation}\label{eq:DO1:Ah}
A_h(x,\mu) 
= \frac{\sqrt{\eps}  \  \ 
{\rm exp}[\frac{1}{2\eps} ((\mu+i\omega_0)^2 - (\mu_0 + i\omega_0)^2)]}
        {\sqrt{ 4 \pi \hat{d} (\mu-\mu_0)}}
\int_R {\rm exp}\left[\frac{-\eps(x-y)^2}{4\hat{d}(\mu-\mu_0)}\right] A_0(y) dy.
\end{equation}
Note that with
$\hat{d} = \eps d$ 
in \eqref{eq:DO1:Ah},
one naturally recovers
the formula for the homogeneous solution
\eqref{eq:Ah-basecase} in the base case.

Next, 
in the scaled variable
$\hat{x}=\sqrt{\eps}x$,
the equation
for $B_p$ is
\begin{equation}\label{eq:Bp-lin-xtilde}
\eps (B_p)_\mu 
= \eps \hat{d} (B_p)_{\hat{x}\hat{x}} 
  + I_a(\hat{x}/\sqrt{\eps}) e^{-\frac{1}{2\eps}(\mu+i\omega_0)^2},
\end{equation}
with $B_p(\hat{x}/\sqrt{\eps},\mu_0)=0$.
The solution is 
\begin{equation}\label{eq:DO1-g}
\begin{split}
B_p\left( \frac{\hat{x}}{\sqrt{\eps}}, \mu \right)
&= \
\frac{1}{\eps} 
\int_{\mu_0}^\mu g(\hat{x},\mu - {\tilde{\mu}}) 
e^{-\frac{1}{2\eps}({\tilde{\mu}}+i\omega_0)^2} d{\tilde {\mu}},
\quad {\rm with}
\\
g(\hat{x}=\sqrt{\eps}x,\mu-\tilde{\mu})
&=\frac{1}{\sqrt{4\pi \hat{d} (\mu - \tilde{\mu}) }}
\int_R 
       e^{\frac{-(\hat{x}-\hat{y})^2}{4{\hat{d}}(\mu-\tilde{\mu})}} 
       I_a\left( \frac{\hat{y}}{\sqrt{\eps}} \right) 
       d {\hat{y}}
=\frac{\sqrt{\eps}}{\sqrt{4\pi \hat{d} (\mu - \tilde{\mu}) }}
\int_R  e^{\frac{-\eps(x-y)^2}{4\hat{d}(\mu-\tilde{\mu})}} I_a(y) dy.
\end{split}
\end{equation}
Compare to \eqref{eq:Bp-def}.

Now, using 
the method of stationary phase 
along the contour
$C = C_1 \bigcup C_2 \bigcup C_3 \bigcup C_4$
(recall Section 3),
one finds that
$B_p=\frac{1}{\sqrt{\eps}} \sqrt{2\pi} g(\sqrt{\eps}x, \mu + i \omega_0) 
+\mathcal{O}(1) $,
with the dominant contributions
again coming from the segments of $C_2$ and $C_3$
near the saddle.
Hence, translating 
$B_p(x,\mu)$ via \eqref{eq:AB} to obtain
$A_p(x,\mu)$, 
one arrives at
the following implicit formula
for the space-time buffer curve:
\begin{equation}\label{eq:DO1-stbc}
\left\{ (x,\mu_{\rm stbc}(x)) \vert  \quad
{\rm Re}\left( 
\ln \left(\sqrt{\frac{2\pi}{\eps}}g(\sqrt{\eps}x,\mu_{\rm stbc}(x)+i\omega_0)\right) 
   + \frac{1}{2\eps}(\mu_{\rm stbc}(x) + i \omega_0)^2 \right)
= 0
\right\}.
\end{equation}
This space-time buffer curve
is spatially flatter than \eqref{eq:stbc},
through the argument $\sqrt{\eps}x$ of $g$.

\begin{figure}[h]
   \centering
\includegraphics[width=5in]{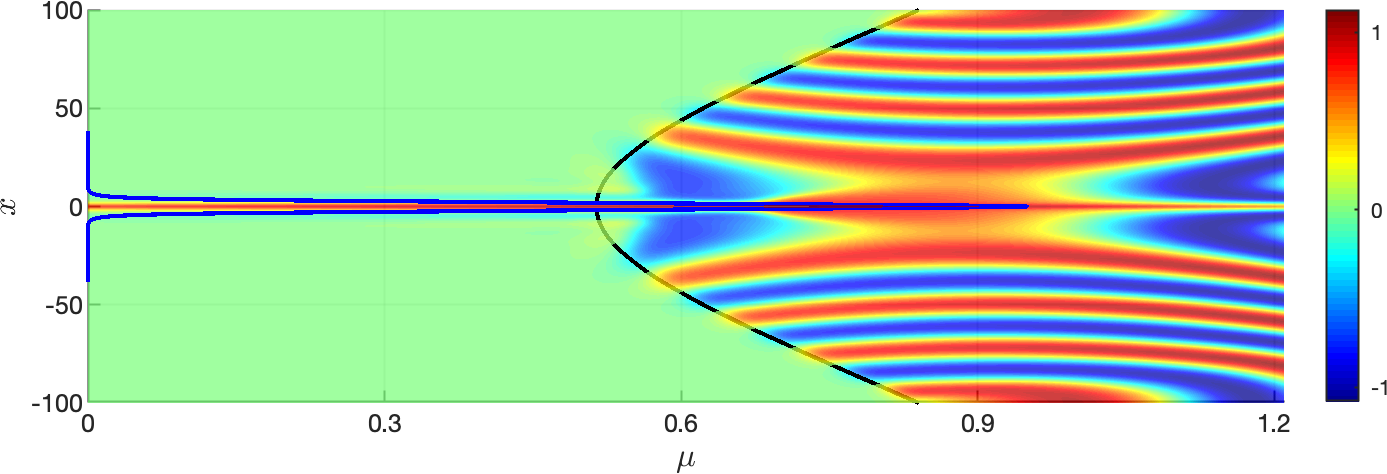}
   \caption{
${\rm Re}(A(x,\mu))$ for the solution of \eqref{eq:full-DO1}
on the spatial interval
$[-100,100]$.
The black space-time buffer curve is
obtained by solving \eqref{eq:DO1-stbc} numerically
with $g$ given by \eqref{eq:DO1-g-Gaussian},
and the blue Hopf bifurcation curve
is obtained from formula \eqref{eq:sec9-Hopf}
and \eqref{eq:dO1-QSS}.
The $\mathcal{O}(1)$ QSS is green 
in most of the region to the left of the space-time buffer curve,
where it has small amplitude;
and the QSS is red in the center,
where it has maximum amplitude.
The parameters are
$\eps=0.01, \omega_0=\frac{1}{2}, {\hat{d}}_R = 1, 
{\hat{d}}_I = 0, \alpha=0.6$,
and $\mu_0=-1$,
with $I_G(x)= e^{-\frac{x^2}{4\sigma}}$, 
and $\sigma=\frac{1}{4}$.
}
   \label{fig:GaussianO1source+diffusivity}
\end{figure}

\bigskip
The above analysis 
of the space-time buffer curve for 
the linearised version of \eqref{eq:full-DO1}
with $\mathcal{O}(1)$ amplitude source term 
and $\mathcal{O}(1)$ diffusivity
may be illustrated using the Gaussian 
source term
$I_G(x)=e^{-\frac{x^2}{4\sigma}}$.
For $\mu$ on $C_4$,
one finds
\begin{equation}\label{eq:DO1-g-Gaussian}
g(x,\mu-{\tilde{\mu}}) 
=  \sqrt{\frac{\eps \sigma}{\hat{d}(\mu+i \omega_0) + \eps\sigma}}
   e^{\frac{-\eps x^2}{4({\hat{d}}(\mu+i \omega_0) + \eps \sigma)}}
   + \mathcal{O}(\eps).
\end{equation}
Then, substitution 
of \eqref{eq:DO1-g-Gaussian}
into \eqref{eq:DO1-stbc}
shows that the leading order space-time buffer curve is
\begin{equation}\label{eq:DO1-spacetime-buffercurve}
\begin{split}
\mu^2 =  \omega_0^2 
&+ \frac{\eps^2 x^2 ({\hat{d}}_R \mu - {\hat{d}}_I \omega_0 + \eps\sigma)}
       {2\left( ({\hat{d}}_R \mu - {\hat{d}}_I \omega_0 + \eps\sigma)^2 
                + ({\hat{d}}_R\omega_0 + {\hat{d}}_I\mu)^2 \right) } \\
& - \eps\ln(2\pi\sigma) 
+ \frac{\eps}{2} \ln \left( ({\hat{d}}_R\mu - {\hat{d}}_I \omega_0 + \eps\sigma)^2
    +  ({\hat{d}}_R \omega_0 + {\hat{d}}_I \mu)^2 \right),
\quad \mu \ge \omega_0.
\end{split}
\end{equation}
An example of the space-time buffer curve
obtained by solving \eqref{eq:DO1-stbc}
numerically with $g$ given by \eqref{eq:DO1-g-Gaussian}
and a Gaussian source term
is plotted
in Figure~\ref{fig:GaussianO1source+diffusivity}.
(Note that with $\hat{d}=\eps d$
it reduces 
to \eqref{eq:stbc-Gaussian}
obtained in the base case with
$\mathcal{O}(\eps)$ diffusivity.)
Compared to
\eqref{eq:stbc-Gaussian},
this space-time buffer curve
for $\mathcal{O}(1)$ diffusivity
defines a spatially flatter 
space-time buffer curve $\mu_{\rm stbc}(x)$,
because the diffusivity is 
one order of magnitude larger here.
Note that the simulation presented 
in Figure~\ref{fig:GaussianO1source+diffusivity}
is on the domain $[-200,200]$,
which is significantly larger than that 
in Figure~\ref{fig:Gaussian-stbc}.
Hence,
the larger the modulus of the diffusivity,
the smaller the magnitude of the spatial contribution,
and the more uniform the delay time becomes.
Also, 
the half-width $\sigma$ of the Gaussian
has a weaker impact.

%----------------------------------------------------------------------------
\section{Conclusions and discussion}		\label{sec:CD}
%----------------------------------------------------------------------------

%----------------------------------------------------------------------------
\subsection{Conclusions}		\label{sec:conclusions}
%----------------------------------------------------------------------------

Considering the prototypical CGL PDE 
\eqref{eq:gen-CGL}
as an equation in its own right,
this article 
has presented a study 
of the phenomenon 
of delayed Hopf bifurcation (DHB)
as the parameter $\mu$ 
increases slowly in time
through an instantaneous Hopf 
bifurcation at $\mu=0$.
It has been shown that 
solutions with initial data
given at $\mu_0 \le -\omega_0$ 
are not only near 
the attracting QSS while $\mu<0$, 
but they remain near the QSS 
as $\mu$ continues to evolve slowly
until well after it has become repelling 
and at least until $\mu$ reaches $\omega_0$.
This analysis 
of the delay of the Hopf bifurcation (DHB) 
was performed 
by directly using the classical methods of stationary phase 
and steepest descents
on the linear PDE,
based on the topography
induced by the saddle point at $\mu=-i\omega_0$,
and followed by using an iterative method 
for solutions of the nonlinear PDE.
Specifically, the nonlinear analysis
is based on an iterative method
for the difference
between the solution $A$ 
of the full cubic PDE
and the particular solution $A_p$ 
of the linear PDE.

Then, with these explicit results,
it was shown that 
there is a competition 
at the heart of DHB
between two exponentially small terms, 
one each
from the particular solution
$A_p(x,\mu)$ 
and the homogeneous solution
$A_h(x,\mu)$, 
to see which component 
first ceases to be exponentially small.
The former stops being exponentially small 
and attains magnitude one
along the space-time buffer curve, 
$\mu_{\rm stbc}(x)$,
and the latter 
along the homogeneous exit time curve,
$\mu_h(x)$. 
Explicit asymptotic formulas 
were derived,
and their properties were illustrated 
with different types of source terms and initial data,
including uni-modal, smoothed step function,
and spatially periodic.
Furthermore, in some of the examples,
it is possible to calculate 
the curves from closed form solutions.

Based on an analysis of different outcomes
of the competition 
between $\mu_{\rm stbc}(x)$ and $\mu_h(x)$,
{\it i.e.,} between which comes first,
several primary cases of DHB were introduced and analyzed.
The first threes cases of DHB are 
for solutions $A(x,\mu)$ 
of \eqref{eq:gen-CGL}
with initial data 
given at $\mu_0 \le -\omega_0$.
Here, Case 1 of DHB arises
when the duration 
of the bifurcation delay
is determined at all points
by $\mu_{\rm stbc}(x)$,
{\it i.e.,} when 
$\mu_{\rm stbc}(x) < \mu_h(x)$ for all $x$.
Case 2 of DHB occurs
when the bifurcation delay 
is determined
at some points by $\mu_{\rm stbc}(x)$
and at others by $\mu_h(x)$.
Case 3 of DHB arises 
when the duration of the delay 
is determined
at all points by $\mu_h(x)$, 
{\it i.e.,}
when $\mu_h(x) < \mu_{\rm stbc}(x)$ 
for all $x$.
Finally, Case 4 of DHB was introduced
for solutions of \eqref{eq:gen-CGL}
with initial data given 
at $-\omega_0 < \mu_0 < -\delta$,
where $\delta>0$ is small but $\mathcal{O}(1)$.
It was shown that, also for these solutions,
the exit time 
from a neighbourhood of the repelling QSS
can be spatially-dependent, as well.

Examples were presented
of the different cases of DHB,
and it was shown how to classify
the DHB for general source terms $I_a(x)$
and various initial data.
The local maxima of the source term and the initial data mark 
the sites at which the solution of the full cubic PDE \eqref{eq:gen-CGL}
first diverges from the repelling QSS,
and where the large-amplitude,
post-DHB oscillations first set in.
The spatial dependence of the DHB and onset of oscillations
was shown to be quadratic
in the case of Gaussians (uni-modal functions),
a smoothed step function in the case of 
source terms given by an error function,
and spatially-periodic 
in the case of spatially-periodic functions.

Finally, extensions 
of the main results
were presented.
Going beyond the main DHB results
established for bounded and positive source terms,
it was shown that 
DHB also occurs in the base case of \eqref{eq:gen-CGL} 
with algebraically-growing
and sign-changing source terms.
The formulas for the space-time buffer 
and homogeneous exit time curves 
(calculated either asymptotically or exactly)
also accurately predict when the oscillations
set in at each point $x$, 
even though the source terms are not bounded
or positive.
Next, for the PDE 
with asymptotically large source terms ($\beta=-\frac{1}{2}$),
it was shown that the instantaneous
Hopf bifurcation curve $\mu_{\rm Hopf}(x)$
can become large, even asymptotically large,
so that the duration of the DHB
can be asymptotically long.
Combined with the information
derived above about how the properties 
of the source terms 
determine the space-time buffer curve,
this provides a high level 
of control or ability to design 
the spatial dependence 
of when the oscillations set in.
Moreover, with large-amplitude source terms,
it was found that 
there is more than one saddle point 
in the complex $\mu$ plane,
and hence the topography 
of the Stokes and anti-Stokes lines
is richer.
A final extension concerns 
the case of $\mathcal{O}(1)$ diffusivity,
for which the space-time buffer curve 
is also derived
and found to be
spatially flatter. 
That the method also extends
to $\mathcal{O}(1)$ diffusivity
enables application
to a broader range of problems,
in which the diffusivities 
are not necessarily small.

There are important considerations about the stability of the
numerical simulations.  As we have demonstrated throughout the
article, the solutions are rapidly oscillating with frequency on the
order of $\frac{\mu \omega_0}{\eps}$.  The numerical stiffness
induced by the combination of rapid oscillations and slow drift in
$\mu$ places an upper bound on the values of $\omega_0$ that can be
used to reliably compute the solution near the repelling QSS.  On the
other hand, our space-time buffer curve predictions require that the
leading order estimate of the spatially-dependent Hopf bifurcation,
$\mu_{\rm Hopf}(x) = \frac{2\eps I_a(x)}{\omega_0^2}$, is small.  As
such, the numerical values of $\omega_0$ that can be used are bounded
from below, for any fixed value of $\eps$.  
Thus, to address the stiffness and satisfy the smallness
of $\mu_{\rm Hopf}(x)$, we have chosen $\omega_0$ values in the range $C
\eps^{1/4} < \omega_0 < 1$, where $C$ is an $\mathcal{O}(1)$
constant. Many of our reported simulations use $0.5 \leq \omega_0
\leq 0.75$ for $\eps=0.01$.

%----------------------------------------------------------------------------
\subsection{Post-DHB spatio-temporal patterns}		\label{sec:Post-DHB}
%----------------------------------------------------------------------------

\begin{figure}
   \centering
\includegraphics[width=6in]{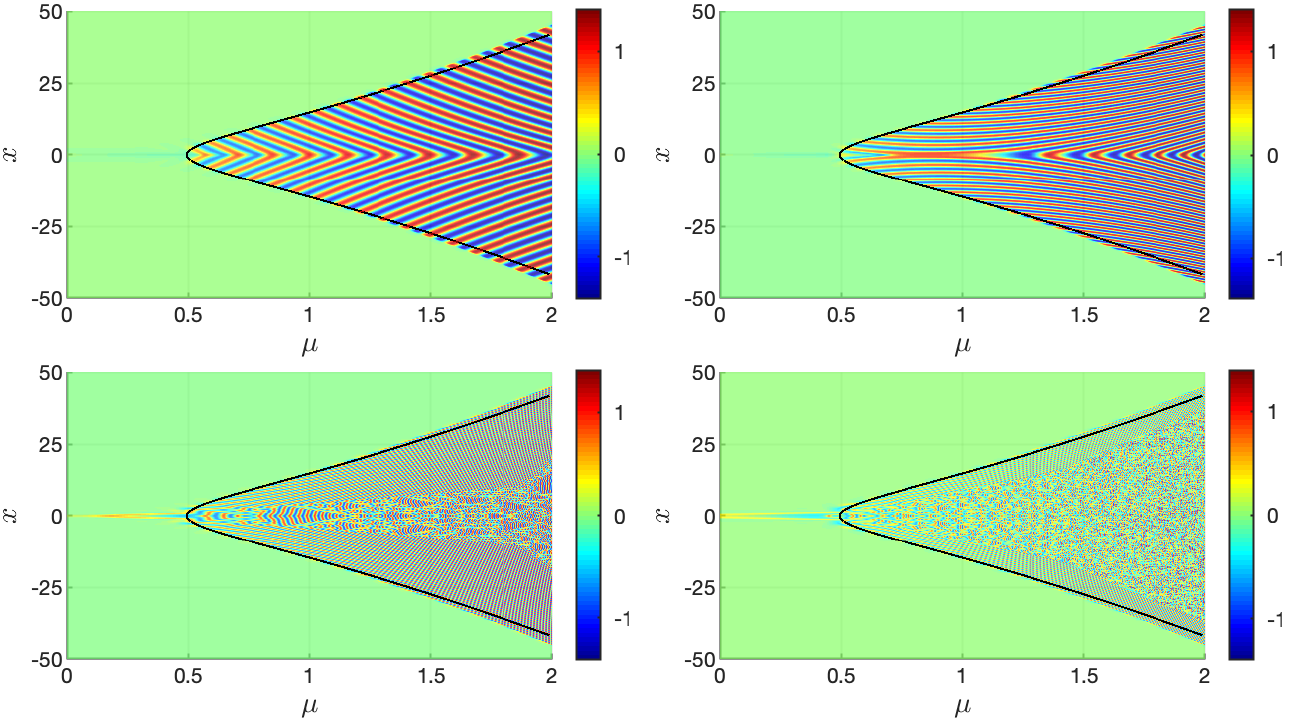}
   \caption{
Brief illustration
of how, beyond the space-time buffer curve,
$\alpha$ affects the direction of propagation
of the maxima and minima of the oscillations
in the base case with a Gaussian source.
(a) $\alpha=0.1$,
(b) $\alpha=0.6$,
(c) $\alpha=2.5$,
(d) $\alpha=5.0$.
In all frames,
the space-time buffer curve 
(black curve)
given by \eqref{eq:stbc-Gaussian}
is superimposed.
The parameters are
$\eps=0.01, \omega_0=\frac{1}{2}, d_R = 1, d_I = 0$,
and $\sigma=\frac{1}{4}$.
The initial data at $\mu_0= -1$ 
is $A_0(x) = -\sqrt{\eps} \frac{I_G(x)}{\mu_0+i \omega_0}$. 
}
   \label{fig:Gaussian-fourvaluesofalpha}
\end{figure}

Post DHB, several types of spatio-temporal patterns are observed
depending on the type of inhomogeneity, the specific case of DHB, and
the dispersion parameters $\alpha$ and $d_I$.  For instance, 
in Case 1 of DHB, Gaussian source terms break the $x$-translation symmetry
in the system, and we find that
it leads to the bifurcation of large-amplitude
periodic wave-trains which organize into stationary, or ``pinned," defects,
see Figures
~\ref{fig:Gaussian-stbc} and
~\ref{fig:Gaussian-fourvaluesofalpha}.

In the nonlinear CGL \eqref{eq:gen-CGL},
with the parameter $\mu>0$ held constant ({\it i.e.,} $\mu_t=0$), 
periodic waves have the explicit form $r e^{i(kx
- \omega t)}$  with amplitude and nonlinear dispersion relations
\begin{equation}\label{eq:nl-disp}
r^2 = \mu - \eps d_R k^2,\qquad\qquad\omega(k) = -\omega_0 +
\alpha \mu - \eps(\alpha d_R - d_I) k^2.
\end{equation}
The amplitude relation shows that periodic patterns 
exist for wavenumbers
$|k|<\sqrt{\mu/\eps d_R}$.
The dispersion relation
shows that the frequency changes sign 
at $\mu = \frac{\omega_0}{\alpha}$, provided $\alpha\ne 0$,
to leading 
order (with $0<\eps \ll 1$ and $k=\mathcal{O}(1)$).
From these relations,
the phase and group velocities for patterns with
wavenumber $k$ take the form
$$
c_p := \omega(k)/k,\qquad\qquad c_g := \omega'(k)  = 2 \eps (
d_I-\alpha d_R) k.
$$
Then, for a dynamic Hopf parameter, we expect amplitudes to vary
adiabatically as $\mu$ is increased, unless a stability boundary is
reached, after which we expect a secondary dynamic bifurcation. 

In  Figure~\ref{fig:Gaussian-stbc} (where $\alpha = 0$, $d_R = 3,$ and 
$d_I =1$), we observe that the heterogeneity induces a symmetric defect
which connects wavetrains of wavenumber $k_+ = k$ and $k_- = -k$ on
$x>0$ and $x<0$, respectively, with $k\approx 1.15$. Inward pointing
phase curves indicate the arrangement $c_{p,-} > 0 > c_{p,+}$ which,
since $\omega(k) <0$ for all $\mu$ and $\eps \ll 1$, indicates a
wavenumber arrangement $k_-<0<k_+$. This then implies that the group
velocities take the form $c_{g,-} < 0 < c_{g,+}$, so that the defect
is a source. 

To leading order, a fixed wavenumber is selected at the space-time
buffer curve which persists for increasing $\mu\in [0.5,2]$ with
fixed oscillation frequency $\omega(k)$ (note that, with $\alpha=0$, 
the frequency $\omega(k)$
has no direct dependence on $\mu$) and increasing amplitude $r(k)$ as
$\mu$ increases. As the space time buffer curve expands in $x$,
additional periods of the wavetrains are added, leading to a
persistent defect solution.  We also remark that in Case 1, where
the DHB is governed by the growth of the particular solution, the
large amplitude pattern is independent of small white noise
perturbations of the initial data. This indicates that such dynamic
Hopf heterogeneities could be used in applications to \emph{select}
organised patterns.

In  Figure~\ref{fig:Gaussian-fourvaluesofalpha} 
(where $\alpha>0$, $d_R = 1, d_I=0$), 
the defect dynamics are somewhat different. 
Both for $\alpha=0.1$ and $\alpha=0.6$
(see the top frames),
the defect at $x=0$ is a symmetric sink.
This may be seen as follows.
For $\alpha = 0.1$, we find that
$\omega(k)<0$ for
$\mu\in[\mu_{\mathrm{stbc}}(x),\frac{\omega_0}{\alpha}]$
(which lies beyond the edge of the figure),
and the phase velocities 
point inward. 
Hence, one readily calculates that, in this case, $c_{g,-}
> 0 >c_{g,+}$, so that the defect is a symmetric sink. 
Then, for $\alpha = 0.6$,
a sink defect is also
formed as $\mu$ increases past $\mu_{\mathrm{stbc}}(x)$ for
$x\sim 0$.  
As $\mu$ increases further, the frequency $\omega(k)$
changes sign at $\mu\approx \frac{\omega_0}{\alpha} =\frac{5}{6}$ 
to leading order in
$\eps \ll 1$ so that phase velocities switch direction,
$c_{p,-}<0<c_{p,+}$. However,
the defect is still a sink with inward
pointing group velocities $c_{g,-} > 0 > c_{g,+}$.

We remark that symmetric sink defects are in general transverse
heteroclinic orbits when viewed in the appropriate spatial dynamics
formulation, and come in one-parameter families parameterised by the
wavenumber $k$ \cite{SS2004,vS1992}.  Hence, we expect them to persist
under the slow drift of $\mu$. We also mention that source defects,
which select the asymptotic wavenumbers $(k_-,k_+)$, are not
transverse heteroclinics and only exist for isolated asymptotic
wavenumber pairs. These solutions are in general robust under the
introduction of a small localised heterogeneity, and it is believed
that they become pinned at the site of the heterogeneity, as we
numerically observe here.  It thus would be interesting to more fully
understand how source and sink defects are selected in the presence
of a spatio-temporal dynamic Hopf bifurcation.

For $\alpha$ yet larger, say $\alpha = 2.5$ and $\alpha=5$, 
we find that core and
wave-train instabilities begin to arise. Since $d_I = 0$, the
Benjamin-Feir condition $1+\alpha d_I>0$ holds and the band of stable
wavenumbers is given, for stationary $\mu$, by 
$$
 |k|^2 < \frac{\mu (1 +\alpha d_I)}{3+\alpha d_I + 2\alpha^2} =
\frac{\mu }{3+ 2\alpha^2}.
$$ 
Hence, an increase in $\alpha$ brings a decrease in the range of
stable wavenumbers. See \cite{AK2002,vS1992} for more information 
on this topic. The effects of adiabatic changes in $\mu$
on these large-amplitude spatio-temporal patterns
is a topic of future research.

%----------------------------------------------------------------------------
\subsection{Open questions}		\label{sec:open}
%----------------------------------------------------------------------------

The analysis and numerical simulations 
presented in this article
raise a number of open questions
and suggest avenues of additional research,
beyond that on the effect of adiabatic
variation in $\mu$ on the post-DHB spatio-temporal patterns,
just discussed in Section~\ref{sec:Post-DHB}.

The formal calculations 
in Sections~\ref{sec:lin-CGL} and \ref{sec:gen-CGL}
of the attracting and repelling QSS 
of the linear and nonlinear systems
for $\mu < 0$ and $\mu > 0$, 
respectively,
show that they are unique to all orders,
but differ in the exponentially small terms.
This suggests studying the Gevrey regularity properties
of the QSS.
Moreover, this question is motivated directly
by the recent demonstration that slow manifolds
for analytic fast-slow systems of ODEs
are Gevrey regular, in the absence of singularities 
in the slow flow, see \cite{dMK2020},
which is a major extension of the classical results
in \cite{F1979}.

Next, we expect that the nonlinear analysis of Section
\ref{sec:gen-CGL}, which is only formal in this article, 
could be made rigorous
and used to prove the convergence of the sequence of functions
$\{b_j\}$ in an appropriate function space. In particular, one would
try to show that 
$$
e^{\frac{1}{2\eps}(\mu+i\omega)^2}(b_{j+1} - b_j) \sim 
\eps^{\frac{2j+3}{2}},
$$
uniformly in $x$ to conclude that $\{b_j\}_{j\in\mb{N}}$ is Cauchy
and hence convergent to a fixed point $b_*$ of the mapping $H$ as
$j\rightarrow+\infty$. This limiting fixed point would, by
construction, yield a solution of the mild formulation
\eqref{e:mild0} with the desired dynamics for $\mu\in
[\mu_0,\omega_0-\tilde\delta)$. Standard bootstrapping and
energy estimate methods \cite{LO1996} could be used to lift this mild
solution to a classical solution.

The results in the base case of \eqref{eq:gen-CGL}
with unbounded
and sign-changing source terms
(Section~\ref{sec:ag+sc})
raise the question about how to generalize
the analysis of DHB beyond the assumptions made here.
Good agreement is found 
between the analytically calculated DHB
and numerically calculated exit times from the neighbourhood
of the repelling QSS and the onset of the oscillations,
even though the general conditions
in Section~\ref{sec:lin-CGL} are not satisfied
with these source terms.

In the study of DHB for solutions of \eqref{eq:gen-CGL},
the analyses of $A_h$, the homogeneous component,
in Section~\ref{sec:Ah-transition}
and of the solution $A$ of the full cubic PDE \eqref{eq:gen-CGL}
in Section~\ref{sec:gen-CGL}
were done working with $\mu$ on the real line.
Analyticity was used in the study 
of the particular solution, $A_p$, of the linear PDE,
to deform the contour in the complex $\mu$ plane
to derive the space-time buffer curve.
It is of interest to determine
whether analyticity is strictly necessary for DHB,
or instead if some slightly weaker assumption might suffice
to observe DHB. 
We refer to sections 2 and 3 of \cite{BB2015}, for example,
for aspects of the general theory of Gevrey regularity
of dissipative PDEs with analytic nonlinearities.
Also, there is a well-developed theory
of Gevrey regularity for solutions
of nonlinear R-D equations on bounded domains,
see for example \cite{P1991},
especially the results for the GL equation in section 2.2 there.
This could be useful for developing
a better understanding 
of DHB in nonlinear
R-D equations, the minimal hypotheses necessary for it,
and the Gevrey regularity of the solutions.

There are a number of possible avenues to explore 
to rigorously determine whether or not 
slow invariant manifolds exist
on $\mu\in[\mu,-\delta)$ and $\mu\in
(-\delta,-\mu_0]$, respectively.
We are presently studying this question for \eqref{eq:gen-CGL}
using scaling variables on the real line and
previous results on slow manifolds; see for example \cite{HK2020} and
references therein.

Finally, the type of competition analyzed herein
between the homogeneous and particular solutions
of the linearised problem,
in terms of which stops being exponentially small first
and transitions through modulus one to becoming exponentially large,
arises in other PDEs.
We are currently exploring it in a number of other
reaction-diffusion systems that exhibit DHB,
see \cite{KV2018}.

%---------------------------------------------------------------------------------
\section*{Acknowledgements}		
%---------------------------------------------------------------------------------
We thank 
Antonella Cucchetti, 
Arjen Doelman, 
Edgar Knobloch,
and Gene Wayne,
for useful comments and questions.
T.K. and T.V. thank the Lorentz Center (Leiden University, NL)
for their hospitality 
and stimulating research environment
during the 2017 Workshop on
Singular Perturbation Theory.
T.K. thanks the organisers of the February 2019
conference on Advances in Pattern Formation (Ben-Gurion University, Israel),
where the initial results about DHB in the CGL equation were presented.
T.K. is grateful to Prof. G.B. Whitham, 
from whom he learned 
the methods of stationary phase and steepest descents.
This work was supported by the US National Science Foundation
[DMS-2006887 to R.G., DMS-1616064 to T.J.K., and DMS-1853342 to T.V.].

%---------------------------------------------------------------------------------------------------------------------------------------------------
%---------------------------------------------------------------------------------

%---------------------------------------------------------------------------------------------------------------------
\appendix
%---------------------------------------------------------------------------------------------------------------------

%---------------------------------------------------------------------------------
\section{Second parametrisation for the integral
$\mathcal{I}_{a2}$ in Section~\ref{sec:lin-CGL-2}}
\label{sec:App-A}
%---------------------------------------------------------------------------------

\noindent
In this appendix, 
we briefly show that $\mathcal{I}_{a2}$
\eqref{eq:Ia2-final}
may be derived also using a second,
implicit parametrisation (by $\sigma$)
of the contour
$C_{a2}$,
\begin{equation}\label{eq:Ca2implicit}
C_{a2}:  \quad
-\frac{1}{2} ({\tilde{\mu}} + i \omega_0)^2 
= -\frac{1}{2} (\mu + i \omega_0)^2
+ \sigma, \quad \sigma \le 0.
\end{equation}
Here, $\sigma$ 
increases monotonically along $C_{a2}$,
starting from $-\infty$,
hitting $-\frac{1}{2}(\omega_0^2 - \mu^2)$ at $q_a$,
and continuing up to zero at the point $\mu$.
Therefore, the integral for $\mathcal{I}_{a2}$
may be written as
\begin{equation}
\mathcal{I}_{a2}
= \frac{e^{-\frac{1}{2\eps}(\mu+i\omega_0)^2}}{\sqrt{\eps}}
\int_{-\infty}^0
g(x,\mu - {\tilde{\mu}}(\sigma) )
e^{\frac{\sigma}{\eps}}
\frac{d{\tilde{\mu}}}{d\sigma}
d\sigma.
\end{equation}

Now,
we fix $\Delta > 0$ sufficiently small,
independent of $\eps$.
The dominant
contribution to the integral
comes from
$\sigma \in (-\Delta,0]$,
since the integrand in $\mathcal{I}_{a2}$
is exponentially small
for all $\sigma < -\Delta$,
due to the exponential term,
where we recall 
$\mu \in [-\omega_0,-\delta]$ is fixed,
$g$ is analytic and hence bounded,
and the differential is bounded 
along $C_{a2}$.
Taylor expanding
$g$ and the differential
about $\sigma=0$,
and using
${\tilde{\mu}}(\sigma) 
= -i\omega_0 
     + (\mu + i \omega_0) 
     \left[ 1 - \frac{2\sigma}{(\mu + i\omega_0)^2}\right]^{\frac{1}{2}},$
we find
\begin{equation}
g(x,\mu-{\tilde{\mu}}(\sigma)) 
= g\left(x,\mu + i\omega_0 
          -(\mu+i\omega_0)
           \left[1-\frac{2\sigma}{(\mu+i\omega_0)^2}\right]^{\frac{1}{2}}\right)
= g\left(x, \frac{\sigma}{\mu+i\omega_0} 
   +\mathcal{O}\left(\frac{\sigma^2}{(\mu+i\omega_0)^3}\right) \right),
\end{equation}
\begin{equation}
\frac{d{\tilde{\mu}}}{d\sigma}
=\frac{-1}{(\mu+i\omega_0)} \left[ 1 - \frac{2\sigma}{(\mu+i\omega_0)^2} \right]^{-\frac{1}{2}}
= \frac{-1}{(\mu+i\omega_0)}
  \left[1 + \frac{\sigma}{(\mu+i\omega_0)^2} 
+ \mathcal{O} \left( \frac{\sigma^2}{(\mu+i\omega_0)^4} \right) \right].
\end{equation}
Hence, 
\begin{equation}
\mathcal{I}_{a2}
= \frac{
-e^{-\frac{1}{2\eps}(\mu + i \omega_0)^2}}
{\sqrt{\eps} (\mu+i\omega_0)} 
\int_{-\Delta}^0
 g\left(x, \frac{\sigma}{\mu+i\omega_0} 
   +\mathcal{O}\left(\frac{\sigma^2}{(\mu+i\omega_0)^3}\right) \right)
 e^{\frac{\sigma}{\eps}}
 \left[ 1 + \mathcal{O}\left( \frac{\sigma}{(\mu+i\omega_0)^2}  \right) \right]
d\sigma.
\end{equation}
Then, 
by \eqref{eq:Bp-def}(b) and standard properties of the 1-D heat kernel,
we observe that
$\lim_{ \chi \to 0^+} g(x,\chi) = I_a(x).$
Therefore, 
we directly obtain
\eqref{eq:Ia2-final},
and hence the same formula for $A_p$.

\bigskip
\noindent
{\bf Remark.}
The two parametrisations 
\eqref{eq:Ca2-muR}
and \eqref{eq:Ca2implicit}
are consistent.
For example, 
expressing \eqref{eq:Ca2-muR}
in terms of small $\vert \sigma \vert$,
we find
${\tilde{\mu}}_R(\sigma)
= \mu - \frac{\mu\sigma}{\mu^2 + \omega_0^2} + \mathcal{O}(\sigma^2)$,
and hence 
${\tilde{\mu}}(\sigma)
= \mu - \frac{\mu-i\omega_0}{\mu^2 + \omega_0^2} \sigma 
+ \mathcal{O}(\sigma^2)$,
via \eqref{eq:Ca2-muR}.
This is identical to the expansion
of \eqref{eq:Ca2implicit},
${\tilde{\mu}}(\sigma) = \mu - \frac{\sigma}{(\mu+i\omega_0)}
+ \mathcal{O}\left( \frac{\sigma^2}{(\mu+i\omega_0)^3} \right).$

%---------------------------------------------------------
\section{The solutions of \eqref{eq:lin-CGL} are also near the 
QSS for $\mu \in (-\delta,\delta)$}
\label{sec:App-B}
%---------------------------------------------------------

In this appendix,
we briefly show how the asymptotic analysis
of the solutions $A_p(x,\mu)$
of the linear CGL \eqref{eq:lin-CGL} in the base case,
which was carried out separately
for $\mu \in [-\omega_0,-\delta]$ 
in Section~\ref{sec:lin-CGL-2}
and for $\mu \in [\delta,\omega_0]$,
in Section~\ref{sec:lin-CGL-3}
extends to 
$\mu \in (0,\delta)$, 
and to $\mu \in (-\delta,0]$,
where we recall that $\delta>0$ 
is $\mathcal{O}(1)$ and small.

First, we consider the case of
$\mu \in (0,\delta)$.
The contour is again 
$C_r= C_{r1} \bigcup C_{r2} \bigcup C_{r3} \bigcup C_{r4}$,
as in the previous subsection.
Here, the point $q_r$ is close to the saddle
at $-i\omega_0$ since $\mu \in (0,\delta)$,
recall Figure~\ref{fig:Contour-Cr}.
The contribution from $C_{r1}$ and $C_{r2}$
is again 
\begin{equation}
A_p(x,\mu) 
= \sqrt{\frac{\pi}{2}} \left( g(x,\mu+i\omega_0)
+\mathcal{O}(\sqrt{\eps}) \right)
e^{\frac{1}{2\eps}(\mu+i\omega_0)^2},
\end{equation}
where the dominant contribution
comes from the final portion of the rise along $C_{r2}$,
the steepest ascents path up to the saddle.
Then, the contributions
along $C_{r3}$ and $C_{r4}$ yield
\begin{equation}
\label{eq:Ap-mu-nearzero}
\begin{split}
A_p(x,\mu) 
&= - \sqrt{\eps} \frac{I_a(x)}{\mu + i \omega_0}
+ \eps^{\frac{3}{2}} \left( \frac{I_a(x)+d(\mu+i\omega_0){I_a}''(x)}
              {(\mu+i \omega_0)^3}
        \right)
+ \mathcal{O}\left(\frac{\eps^{\frac{5}{2}}}{(\mu+i\omega_0)^5}\right) \\
&+ \left( \sqrt{\frac{\pi}{2}} + c(\mu) \right)
\left( g(x,\mu+i\omega_0)
+\mathcal{O}(\sqrt{\eps}) \right)
e^{\frac{1}{2\eps}(\mu+i\omega_0)^2}.
\end{split}
\end{equation}
Here,
$c(\mu) = (1+i) \int_0^{\sqrt{\frac{\omega_0 \mu}{\eps}}}
e^{-i {\tilde \sigma}^2} d {\tilde \sigma}$,
which comes from the contribution along $C_{r3}$,
the stationary phase path from emerging from the saddle; 
recall \eqref{eq:Ir3-asymptotics}.
The function $c(\mu)$
increases monotonically
from zero in the limit $\mu \to 0^+$
to $\sqrt{\frac{\pi}{2}}$
for $\delta > 0$ and for other 
${\cal O}(1)$ values of $\mu>0$.

Next, we consider the case 
of $\mu \in (-\delta,0)$.
We again use the contour $C_a$ 
from Section~\ref{sec:lin-CGL-2},
recall Figure~\ref{fig:Contour-Ca}.
Then, $A_p(x,\mu)$
is given by \eqref{eq:Ap-QSS-a},
since also for $\mu \in (-\delta,0)$
the dominant part comes from the final segment 
up to $\mu$ on the real axis.

Finally, we consider the case of $\mu=0$.
The contour of integration
is the union of the line segments
\begin{equation}
\mathcal{L}_1 = \{ 
{\tilde{\mu}}_I = - \omega_0, 
{\tilde{\mu}}_R  \le 0 \} 
\quad
{\rm and} 
\quad
\mathcal{L}_2 = \{ 
{\tilde{\mu}}_R = 0, 
-\omega_0 \le {\tilde{\mu}}_I  \le 0 \}.
\end{equation}
These line segments 
lie on the horizontal and  vertical
anti-Stokes lines 
$\psi=0$, respectively, through the saddle.
Application of the method of steepest descents
along $\mathcal{L}_1$ yields
$\left( \sqrt{\frac{\pi}{2}} g(x,i\omega_0)  + \mathcal{O}(\sqrt{\eps})\right)
e^{-\frac{\omega_0^2}{2\eps}}$.

Then, the contribution from the steepest ascents
path $\mathcal{L}_2$ consists of two pieces,
one from the lower limit of integration
and a second piece from the upper limit of integration
at ${\tilde{\mu}}=0$.
The former cancels out 
the contribution from $\mathcal{L}_1$.
Hence, for $\mu=0$, we find
\begin{equation}
\label{eq:Ap-mu=0}
A_p(x,0)
= \sqrt{\eps} \frac{i I_a(x)}{\omega_0} 
+ \eps^{\frac{3}{2}} \left( 
\frac{i I_a(x) - d \omega_0 I_a''(x)}{\omega_0^3} 
\right) 
+\mathcal{O}(\eps^{\frac{5}{2}}).
\end{equation}
Therefore, the solution is also the QSS to all orders
at $\mu=0$,
just as 
for $\mu \in (-\delta,0)$
and for $\mu \in (0,\delta)$.

%=========================================================================================
\end{document}